\documentclass[aop,preprint]{imsart}

\RequirePackage[OT1]{fontenc}
\RequirePackage{amsthm,amsmath}
\RequirePackage[numbers]{natbib}
\RequirePackage[colorlinks,citecolor=blue,urlcolor=blue]{hyperref}

\arxiv{arXiv:1205.4749}
\usepackage{pkgfile}

\newcommand{\odeg}{\mu \langle \Delta_o \rangle}
\newcommand{\odege}{\mu^e \langle \Delta_o \rangle}
\newcommand{\cSG}{\cG}
\newcommand{\ptheta}{p}
\newcommand{\AM}{M}
\renewcommand{\blue}{}           
\renewcommand{\bm}{\ul}
\newcommand{\Jsp}{{\rm J}}       
\newcommand{\Amg}{{\rm M}}       
\newcommand{\Fsp}{\bar{\Jsp}}    
\newcommand{\rmg}{m}             
\newcommand{\ola}{\bar{a}}       
\newcommand{\HF}{F}              
\newcommand{\gth}{\gamma}        
\newcommand{\hf}{f}              
\newcommand{\cUs}{\cU_*}         
\newcommand{\olS}{W}             
\newcommand{\Pb}{{\sf P}}        
\newcommand{\oPb}{\ol{\sf P}}    
\newcommand{\opsi}{\ol{\varphi}} 
\newcommand{\tpsi}{\wt{\varphi}} 
\newcommand{\tgm}{\wt{\gm}}      
\newcommand{\sfp}{p_\star}       

\begin{document}

\begin{frontmatter}
\title{Ferromagnetic Ising Measures \\ on Large Locally Tree-Like Graphs}
\runtitle{Ferromagnetic Ising Measures on Large Locally Tree-Like Graphs}

\begin{aug}
\author{\fnms{Anirban} \snm{Basak}\thanksref{t1,m1}\ead[label=e1]{anirbanb@math.duke.edu}}
\and
\author{\fnms{Amir} \snm{Dembo}\thanksref{t1,t2,m2}\ead[label=e2]{adembo@stanford.edu}}

\thankstext{t1}{Part of the work was done when the authors participated in the MSRI program on Random Spatial Processes.}
\thankstext{t2}{Research partially supported by NSF grant DMS-1106627}
\runauthor{A.~Basak and A.~Dembo}

\affiliation{Duke University\thanksmark{m1} and Stanford University\thanksmark{m2}}

\address{Department of Mathematics\\
Duke University\\
\printead{e1}}

\address{Department of Mathematics and Department of Statistics\\
Stanford University\\
\printead{e2}}
\end{aug}

\begin{abstract}
We consider the ferromagnetic Ising model on a sequence of graphs 
$\Graph_n$ converging locally weakly to a rooted random tree. Generalizing \cite{mossel_mont_sly}, under an appropriate ``continuity'' property, we show that the Ising measures on these graphs converge locally weakly to a measure, which is obtained by first picking a random tree, and then the symmetric mixture of Ising measures with $+$ and $-$ boundary conditions on that tree. Under the extra assumptions that 
$\Graph_n$ are edge-expanders, 
we show that the local weak limit of the Ising measures conditioned on positive magnetization, is the Ising measure with $+$ boundary condition on the limiting tree. 
The ``continuity'' property holds except possibly for 
countable many choices of $\be$, which for 
limiting trees of minimum degree at least three,
are all within certain explicitly specified compact interval. 
We further show the edge-expander property for 
{(most of)}
the configuration model graphs corresponding to limiting (multi-type) Galton Watson trees.
\end{abstract}



\end{frontmatter}

\section{Introduction}
\noindent
The \emph{ferromagnetic Ising model} on a finite undirected graph {$\blue{\Graph}=(V,E)$, 
 is the probability distribution over $\underline{x}= \{x_i  :  i \in V \}$ with $x_i \in \{-1,+1\}$, for some $\be \ge 0$ ({\em inverse temperature} parameter), $B \in \R$ ({\em external magnetic field}),  given by
\begin{equation} \label{eq:Ising_model_descrip}
\nu_\Graph^{\be,B}\lp \ul{x} \rp = \frac{1}{Z_\Graph(\beta, B)} \exp \Big\{\be \sum_{( i,j) \in E} x_i x_j + B \sum_{i \in V} x_i\Big\},
\end{equation}
where $Z_\Graph(\be ,B)$} is the normalizing constant 
(also known as \emph{partition function}).

The Ising model is a paradigm model in statistical physics \cite{niss}, with much recent interest also in the Ising model on \emph{non-lattice complex networks} (see \cite{newmann}, and the references therein). In this paper we focus on sparse 
graph sequences $\{\Graph_n\}_{n \in \N}$ converging locally weakly to (random) trees (see
Definition \ref{dfn:lwc_graph}).
The study of statistical physics models on such graphs is motivated by numerous examples from combinatorics, computer science and statistical inference (c.f. \cite{dembo_montanari_survey, mezard_montanari}). The key to such studies is the asymptotics of the log partition function, appropriately scaled,  as derived for example in  \cite{dembo_mont, gerschenfield_montanari, sanctis_guerra}. In particular, \cite{dembo_montanari_sun} shows that for any sequence of graphs 
{$\Graph_n=(V_n,E_n)$, with $V_n$ of size $n$},
that converges
locally weakly to random trees, the asymptotic {\em free entropy density} of the ferromagnetic Ising {models 
\begin{equation} \label{eq:model_descrip}
\nu_n^{\be,B}\lp \ul{x} \rp = \frac{1}{Z_n( \beta, B)} \exp \Big\{\be \sum_{( i,j) \in E_n} x_i x_j + B \sum_{i \in V_n} x_i\Big\}\,,
\end{equation}
}
exists, i.e.,
{\beq
\phi(\be, B) := \lim_{n \ra \infty}  \phi_n (\be, B),
\eeq
where $\phi_n(\be,B):=\f{1}{n} \log Z_n (\be, B)$.}
Beyond that, perhaps the most interesting feature of the distribution in (\ref{eq:Ising_model_descrip}) is its ``phase transition'' phenomenon. Namely, for a wide class of graphs, the Ising measure for large enough $\be$ and $B=0$ decomposes into convex combination of well-separated simple components. 
This has been shown for the complete graph \cite{ellis_newman}, and for grids \cite{aizenman,  bodineau, dobrushin_shlosman, georgii_higuchi}. 

In the context of tree-like graphs $\Graph_n$, where the neighborhood of a typical vertex has, for large $n$, approximately the law of the neighborhood of the root of a randomly chosen limiting tree, this picture is only proven for a $k$-regular limit, see Montanari, Mossel and Sly \cite{mossel_mont_sly}. We show here the {\em universality} of this phenomenon, applicable for a general sequence of locally tree-like graphs, including in particular, 
Erd\"{o}s-R\'{e}nyi graphs, random uniform 
{$q$-partite} graphs, and random graphs of a given degree distribution.  More precisely, one expects that the marginal distribution of $\nu_n^{\be,B} (  \cdot  )$ converges to the marginal distribution on a neighborhood of the root for some Ising Gibbs measure on the limiting tree $\Tree$. Denoting by $\nu_{\pm,\Tree}^{\be,B}$ the Ising Gibbs measures on $\Tree$, corresponding to plus and minus boundary conditions, for $B >0$ it easily follows from \cite{dembo_montanari_sun} that, the limiting measure is given by first picking the random tree $\Tree$, and then conditioned on $\Tree$, using the Ising Gibbs measure $\nu_{+,\Tree}^{\be,B}$ (the same applies for $B <0$ with $\nu_{+,\Tree}^{\be,B}$ is replaced by $\nu_{-,\Tree}^{\be,B}$). Recall that for $B=0$ and $\be$ large, there are uncountably many Ising Gibbs measures, hence the convergence to a  particular Gibbs measure is not at all clear, as is the choice of the correct Gibbs measure. As demonstrated in \cite{mossel_mont_sly}, for $k$-regular trees, the plus/minus boundary conditions play a special role. Indeed, it is shown in \cite{mossel_mont_sly} that if $\Graph_n$'s converge locally weakly to $k$-regular trees $\Tree=\Tree_k$ then, for any $\be > 0$ and $B=0$,
\beq \label{eq:mms_result}
\nu_n^{\be,0} ( \cdot  ) \ra \f{1}{2} \nu_{+,\Tree}^{\be,0}(  \cdot  ) + \f{1}{2} \nu_{-,\Tree}^{\be,0}( \cdot  ).
\eeq 
It is further shown there that, when the graphs 
$\{\Graph_n\}_{n \in \N}$ are {\em edge-expanders} ,
\beq \label{eq:mms_result2}
\nu_{n, \pm}^{\be,0} (  \cdot  ) \ra \nu_{\pm, \Tree}^{\be,0} ( \cdot  ), 
\eeq
where $\nu_{n,+}^{\be,0}( \cdot  )$ and $\nu_{n,-}^{\be,0} ( \cdot  )$ are the measures (\ref{eq:model_descrip}) conditioned to, respectively, 
$\sum_{i} x_i {\ge} 0$ 
and $\sum_{i} x_i {\le} 0$
{(when $n$ is odd, see Remark \ref{rmk:even} on
slight modification usually taken for even $n$).} 
The latter sharp result provides a better understanding of $\nu_n ( \cdot  )$, and is much harder to prove than (\ref{eq:mms_result}). For genuinely random limiting trees, one expects (\ref{eq:mms_result}) and (\ref{eq:mms_result2}) to apply where now $\Tree$ is chosen according to the limiting tree measure. {As we focus on the case $B=0$, hereafter 
{we write $\nu_n^\be(\cdot):=\nu_n^{\be,0}(\cdot)$
and adopt the convention of using $\nu_n^B(\cdot)$ (or 
just $\nu_n$, in case $B=0$), when the value of $\be$ 
is either arbitrary, or clear from the context. 
Similar notations apply} for Ising measures on 
the limiting trees.}

{It is well known (see \cite{lyons_ising}) that there exists a value of $\be$, denoted here by $\be_c$, such that for $\be < \be_c$ there is a unique Ising Gibbs measure, and for $\be > \be_c$ there are multiple Ising Gibbs measures.}
{In the more interesting case of $\be \ge \be_c$}, key estimates in the proof of (\ref{eq:mms_result}) and (\ref{eq:mms_result2}) in \cite{mossel_mont_sly}, involve explicit calculations which crucially rely on the regularity of both graph sequence, and the limiting tree. Several new ideas are necessary in the absence of such regularity. For example, the key to the proof of (\ref{eq:mms_result}) in \cite{mossel_mont_sly} is the continuity, for $k$-regular infinite trees, of root magnetization under $\nu_{+,\Tree_k}(\cdot)$, obtained there out of its representation as the largest zero of a real analytic function. {While no such representation is known for any other possible 
limiting tree measure, in case it a.s. has minimum degree 
$d_\star > 2$, we prove here the continuity 
of root magnetization under $\nu_{+,\Tree}(\cdot)$ 
for all 
$\be  > \atanh[(d_\star-1)^{-1}] > \be_c$} 
(see Section \ref{section:U_beta_continuous}).\footnote{for $\be=\be_c$ one may use the equivalent capacity criterion provided in \cite{pemantle_peres}.} 
The proof of (\ref{eq:mms_result2}) relies on choosing functionals $\Fsp^l(\cdot)$ 
of the spin configurations on $\Graph_n$, which approximate the indicator on the vertices that are in ``$-$ state'', and whose values concentrate as $n,l \ra \infty$. The regularity of the graphs $\Graph_n$, and that of their limit, provide for such functionals, and allows explicit computations involving them, both of which fail as soon as we move away from the regular regime. At the level of generality of our setting the only tools are {\em unimodularity} of the law of the 
limiting tree (see Definition \ref{def:unimod}), and properties of 
simple random walk on it. Hence, a completely different choice of functionals is required here. With $\Fsp^l(\cdot)$ defined 
via {\em average occupation measure} of the {variable speed continuous time simple random walk 
(\abbr{vsrw})} 
on the tree, we show here that (\ref{eq:mms_result2}) holds under the same continuity property, for {\em any} edge-expander $\Graph_n$'s 
(see Theorem \ref{thm:average_+}). 
We also confirm the {root magnetization 
continuity property at $\be=\be_c$ for
{\em multitype Galton Watson}  
(\abbr{MGW}) trees
which arise as the limit of many natural locally tree-like graph ensembles, and show} that 
subject to minimal degree at least $3$, the 
corresponding configuration models  
are edge-expanders (see Section \ref{section:U_beta_continuous}). Thus, our theorem applies for most naturally appearing locally tree-like graphs.

An interesting byproduct of our results is the continuity of percolation probability for {\em random cluster model}, with $q=2$, and {\em wired boundary condition} (see \cite{grimmett} for details on \abbr{RCM}, and its connection with Ising model). Another interesting byproduct of this work is the uniqueness of the {\em splitting Gibbs measure} (for a definition see \cite[Chapter 12]{georgii}), for {large $\be$}, $B=0$ and any boundary condition strictly larger than the free boundary condition (see Lemma \ref{lem:recursion} and Remark \ref{rmk:splitting}). Many of the techniques developed here should extend to more general settings, e.g. the Potts model.

\subsection{Graph preliminaries and local weak convergence}

\noindent
In a connected undirected graph $\Graph=( V,E)$ the
\emph{distance} between two vertices $v_1$ and $v_2$ is defined to be the length of the shortest path between them. For each vertex $v \in V$, we denote by $\blue{\Ball_v( r)}$ the ball of radius $r$ around $v$, i.e. the collection of all vertices whose distance from $v$ in 
$\Graph$ is at most $r$. 
The set $\Ball_v(1) \setminus \{v\}$ of all vertices adjacent 
to $v$ is also denoted by $\partial v$, with 
$\Delta_v := |\partial v|$, 
denoting its size, namely, the degree of $v$ in $\Graph$.

\noindent
A \emph{rooted graph} $( \Graph,o)$ is a graph $\Graph$ with a specified vertex $o \in V$, called the \emph{root}, and 
a \emph{rooted network} $(\ol{\Graph},o)$
is a rooted graph $(\Graph,o)$
with vector $\ul{x}_\Graph$ of $\cX$-valued 
marks on each of its vertices {(for Ising 
models $\cX=\{-1,1\}$, more generally $\cX$  
assumed throughout to be a fixed finite set)}.
A \emph{rooted isomorphism} of rooted graphs (or networks) is a graph isomorphism which maps the root of one to that of another
(while preserving the marks in case of networks), with 
{$[\Graph,o]$ denoting the collection of all rooted graphs that are isomorphic to $(\Graph,o)$ (and 
$[\ol{\Graph},o]$ denoting the collection of 
all rooted networks isomorphic to $(\ol{\Graph},o)$}). 

\noindent
Let $\cG_*$ be the space of rooted isomorphism classes of rooted \emph{connected} locally finite graphs. Similarly, for rooted networks let $\overline{\cG}_*$ denote the space of rooted isomorphism classes of rooted connected locally finite networks.
Setting the distance between $[\Graph_1,o_1]$ and $[\Graph_2,o_2]$ (and the same between $[\ol{\Graph}_1,o_1]$ and $[\ol{\Graph}_2,o_2]$) to be $1/(\alpha+1)$, where $\alpha$ is the supremum over $r\in \N$ such that there is a rooted isomorphism of balls of radius $r$ around the roots of $\Graph_i$ (and marks in those balls are same),
results with $\cG_*$ and $\ol{\cG}_*$ 
which are complete separable metric spaces (see \cite{aldous_steele, benjamini_schramm}). We use hereafter this metric topology, 
denoting by $\sC_{\cG_*}$ and $\sC_{\ol{\cG}_*}$
the corresponding Borel $\sigma$-algebras on $\cG_*$ and $\ol{\cG}_*$, 
respectively \blue{(but forgo the conversion $r \mapsto 1/(r+1)$,
letting $\Ball_{\Graph}(r)$ stand throughout for the $\cG_*$-metric 
ball of radius $1/(r+1)$ around $\Graph$, namely those rooted graphs 
$(\Graph',o')$ having $\Ball_{o'}(r)$ isomorphic to 
$\Ball_o(r) \subset \Graph$).} 
Similarly, we equip the spaces $\cT_*$ and 
$\ol{\cT}_*$ of all rooted isomorphism classes 
of locally finite trees (and marked trees, respectively), 
with the metric topology and Borel $\sigma$-algebra 
induced by $\cG$ and $\ol{\cG}_*$, respectively
\blue{(while using as before $\Ball_{\ol{\Graph}}(r)$,
$\Ball_{\Tree}(r)$ and $\Ball_{\ol{\Tree}}(r)$ for the 
metric balls of radius $1/(r+1)$ in $\ol{\cG}_*$,
$\cT_*$ and $\ol{\cT}_*$, respectively).}

\begin{dfn}
For $\ze_n$ and $\mu$ Borel probability measures on $\cG_*$ (or $\ol{\cG}_*$), we write $\ze_n \Ra \mu$ when $\ze_n$ converges weakly to $\mu$ with respect to the metric on $\cG_*$ (or $\ol{\cG}_*$) and for 
any $\Graph \in \cG_*$ we denote
by $\delta_\Graph$ the probability measure on $\cG_*$ 
assigning point mass at $\Graph$.
\end{dfn}

\noindent
For probability measure $\nu$ 
on $(\cX_1,\cB_1)$ and measurable map 
$f : (\cX_1,\cB_1) \mapsto (\cX_2,\cB_2)$ 
we let $\nu \circ f^{-1}$ denote
the probability measure on $(\cX_2,\cB_2)$ such that 
$\nu \circ f^{-1} (\cdot) = \nu (f^{-1}(\cdot))$, 
and in case $f$ is real-valued, use  
the shorthand $\nu[f]$ or $\nu \langle  f \rangle$ 
for the $\nu$-expected 
value of $f$ (i.e.~$\int f d\nu$), using also 
$\langle f \rangle$ when the choice of $\nu$ is clear 
form the context. Equipped with these notations
we proceed to define the
\emph{local weak convergence} of graphs.
{
\begin{dfn} \label{dfn:lwc_graph}
For a sequence of graphs $\{\Graph_n\}_{n \in \N}$
having vertex sets $[n]$, let
$\mu_n$ denote the law of $({\Graph_n},I_n)$ in $\cG_*$ 
for $I_n$ chosen uniformly over $[n]:= \{1,2,\ldots,n\}$.
We call such $\{\Graph_n\}$ {\em uniformly sparse}, 
if $\Delta_o$ is uniformly integrable under $\{\mu_n\}$. 
That is, if
\beq\label{eq:unif-sparse}
\lim_{k \ra \infty} \limsup_{n \ra \infty} \frac{1}{n} \sum_{i \in [n]} \Delta_i(\Graph_n) \bI (\Delta_i(\Graph_n) \ge k) =0.
\eeq
If in addition $\mu_n \Ra \mu$, a probability measure on $\cG_*$,
we say that the uniformly sparse collection 
$\{\Graph_n\}$ converges locally weakly to $\mu$,  
denoted by $\Graph_n \stackrel{\text{\abbr{LWC}}}{\Lra} \mu$. {In particular, due to uniform sparseness
$\odeg
$ is finite 
for any such limit.}
\end{dfn}
}

\noindent
Similarly to the space $\cG_*$, one defines 
$\cG_{**}$ as the space of all isomorphism classes of locally finite connected graphs with an ordered pair of distinguished vertices and the corresponding topology thereon, where a function $f$ on $\cG_{**}$ is written
as $f(\Graph,x,y)$, to indicate the distinguished pair of 
vertices $(x,y)$.
{In \cite{benjamini_schramm} it is shown that any \abbr{LWC} limit point must be {\em involution invariant},
a property that was found in \cite{aldous_lyons} to 
be equivalent to the following property of {\em unimodularity}}.
\begin{dfn}\label{def:unimod}
A Borel probability measure $\mu$ on $\cG_*$ is called \em{unimodular} 
if for any Borel function $f : \cG_{**} \ra [0, \infty]$,
\beq 
\int \sum_{x \in V( \Graph)} f( \Graph,o,x) d \mu( [\Graph,o]) = \int \sum_{x \in V ( \Graph)} f( \Graph,x,o) d \mu ( [\Graph,o]) \,. \label{eq:mass_transport}
\eeq
We denote by $\cU$ the collection of all 
unimodular probability measures $\mu$ on $\cG_*$ 
for which $\odeg$ is finite \blue{and by $\cUs$ those $\mu \in \cU$ 
having $\mu(\cT_*)=1$.}
\end{dfn}

\noindent
{We consider throughout \emph{tree-like graphs},
namely $\Graph_n \stackrel{\text{\abbr{LWC}}}{\Lra} \mu$
with a limiting object which is a (random) tree, namely
having $\mu \in \cUs$.
This assumption, and the fact that any 
\abbr{LWC} limit points is in $\cU$ are both
key for our results,}  
with (\ref{eq:mass_transport}) being utilized 
in several proofs. 

\subsection{Local weak convergence of Ising measures}  
The space of all probability measures on $(\ol{\cG}_*, \sC_{\ol{\cG}_*})$ will be denoted by $\cP(\ol{\cG}_*)$. 
For example, upon choosing a root, Ising measures on 
connected, locally finite graphs can be considered 
elements of $\cP(\ol{\cG}_*)$. {Considering the 
$\cG_*$-projection $[\ol{\Graph},o] \mapsto [\Graph,o]$ from rooted
networks in $\ol{\cG}_*$ to rooted graphs in $\cG_*$,
we let $\mu \otimes \nu_\Graph$ denote an element of $\cP(\ol{\cG}_*)$, whose marginal distribution on $\cG_*$ is $\mu \in \cP(\cG_*)$, and given any fixed 
$\Graph \in \cG_*$ has the (conditional) distribution 
$\nu_\Graph$ on the corresponding mark space $\cX^{\Graph}$.} 

\noindent
{For any positive integer $t$, 
the subgraph $(\Graph,o)(t)$ of $(\Graph,o)$ induced by the 
vertices $\Ball_o(t)$, is called the graph truncated at {\em height} $t$, with the corresponding definition for a rooted network. We further use the notations $\Graph(t)$ and 
$\ol{\Graph}(t)$, when the choice of root is  
clear from the context. For example, $\Tree(t)$ denotes
the first $t$ generations of a tree $\Tree$ (i.e. 
the subtree induced by the vertices of $\Tree$ 
of distance at most $t$ from its root)}. Accordingly, 
for each $t$ 
we let $\ol{\cG}_*(t)$ denote the space of rooted isomorphism classes of rooted connected locally 
finite networks truncated at height $t$, with  $\sC_{\ol{\cG}_*(t)}$ the corresponding Borel 
$\sigma$-algebra, yielding for each 
$\ol{\nu} \in \cP(\ol{\cG}_*)$ the 
probability measure $\ol{\nu}^t$ induced on 
$(\ol{\cG}_*(t), \sC_{\ol{\cG}_*(t)})$ by such 
truncation (of the network), and for each probability 
measure ${\ol{\gm}}$ \emph{on} $\cP(\ol{\cG}_*)$ the correspondingly
induced probability measure $\ol{\gm}^t$ on
$\cP(\ol{\cG}_*(t), \sC_{\ol{\cG}_* (t)})$. 

\noindent
We next adapt \cite[Definition 2.3]{mossel_mont_sly} 
to the case of non-deterministic graph limits. 
\begin{dfn}\label{dfn:conv_weakly}
Given a sequence of graphs $\{\Graph_n\}_{n \in \N}$ having
vertex sets $[n]$, and probability measures 
$\ze_n$ on {$\cX^{V_n}$,} {for any
positive integer $t$} let $\oPb_n^t(i) \in 
\cP(\ol{\cG}_*(t), \sC_{\ol{\cG}_* (t)})$ denote 
the law of the pair $((\Ball_i(t),i),\ul{x}_{\Ball_i(t)})$ for 
$\ul{x}$ drawn according to $\ze_n$ and $i \in [n]$ 
some vertex of $\Graph_n$. 

\noindent
When combined with the 
uniform measure over the choice of 
random vertex $I_n \in [n]$, this results 
with the random distributions $\oPb_n^t(I_n)$, and we
say that {$\{(\Graph_n,\ze_n)\}_{n \in \N}$ 
(or in short $\{\ze_n\}$),} 
converges 
{\em locally weakly} to a probability measure $\ol{\gm}$ on 
$\cP(\ol{\cG}_*)$, if the law of $\oPb_n^t(I_n)$ 
converges weakly to $\ol{\gm}^t$, {as $n \to \infty$},
for each $t \in \N$.
\end{dfn}

\noindent
Notions of convergence similar to Definition \ref{dfn:conv_weakly}, and the weaker form of convergence of  Definition \ref{dfn:conv_average} were studied under the name of {\em metastates for Gibbs measures} (see \cite{aizenman_wehr, kulske, newman_stein}).

\noindent
We proceed to formally define the relevant limiting
Ising Gibbs measures $\nu^{\be,B}_{\pm,\Tree}$. 
\begin{dfn}\label{dfn:ising-pm}
For each $t$, consider the following Ising measures on $\Tree(t)$: 
\begin{eqnarray}
\nu^{\be,B,t}_{ +,\Tree}(\ul{x}) & :=  & \frac{1}{Z^{t,+}} \exp \Big\{ \beta \sum_{(i,j) \in E(\Tree(t))} x_i x_j + B \sum_{i \in V(\Tree(t))}  x_i \Big\} \bI \big(\ul{x}_{\Tree \setminus \Tree(t-1)} = (+)_{\Tree \setminus\Tree (t-1)}\big), \notag\\
\nu^{\be,B,t}_{ -,\Tree}(\ul{x}) & :=  & \frac{1}{Z^{t,-}} \exp \Big\{ \beta \sum_{(i,j) \in E(\Tree(t))} x_i x_j + B \sum_{i \in V(\Tree(t))}  x_i \Big\} \bI \big(\ul{x}_{\Tree \setminus \Tree(t-1)} = (-)_{\Tree \setminus \Tree (t-1)}\big), \notag
\end{eqnarray}
where for any $W \subseteq V(\Tree)$, we denote by $(+)_W$  
the vector $\{x_i=+1, i \in W\}$, and 
by $(-)_W$ the vector $\{x_i=-1, i \in W\}$, respectively. It is well known that as $t \ra \infty$ both $\nu_{+,\Tree}^{\be,B,t}$ and $\nu_{-, \Tree}^{\be,B,t}$ converge to probability measures on $\{-1,+1\}^\Tree$, denoted as $\nu_{+,\Tree}^{\be,B}$ (plus measure) and $\nu_{-,\Tree}^{\be,B}$ (minus measure), respectively (see \cite[Chapter IV]{liggett}). 
\end{dfn}

\noindent
For any $\be,B \ge 0$ and 
{$\mu \in \cU$ supported on the collection of 
rooted trees 
$(\Tree,o) \in  \cT_*$,} let
{
\beq
\U (\be,B) := \f{1}{2} 
\mu \Big[ 
\sum_{i \in \partial o} \nu_{+,\Tree}^{\be,B} \langle x_o  x_i \rangle 
\Big]
\,.
\eeq
} 
 
\noindent
Our first result generalizes \cite[Theorem 2.4.I]{mossel_mont_sly}, namely the limit (\ref{eq:mms_result}),
to any limiting measure $\mu$ 
supported on $\cT_*$ subject to a 
mild continuity assumption on $\U (\cdot,0)$.
\begin{thm} \label{thm:average_free}
Suppose $\Graph_n \stackrel{\text{\abbr{LWC}}}{\Lra} \mu$ {for some $\mu \in \cUs$}.
Then, {at any continuity point $\be\ge 0$ 
of the bounded, non-decreasing, right-continuous function $\U (\be, 0)$}, the Ising measures 
$\nu_n^{\be}$ on $\Graph_n$ converge locally weakly to 
$\blue{\ol{\gm}}= \mu \circ \blue{\opsi}^{-1}$, where $\opsi: \cT_* \ra \cP(\ol{\cT_*})$ with $\opsi(\Tree)= \delta_\Tree \otimes (\f{1}{2} \nu^{\be}_{+,\Tree} + \f{1}{2} \nu^{\be}_{-,\Tree})$. 
\end{thm}

\noindent
Our generalization of (\ref{eq:mms_result2}), namely 
\cite[Theorem 2.4.II]{mossel_mont_sly}, 
to all limiting tree measures, requires that 
the graph sequence has certain edge-expansion 
property related to the following definition.
\begin{dfn}
A finite graph $\Graph=(V,E)$ is a $(\delta_1,\delta_2, \lambda)$ {\em edge-expander} if, for any set of vertices $S \subseteq V$, with $\delta_1 |V| \le |S| \le  \delta_2 |V|$, we have $|\partial S| \ge \lambda |S|$, where 
$| \cdot |$ denotes the cardinality of a set and 
$\partial S$ denotes the collection of edges between $S$ and $S^c$.
\end{dfn}
 
\begin{thm} \label{thm:average_+}
Suppose $\{\Graph_n\}_{n \in \N}$ are  
$(\delta,1/2, \lambda_\delta)$ edge-expanders for 
all $0< \delta <1/2$ and some $\lambda_\delta >0$
(which is independent of $n$). If 
$\Graph_n \stackrel{\text{\abbr{LWC}}}{\Lra} \mu$ {for some $\mu \in \cUs$}, then
{at any continuity point of $\be \mapsto \U (\be, 0)$},
the measures $\{\nu^{\be}_{n, +}\}$ converge locally weakly to 
$\blue{\ol{\gm}_+}= \mu \circ \blue{\opsi_+}^{-1}$ where 
$\opsi_+ : \cT_* \ra \cP(\ol{\cT}_*)$ with $\opsi_+(\Tree)= \delta_\Tree \otimes \nu^{\be}_{+,\Tree}$. 
\end{thm}

{
\begin{rmk}\label{rmk:left-cont}
Theorems \ref{thm:average_free} 
and \ref{thm:average_+} apply except for possibly 
countable set of discontinuity points of 
$\be \mapsto \U (\be,0)$. Note  
that $\U (\be,B)$ is uniformly bounded for any $\mu \in \cU$,
and while proving Lemma \ref{lem:pair_exp_conv} we see that 
it is non-decreasing, right-continuous at all $\be,B \ge 0$,
and continuous whenever $B>0$. Further, in proving both 
theorems, left-continuity of $\U (\be,0)$ is only required  
for relating it to the limiting correlation 
$\nu_n^{\be,0} \langle x_i x_j \rangle$ across 
a uniformly chosen edge of $\Graph_n$ 
(see Lemma \ref{lem:pair_exp_conv}).
\end{rmk}
}

\begin{rmk}\label{rmk:even}
With $B=0$, for $n$ odd the probability 
measures $\nu^{\be}_{n,\pm}$ supported on 
$\pm \sum_i x_i \ge 0$ are uniquely determined by
the identity
$\nu^{\be}_n = \frac{1}{2} \nu^{\be}_{n,-} 
+ \frac{1}{2} \nu^{\be}_{n,+}$. 
To circumvent non-essential technical issues, 
one slightly modifies $\nu^{\be}_{n,\pm}$ 
in case $n$ is even to retain this property, 
as well as having 
$\nu^{\be}_{n,+}(\ul{x})=\nu^{\be}_{n,-}(\ul{x})$
whenever $\sum_i x_i =0$.
\end{rmk}

\begin{rmk}
Recall the example in  \cite[\textsection 2.3]{mossel_mont_sly}, where it is shown that 
{even in case of $k$-regular tree limits 
one can not completely dispense of the
expander-like condition when dealing with
the convergence of $\nu^{\be}_{n,+}$.}
\end{rmk}

\subsection{Configuration models and multi-type 
Galton-Watson (\abbr{MGW}) trees} 
We proceed to verify that our results apply 
for a general class of random graphs from the 
{\em configuration model}, for which the 
limiting tree follows a \abbr{MGW} distribution, 
starting with the definition of the configuration 
model we consider.

\begin{dfn}\label{dfn:config} 
{Fix a strictly positive} 
{
probability measure $\blue{\ptheta(\cdot)}$ on 
some finite (type) space $\cQ$.
Let $\Z_\ge$ denotes the set of all non-negative integers
and $\Z_{\ge}^{|\cQ|}:= \{ \ul{k}=(k_1,k_2,\ldots,k_{|\cQ|}):  k_j \in \Z_\ge, \, j =1,2,\ldots, |\cQ|\}$. 
Consider a (finite) collection of probability 
measures $P_i(\cdot)$, $i \in \cQ$
on $\Z_{\ge}^{|\cQ|}$, such that 
{for all $i,j \in \cQ$,
\begin{align}
\blue{\AM(i,j)} &:= \sum_{\ul{k}} P_i(\ul{k}) k_j  < \infty \,,
\label{eq:bd-degree}
\\
\ptheta(i) \AM(i,j) &= \ptheta(j) \AM (j,i)\,.
\label{eq:balance}
\end{align}
}

\noindent
For $m \in \N$, let an {\em $m$-star} denote the connected graph of $(m+1)$ vertices, with one vertex of degree $m$ and all others having degree one. {Such $m$-star has 
two {\em ends}, one end being its vertex of degree $m$, with the other end consisting of the remaining $m$ 
degree one vertices of the $m$-star.} 
Now for each $n$ we define 
the random graph $\Graph_n= (V_n,E_n)$ as follows.  
For every $i \in \cQ$ and $\ul{k} \in \Z_\ge^{|\cQ|}$, we create $\lfloor n \ptheta(i) P_i(\ul{k})\rfloor$ many $(\sum_j k_j+1)$-{stars} with types, such that the end of each of the stars with one vertex has type $i \in \cQ$ and the other end consists of $\sum_j k_j$ vertices, of which exactly 
$k_j$ have type $j$, for each $j \in \cQ$. 

\noindent
Edges in a star will be termed as {\em half-edges}, and we use the generic notation $(v,e_v)$ to denote a half-edge with $v$ being the single vertex at one end of the star, 
and $e_v$ being one of the vertices present in the other 
end of the star. The vertex $v$ here will be called a {\em permanent} vertex, whereas the vertices like $e_v$ will be termed as {\em floating} vertices.} 
{We denote half-edges $(v,e_v)$ having a permanent
end $v$ of type $q(v)=i$ and a floating end $e_v$ of 
type $q(e_v)=j$
by $\overrightarrow{(i,j)}$. Due to condition (\ref{eq:balance}), if not for the integer
truncation effects, for any $i,j \in \cQ$ 
the number of half-edges of type $\overrightarrow{(i,j)}$ 
would match that of type $\overrightarrow{(j,i)}$.
We thus achieve such equality between the numbers of  
$\overrightarrow{(i,j)}$ and
$\overrightarrow{(j,i)}$ half-edges, upon adding to 
$\Graph_n$ at most 
$$
2 \sum_{i,j} \sum_{\ul{k}} \{ n \ptheta(i) P_i(\ul{k}) \} k_j 
$$
half-edges. This amounts to adding only $O(1)$ 
half-edges to the stars (since $\sum_{i,j} \AM(i,j)$ is
finite, due to (\ref{eq:bd-degree})).}

{
\noindent
Thereafter for every $i,j \in \cQ$ we perform a uniform matching between half-edges with type $\overrightarrow{(i,j)}$ and half-edges with type $\overrightarrow{(j,i)}$. Once we have obtained a matching between these half-edges we throw out the floating vertices and join the permanent vertices of those half-edges, which have been matched, to get a graph with types (For example, if in a matching the half-edge $(v,e_v)$ of type $\overrightarrow{(i,j)}$ matches with the half-edge $(w,e_w)$ of type $\overrightarrow{(j,i)}$ then we join $v$ and $w$, and $q(v)=i, \, q(w)=j$). This completes 
the recipe for generating the random graph 
$\Graph_n=(V_n,E_n)$.
}
\end{dfn}

{We associate with each $\ptheta(\cdot)$ and  
collection of probability measures $P_i(\cdot)$
that satisfies the conditions of Definition \ref{dfn:config}, a unimodular version of the \abbr{MGW} law, to be denoted hereafter by \abbr{UMGW}.} 
\begin{dfn}\label{dfn:umgw}
For each $\ptheta, \, \{P_i(\cdot), i \in \cQ\}, \text{ and }\AM(\cdot, \cdot)$ satisfying  (\ref{eq:balance}) 
and
(\ref{eq:bd-degree}), let 
{$\cQ_\AM:=\{(i,j): \AM (i,j)>0\} 
\subseteq \cQ \times \cQ$ and $\wh{P}_{i,j}( \cdot )$ for  
$(i,j) \in \cQ_\AM$,} be the probability measures 
on $\Z_{\ge}^{|\cQ|}$ given by 
\beq
\blue{\wh{P}_{i,j}}(\ul{k}) = P_i(\ul{k}+e_j)\f{k_j+1}{\AM (i,j)}\,, \notag
\eeq
where $e_j$ denotes the vector with $1$ at $j^{th}$ co-ordinate and $0$ elsewhere, and we assume that $\wh{P}_{i,j}(\ul{k}) > 0$ for some $(i,j)$ and $\|\ul{k}\| :=\sum_j k_j\ne 1$ (in the branching processes literature 
this property is called {\em non-singularity},
c.f. \cite[pp. 184]{athreya_ney}).

\noindent
{We assume that 
the mean matrix $\blue{\wh{\AM}}$ for
the kernel $\wh{P}$ over $\blue{\cQ_\AM}$, which is given by} 
\beq
\wh{\AM}((i_1,j_1),(i_2,j_2)):= 
{\bI_{j_2=i_1}}
\sum_{\ul{k}} \wh{P}_{i_1,j_1}(\ul{k})
{k_{i_2}}\,, \notag
\eeq
{is {\em positive regular}. That is, 
we require 
that for some finite positive integer $r$ all 
entries of $(\wh{\AM})^r$ be strictly positive 
(possibly infinite, and when multiplying matrices we 
adopt the convention that $\infty \times 0 = 0$)}.

\noindent
The \abbr{UMGW} measure on the trees with types is the following: Type of the root is chosen according to $\ptheta(\cdot)$, and conditional on the type of the root, say $i_0$, it's off-spring number and types are chosen according to $P_{i_0}(\cdot)$. From the next generation onward, the off-spring numbers and types are chosen independently at each vertex according to $\wh{P}_{i,j}$ where $i$ is the type of the current vertex and 
$j$ being the type of its parent.
\end{dfn}

{
\begin{rmk}
In the special case $|\cQ|=1$, there are no types in
the random graphs $\Graph_n$ of Definition \ref{dfn:config}, 
neither in the random \abbr{UGW} {(\abbr{UMGW})} tree of 
Definition \ref{dfn:umgw}. The condition (\ref{eq:balance}) 
{and positive regularity then trivially hold, while 
non-singularity and (\ref{eq:bd-degree}) amount to having $P(1)<1$ and}
finite average degree $\sum_k k P(k)$. In this setting
$\Graph_n$ is the configuration model corresponding to  
uniformly chosen random graphs subject to given degree
distribution $P(\cdot)$ (c.f. \cite[Section 1.2.4]{dembo_montanari_survey}), which is uniformly sparse
and converges weakly to the corresponding \abbr{UMGW}
measure of Definition \ref{dfn:umgw} (see
\cite[Proposition 2.5]{dembo_montanari_survey}).
The latter is precisely the 
\abbr{UGW} tree measure of \cite[Example 1.1]{aldous_lyons}, 
and \cite[Section 2.1]{dembo_mont}. 

\noindent
In particular, taking $P(\cdot)$ a Poisson law of parameter $2\alpha$, results with $\wh{P}(k)=P(k)$ (i.e., here the \abbr{UGW} measure coincides with the usual \abbr{GW} law). The configuration model is then closely related to 
Erd\"{o}s-R\'{e}nyi random graph ensembles 
of $n^{-1} |E_n| \to \alpha$ which also have the 
\abbr{UGW} measure as their a.s. \abbr{LWC} limit (see 
\cite[Proposition 2.6 and Lemma 2.3]{dembo_montanari_survey}). 

\noindent
For $|\cQ|>1$ the uniform sparseness of $\{\Graph_n\}$ of 
Definition \ref{dfn:config} is an immediate 
consequence of finiteness of
$\cQ$ and $\sum_{i,j} \ptheta(i) \AM(i,j)$,
while its local weak convergence to the 
corresponding \abbr{UMGW} measure follows along the
lines of \cite[Proof of Proposition 2.5]{dembo_montanari_survey} (
from the latter convergence we know
that each \abbr{UMGW} measure 
of Definition \ref{dfn:umgw} is unimodular). 
One concrete example is the configuration model 
$\{\Graph_n\}$ and \abbr{UMGW} for random uniform 
$q$-partite, $q \ge 2$, graphs (of 
$\lfloor \alpha n \rfloor$ edges), which fit within
our framework upon taking $\ptheta$ uniform on 
$\{1,\ldots,q\}$ and $P_i(\ul{k})=\prod_{\ell \ne i} P(k_\ell)$, with $P(\cdot)$ the Poisson law of parameter
$2 \alpha q/(q-1)$. 
\end{rmk}
}

\begin{lem}\label{lem:assumption_gw_mgw}
If $\mu$ is any of the \abbr{UMGW} measures 
of Definition \ref{dfn:umgw}, with minimum degree 
$d_\star >2$, one has that 
$\be \mapsto \U(\be,0)$ is continuous except for possibly countably many values of $\be \in (\be_c,\be_\star]$, where $\be_\star= \atanh[(d_\star-1)^{-1}]$.
\end{lem} 

Thus, upon applying Theorem \ref{thm:average_free} and Theorem \ref{thm:average_+} we immediately obtain that:
\begin{cor}\label{thm:convergence_gw_mgw}
Suppose $\Graph_n \stackrel{\text{\abbr{LWC}}}{\Lra} \mu$
with $\mu$ a \abbr{UMGW} measure as in Definition 
\ref{dfn:umgw}, having a.s. minimum degree $d_\star >2$. Then, except for a possibly countably many values of $\be \in (\be_c,\be_\star]$,

\noindent
(a) $\nu_n$ converges locally weakly to $\ol{\gm}= \mu \circ \opsi^{-1}$, for $\opsi$ as in Theorem \ref{thm:average_free}.

\noindent 
(b) If in addition $\{\Graph_n\}_{n \in \N}$ are  
$(\delta,1/2, \lambda_\delta)$ edge-expanders for 
all $0< \delta <1/2$ and some $\lambda_\delta >0$
(independent of $n$), then
$\nu_{n,+}$ converges locally weakly to $\ol{\gm}_+= 
\mu \circ \opsi_{+}^{-1}$, for $\opsi_+$ as in Theorem \ref{thm:average_+}. 
\end{cor}

\noindent
Examples of expander graphs are abundant in literature. {Specifically, it is well-known that a uniformly chosen random $d$-regular graph is an expander with 
probability tending to $1$ as its size $n \to \infty$.
Further, the edge-expander requirement of Corollary 
\ref{thm:convergence_gw_mgw}(b) holds for
the configuration models of Definition 
\ref{dfn:config}, subject only to 
uniformly bounded degree and minimal 
degree at least three. That is,}
\begin{lem}\label{lem:expander_mgw}
Suppose (\ref{eq:balance}) holds for 
$\ptheta(\cdot)$ strictly positive and $\{P_i, i \in \cQ\}$ 
of bounded support, such that $P_i (\ul{k})=0$
whenever $\|\ul{k}\| := \sum_j k_j \le 2$. Then, for 
any $0 < \delta < 1/2$ there exists $\lambda_{\delta}>0$, such that with probability tending to $1$ as $n \to \infty$,
the random graph $\Graph_n$ of Definition \ref{dfn:config}
is an  $(\delta, 1/2, \lambda_{\delta})$ edge-expander. 
\end{lem}
\noindent
In particular, 
Corollary \ref{thm:convergence_gw_mgw} holds 
for such configuration models without 
the edge-expander assumption. 

The following by-product of our 
{proof of Lemma \ref{lem:assumption_gw_mgw}} 
is of independent interest.
\begin{lem}\label{lem:recursion}
Fix {\abbr{UGW} measure with off-spring distribution $P$} of finite mean, such that 
$P([0,d_\star))=0$
for some $d_\star \ge 3$ and let 
$\blue{\widehat{\Delta}}$ be distributed on $\N$ according to 
$\blue{\widehat{P}_k}:=(k+1) P_{k+1}/\sum_j j P_j$, $k \ge 0$. 
For any fixed
$\be > \be_c$ consider the recursion over $t \ge 0$,
\beq
h^{(t+1)} \stackrel{d}{=} \sum_{\ell=1}^{\widehat{\Delta}} 
\atanh [ \tanh (\beta) \tanh (h_{\ell}^{(t)}) ]\,, \label{eq:recursion_h}
\eeq
where $h_{\ell}^{(t)}$ are i.i.d. copies 
of $h^{(t)}$ which are further independent of $\widehat{\Delta}$.
Denote by $h^{\be,+}$ its limit in law 
when $t \to \infty$ and starting at $h^{(0)}=\infty$. 
Then, fixing any
$\be \ge \be_0  > \be_\star$ and starting this recursion at 
a stochastically dominating 
{$h^{(0)} \succeq h^{\be_0,+}$}, yields a
sequence $\{h^{(t)}\}$ that converges in 
law to {$h^{\beta,+}$}.
\end{lem}

\begin{rmk}\label{rmk:splitting}
{Fixing $\be > \be_c$,
recall that any Ising Gibbs measure arising out of a 
fixed point of (\ref{eq:recursion_h})}
is a {\em splitting Gibbs measure} (see 
\cite[Remarks 1.13 and 2.6]{dembo_montanari_sun}).
Hence, Lemma \ref{lem:recursion} implies 
that there is only one {\em Bethe Gibbs measure} 
(see \cite[Remark 2.6]{dembo_montanari_sun}),  
{that corresponds to some} $h \succeq h^{\be_0,+}$, 
$\be_0 \in (\be_\star,\be)$, with a similar conclusion 
for the \abbr{UMGW} measures of Definition \ref{dfn:umgw}.

\noindent
We expect both Lemma \ref{lem:assumption_gw_mgw} and Lemma \ref{lem:recursion} to hold for \abbr{UGW} and \abbr{UMGW} measures at all $\be$ (and without a minimum degree assumption). 
However, {the
non-regularity of $\Tree$ under genuinely random 
\abbr{UGW} and \abbr{UMGW} measures yields for 
$\be \in (\be_c,\be_\star]$ a} 
technical difficulty which we can not overcome 
(c.f. Remark \ref{rmk:U_cont_difficulty}).
\end{rmk}

\noindent
{\bf Outline of the paper.} 
\begin{itemize}
\item As shown in \textsection \ref{section:conv_ising_gibbs}, 
weak convergence of $\mu_n$ {(of Definition \ref{dfn:lwc_graph})
implies that the corresponding measures} 
 $\{\nu_n\}$ and $\{\nu_{n,+}\}$ have 
sub-sequential local weak limit points
{(see Lemma \ref{lem:gm}),} 
which {subject to uniform sparseness}
are supported on the set of Ising Gibbs measures
{(see Lemma \ref{lem:conditioned_gibbs}). Both
results neither require an Ising 
model nor tree-like graphs.}

\item Relying upon the \abbr{LWC} of $\Graph_n$ to 
a law $\mu$ supported on $\cT_*$, we find
in {Lemma \ref{lem:pair_exp_conv}}
that {at its continuity points} $\U (\be,0)$ is the 
limit of both the $\nu_n$-expected values and 
$\nu_{n,+}$-expected values, 
of certain functionals of $\ul{x}$. Extending 
{(in Lemma \ref{lem:ineq_extremal}),} the result of
\cite[Lemma 3.2]{mossel_mont_sly}, we deduce 
in Lemma \ref{lem:olnu+} 
that the weak limit points of \textsection \ref{section:conv_ising_gibbs} must be convex 
combinations of $\nu_{\pm,\Tree}$
{and} 
get Theorem \ref{thm:average_free}
by the symmetry relation $\nu_n(\ul{x})=\nu_n(-\ul{x})$.
\item In \textsection \ref{section:mu_n_+} we prove 
Theorem \ref{thm:average_+}. First we deduce {in 
Lemma \ref{lem:loc_fns_conv}} 
out of \abbr{LWC} of $\Graph_n$ that the $\nu_{n,+}$-expected 
values of suitable functionals converge in expectation 
to {the corresponding values for} the limiting tree. Then, {using {in Lemma \ref{lem:boundary_ineq} and 
Lemma \ref{lem:plusminusdiff}} properties 
of \abbr{SRW} on trees}, {the assumed edge-expander
condition for $\Graph_n$ eliminates
all but one choice} for the 
convex combination of $\nu_{\pm,\Tree}$
(thus proving the theorem). 

\item In \textsection \ref{section:U_beta_continuous} we 
{deal with continuity of $\be \mapsto \U(\be,0)$.
Constructing in Lemma \ref{lem:h_cont} a suitable sequence of random variables that increases to the root 
magnetization under $\nu_{+,\Tree}$, we 
establish in Lemma \ref{lem:U_cont}
such continuity at any $\be > \be_\star$. Further,
Lemma \ref{lem:recursion} follows upon 
specializing Lemma \ref{lem:h_cont} to the context of \abbr{UGW} measures, and we provide in Lemma \ref{lem:U_cont_critical} a capacity criterion for  
continuity of $\be \mapsto \U(\be,0)$ at $\be = \be_c$, 
which we verify for \abbr{UMGW} measures.
Lastly, while Lemma \ref{lem:expander_mgw} is well known, 
for completeness we outline its proof.
}
\end{itemize}

\noindent
{\bf Acknowledgement.} We thank Allan Sly for {suggesting the weighted averages of (\ref{eq:ydef})},
Andrea Montanari for a key idea in the proof of Lemma \ref{lem:U_cont}, and Yuval Peres for helpful discussions {about Remark \ref{rem:yuval}} and the proof of Lemma \ref{lem:U_cont_critical}. We also thank Noga Alon, Russell Lyons and Nike Sun for many helpful conversations. We thank the anonymous referee for her/his helpful suggestions on improving the presentation of the paper.

\section{Convergence to Ising Gibbs measure}\label{section:conv_ising_gibbs}

We start with a general lemma about 
existence of sub-sequential local weak limits
{(based only on weak convergence of $\mu_n$
and having marks from a finite set $\cX$)}. 
\begin{lem}\label{lem:gm}
Suppose $\{\mu_n\}$ {of Definition \ref{dfn:lwc_graph}}
converges weakly in $\cP(\cG_*)$. 
Then for any probability measures $\ze_n$ on $\cX^{[n]}$  
and any sub-sequence $\{n_\ell\}_{\ell \in \N}$ there exists 
a further sub-sequence $\{n_{\ell_k}\}_{k \in \N}$ 
such that $\{\ze_{n_{\ell_k}}\}$ converges locally 
weakly to a limit $\ol{\gm}$ (which may depend on {$\{n_{\ell_k}\}$).} 
\end{lem}

\noindent
\emph{Proof}: 
{Fixing $\{\ze_n\}$ and $t \in \N$ recall that $\mu_n^t$ are such that 
\beq
\mu_n^t(\Graph) := 
{\f{1}{n}}\sum_{i=1}^n \bI (\Ball_i(t) \simeq \Graph)\,,
\notag\\
\eeq
for each $\Graph \in {\cG}_*(t)$ and the balls
$\Ball_i(t)$ in $\Graph_n$.
The assumed convergence of $\{\mu_n\}_{n \in \N}$ in $\cP(\cG_*)$ implies the convergence of 
$\{\mu_n^t\}$ in $\cP(\cG_*(t))$, so by Prohorov's 
theorem $\mu_n^t$ are uniformly tight. With 
$\cG_*(t)$ a discrete space, any compact subset
of $\cG_*(t)$ is finite, hence for any $\vep>0$ we 
have a finite set $\cG_\vep 
{(t)} \subset \cG_*(t)$, such that 
\beq\label{eq:tight-sets}
{\liminf_{n \to \infty}} \, \mu_n^t ( \cG_\vep (t) ) \ge 1- \vep
\,.
\eeq
Further, per $\Graph \in \cG_*(t)$ 
the space of marks $\cX^{\Graph}$ is finite, so 
the set $\ol{\cG}_\vep 
{(t)} := 
\{(\Graph,\ul{x}_{\Graph}) :  \Graph \in \cG_\vep {(t)}, 
\ul{x}_\Graph \in \cX^\Graph\}$ is also finite, and by Prohorov's theorem
the collection of all probability measures on 
$\ol{\cG}_\vep (t)$ is compact. In particular, 
$\cM_\vep {(t)} :=\{\delta_\Graph \otimes \nu_\Graph : 
\Graph \in \cG_\vep (t), 
\nu_\Graph \in \cP(\{-1,1\}^\Graph)\}$ is a pre-compact collection 
of probability measures on $\cP(\ol{\cG}_*(t))$.
Since $\oPb_n^t(I_n) \in \cM_\vep {(t)}$ with 
probability $\mu_n^t(\cG_\vep (t) )$, it thus follows 
that for each $t \in \N$, the laws of 
$\oPb_n^t(I_n)$ are uniformly tight, hence relatively compact. 
Consequently, there exists a diagonal sub-sequence along which the random probability measures $\oPb_n^t(I_n)$ converge in law, to say $\ol{\gm}_t$, simultaneously 
for all $t \in \N$. By the obvious 
embedding of $\ol{\cG}_*(t)$ within $\ol{\cG}_*(t+1)$, 
each ${\ol{\nu}}_{t+1} \in \cP(\ol{\cG}_*(t+1))$ 
induces a marginal probability 
measure on $\ol{\cG}_*(t)$, denoted $\pi_t({\ol{\nu}}_{t+1})$.
By definition, {$\pi_t(\oPb_n^{t+1}(I_n))=\oPb_n^t(I_n)$ 
for all $t,n \in \N$}. This implies the relation
{$\ol{\gm}_t = \ol{\gm}_{t+1} \circ \pi_t^{-1}$}
between the corresponding weak limits. That is, the
sequence $\{\ol{\gm}_t\}$ of probability measures
on the Polish spaces $\cP(\ol{\cG}_*(t))$ is consistent
with respect to the {projections $\pi_t$}. This completes 
the proof, since by Kolmogorov's extension theorem there exists a probability measure $\ol{\gm}$ on 
$\cP(\ol{\cG}_*)$ such that $\ol{\gm}_t = \ol{\gm}^t$ 
for all $t$.} 
\qed

\vskip 10pt

\noindent
Fixing {$\be \ge 0$} and $B=0$, 
with $\{\nu_n\}_{n \in \N}$ 
being Ising Gibbs measures on finite graphs $\Graph_n$, we 
wish to identify their sub-sequential limits in terms 
of Ising Gibbs measures on $\ol{\cG}_*$, 
which we define next. First recall that probability
measure $\nu_\Graph$ on {$\cX^\Graph$} for a fixed infinite
graph $\Graph \in \cG_*$ is an Ising Gibbs measure 
iff $\nu_\Graph$ satisfies the relevant {\abbr{DLR}} condition.
That is, setting $\Graph(\infty)=\Graph$, $\Graph(-1)=\emptyset$  
and $\Graph(t,\ol{t})=\Graph(\ol{t})\setminus \Graph(t)$ for 
$t<\ol{t} \le \infty$, one requires that 
for {$\ol{t}=\infty$, any $t \in \N$
and $\nu_\Graph$-a.e.
$\ul{x}_{\Graph(t,\ol{t})}$,
}
{
\beq
{\nu}_\Graph \big( \ul{x}_{\Graph (t)} \,|\, \ul{x}_{\Graph (t,\ol{t})}\big) 
= \widetilde{\nu} \big(\ul{x}_{\Graph (t)}\,|\,\ul{x}_{\Graph (t,t+1)},\Graph (t+1)\big)\,, 
\label{eq:dlr}
\eeq}
where for any finite $\Graph'=(V',E') \in \cG_*$ and 
$W \subseteq V'$, 
\beq\label{eq:fin-ising}
\widetilde{\nu} \big(\ul{x}_{W}\,|\,\ul{x}_{V' \setminus W},\Graph'
\big) 
:= \frac{\exp\Big\{\beta \displaystyle{\sum_{(i,j) \in E'} }x_i x_j\Big\}}
{\displaystyle{\sum_{\{\ul{x}'_{V'}: \ul{x}'_{V' \setminus W}=\ul{x}_{V' \setminus W}\}}}
\exp\Big\{\beta \displaystyle{\sum_{(i,j) \in E'} }x'_i x'_j\Big\}} 
\eeq
denotes the Ising measure on $W$, given boundary 
values at $V' \setminus W$ (see \cite[Chapter 2]{georgii}). 

\noindent
Next, for any 
$t \ge -1$ and $\ol{t}=t+1,\ldots,\infty$, 
fixing $r \in \N$, $\Graph \in \cG_*$ and the 
marks $\ul{x}_{\Graph(r \wedge \ol{t}) \setminus \Graph (r \wedge t)}$ we denote by
$\Ball^{(t,\ol{t})}_{(\Graph,\ul{x}_{\Graph})}(r)$ 
the union over all possible mark values $\ul{x}_{\Graph(r \wedge t)}$
of the \blue{$\ol{\cG}_*$-metric balls $\Ball_{\ol{\Graph}}(r \wedge \ol{t})$
centered at $\ol{\Graph} = (\Graph,\ul{x}_\Graph)$.}
Considering the  
sub-$\sigma$-algebras 
\beq\label{dfn:cyl-sig-alg}
\sC_{\ol{\cG}_*(t,\ol{t})}
:= \sigma \big(
\Ball^{(t,\ol{t})}_{(\Graph,\ul{x}_\Graph)}(r), 
\;\; (\Graph,\ul{x}_\Graph) \in \ol{\cG}_*,\, 
r \in \N 
\big) \,,
\eeq
{generated by these sets, note that 
$\sC_{\ol{\cG}_*(t,\ol{t})}$
are non-decreasing in $\ol{t}$ and 
non-increasing in $t$, where in} particular,
$\sC_{\ol{\cG}_*(\ol{t})} = \sC_{\ol{\cG}_*(-1,\ol{t})}$ 
and $\sC_{\cG_*} = \sC_{\ol{\cG}_*(\infty,\infty)}$ 
(as a ball $\Ball_\Graph(r) \subset \cG_*$ of radius $r$ and 
center $\Graph$ is the $\cG_*$-projection of the union 
over all $\ul{x}_\Graph \in \cX^\Graph$ 
of the corresponding balls $\Ball_{(\Graph,\ul{x}_\Graph)}(r)$
in $\ol{\cG}_*$).
%
Since $\ol{\cG}_*$ is a Polish space, 
the regular conditional probability measure
$\ol{\nu} (\cdot | \sC_{\cG_*})$
is thus well defined for any $\ol{\nu} \in \cP(\ol{\cG}_*)$
(see \cite[\textsection 9.2]{stroock}),
and we lift the notion of Ising Gibbs 
measure to $\cP(\ol{\cG}_*)$, by considering 
the \abbr{DLR} condition (\ref{eq:dlr}) with
this conditional measure playing the 
role of $\nu_\Graph$. In this setting, per $t \in \N$ 
what one has in the left-side of (\ref{eq:dlr})
amounts to the restriction to $\ul{x}_{\Graph (t)}$ 
of the regular conditional probability measure 
$\ol{\nu} ( \cdot | \sC_{\ol{\cG}_*(t,\infty)}\,)$,
resulting with the following definition.
\begin{dfn}\label{dfn:ising}
A probability measure $\ol{\nu} \in \cP(\ol{\cG}_*)$ is 
called an Ising Gibbs measure, denoted by
$\ol{\nu} \in {\cI}$, 
if for any $t \in \N$, $\ol{\nu}$-a.e. 
\beq\label{eq:dlrm}
\ol{\nu} \big(\ul{x}_{\Graph(t)} \,\big|\, 
\sC_{\ol{\cG}_*(t,\infty)} \big) 
= \widetilde{\nu} \big(\ul{x}_{\Graph (t)}\,|\,\ul{x}_{\Graph (t,t+1)},\Graph (t+1)\big)
\,,
\eeq
which we interpret as point-wise identities
in the discrete countable space $\ol{\cG}_*(t+1)$.
\end{dfn}

\noindent
\begin{rmk}
It is easy to verify from (\ref{dfn:cyl-sig-alg})
that $\sC_{\ol{\cG}_*(t,\ol{t})} \uparrow 
\sC_{\ol{\cG}_*(t,\infty)}$ as $\ol{t} \uparrow \infty$.
Thus, from L\'evy's upward theorem (applied 
point-wise on $\ol{\cG}_*(t+1)$), we have that
$\ol{\nu} \in \cI$ iff for $\ol{\nu}$-a.e.
and any $t < \ol{t} \in \N$,
{\beq
\ol{\nu} \big(\ul{x}_{\Graph (t)} \,\big|\, 
\sC_{\ol{\cG}_*(t,\ol{t})} \big) 
= \widetilde{\nu} \big(\ul{x}_{\Graph (t)}\,|\,\ul{x}_{\Graph (t,t+1)},\Graph (t+1)\big) \,.
\label{eq:dlr2} 
\eeq}
{We focus hereafter on the subset $\cI_*$
of all Ising Gibbs measures of the form 
$\ol{\nu} = \delta_\Graph \otimes \nu_\Graph$, with $\nu_\Graph$ being 
an Ising Gibbs measure for the \emph{fixed} 
graph $\Graph \in \cG_*$. Denoting by $\cI_{(t,\ol{t})}$ those 
$\ol{\nu}=\delta_\Graph \otimes \nu_\Graph$ in 
$\cP(\ol{\cG}_*)$ with $\nu_\Graph$ 
satisfying (\ref{eq:dlr}) per fixed $t < \ol{t}$ finite, we see that 
\beq\label{eq:ising-char}
\cI_* = \bigcap_{t<\ol{t}} \cI_{(t,\ol{t})} \,.
\eeq
Further, since $\ol{\cG}_*(\ol{t})$ is a discrete 
countable space, 
$\sC_{\ol{\cG}_*(t,\ol{t})}$ being a subset of its
Borel $\sigma$-algebra, 
is countably generated and the collection 
$\cI_{(t,\ol{t})}$ is completely 
determined in terms of the
marginals $\ol{\nu}^{\ol{t}}$ of probability
measures $\ol{\nu}$ on $\ol{\cG}_*$. For that reason we
hereafter take the liberty of using $\cI_{(t,\ol{t})}$
also for the subset of $\cP(\ol{\cG}_*(\ol{t}))$ consisting of the corresponding collection of 
marginals $\ol{\nu}^{\ol{t}}$.  
}
\end{rmk}

\noindent
Considering (\ref{eq:dlr}) at fixed $\ol{t}>t$ 
for $\nu_n$ and $\nu_{n,+}$, we next characterize
the sub-sequential local weak limits of $\{\nu_n\}$ 
and $\{\nu_{n,+}\}$ in terms of {certain}
Ising Gibbs measures.
\begin{lem}\label{lem:conditioned_gibbs}
Suppose $\mu_n \Ra \mu$, 
{for $\mu_n$ as in Definition \ref{dfn:lwc_graph}.}
Then,

\noindent
(a) Any sub-sequential local weak limit {$\ol{\gm}$} of 
$\{\nu_n\}$ is supported on the collection $\cI_*$ 
of Ising Gibbs measures and restricted to $\cP(\cG_*)$ 
it has the marginal $\tgm=\mu \circ \tpsi^{-1}$, where ${\tpsi}
(\Graph)=\delta_\Graph$ for any $\Graph \in \cG_*$.

\noindent
(b) The same holds for sub-sequential limits
{$\ol{\gm}_+$} of $\{\nu_{n,+}\}$, provided
$\{\Graph_n\}$ is uniformly sparse.
\end{lem}

\noindent
\emph{Proof}: 
Fix a sub-sequence $n_\ell$ along which $\{\nu_n\}$
(or $\{\nu_{n,+}\}$), converges locally weakly 
to some $\ol{\gm}$. 
{Then, for each $t \in \N$ the  
$\cP(\cG_*(t))$-restriction $\Pb_n^t(I_n)$ 
of $\oPb_n^t(I_n)$ converges in law to $\tgm^t$. 
Thus, for any fixed $\Graph \in \cG_*$,
\beq
\tgm^t(\delta_{\Graph (t)})=\lim_{\ell \ra \infty}  
\f{1}{n_\ell} \sum_{i=1}^{n_\ell} \bI (\delta_{\Ball_i(t)} = \delta_{\Graph (t)}) = \lim_{\ell \ra \infty} \f{1}{n_\ell} \sum_{i=1}^{n_\ell} 
\bI(\Ball_i(t) \simeq \Graph (t)) 
\notag
  =  \mu^t (\Graph (t)),
\eeq
where, denoting by $\mu^t$ the probability measure on $\cG_*(t)$ induced by $\mu$, the last equality follows from
the weak convergence of $\mu_n$ to $\mu$ in $\cG_*$.
Thus, for any $t \in \N$ the measure $\tgm^t$ is supported on the set of atomic
measures $\{\delta_{\Graph (t)}: \Graph \in \cG_*\}$ 
and coincides with $(\mu \circ \tpsi^{-1})^t$.
Since any probability measure 
$\gm$ on $\cP(\cG_*)$ is uniquely determined
by the collection $\{\gm^t : t \in \N\}$,
we conclude that $\tgm = \mu \circ \tpsi^{-1}$.
As for proving that $\ol{\gm} \in \cI_*$, in view of 
(\ref{eq:ising-char})}
it suffices to show that for any finite $\ol{t} > t$,
\beq\label{eq:supp-ising}
\ol{\gm}^{\ol{t}} (\cI_{(t,\ol{t})})=1 \,.
\eeq
(a) Considering first the measures $\{\nu_n\}$, recall
Definition \ref{dfn:conv_weakly} 
that $\oPb_n^{\ol{t}}(I_n)$ is supported for each $n$ 
on the collection $\{\delta_{\Ball_i(\ol{t})} \otimes \nu_{n,\Ball_i(\ol{t})} : i \in [n] \}$, where 
the restriction $\nu_{n,\Ball_i(\ol{t})}$ to 
$\Ball_i(\ol{t})$ of the Ising Gibbs measure 
$\nu_n$,
is also an Ising Gibbs measure. 
Next, per $\vep>0$ recall the finite set of graphs
$\cG_\vep(\ol{t}+1)$ we defined while proving Lemma \ref{lem:gm}, and let 
$\cG_\vep^+(\ol{t}) := \{ \Graph (\ol{t}) : 
\Graph \in \cG_\vep(\ol{t}+1) \}$, denote the corresponding
collection of one generation truncations.
Based on it, define for each $\delta \in [0,1)$,   
\beq
\label{eq:def-eps-Ising}
\cI^{\vep,\delta}_{(t,\ol{t})} 
:=
\big\{ \delta_\Graph \otimes \nu_\Graph : 
\, \Graph \in \cG^+_\vep(\ol{t}), \; 
1-\delta \le 
\frac{\nu_\Graph (\ul{x}_{\Graph (t)} \,|\, \ul{x}_{\Graph (t,\ol{t})})}
{\widetilde{\nu}(\ul{x}_{\Graph (t)}\,|\,\ul{x}_{\Graph (t,t+1)},
\Graph (t+1))} \le \frac{1}{1-\delta}
\, \big\} 
\,,
\eeq
a closed subset of $\cP(\ol{\cG}_*(\ol{t}))$.
Now, if $\Ball_i(\ol{t}+1) \simeq \Graph$ for some 
$\Graph \in \cG_\vep(\ol{t}+1)$, then   
$$
\nu_{n,\Ball_i(\ol{t})} (\ul{x}_{\Ball_i(t)}\,|\, \ul{x}_{\Ball_i(t,\ol{t})})
=\widetilde{\nu}(\ul{x}_{\Graph (t)}\,|\,\ul{x}_{\Graph(t,t+1)},\Graph
(t+1))
$$
and consequently $\oPb_n^{\ol{t}}(i) \in 
\cI^{\vep,0}_{(t,\ol{t})}$. Clearly, for any $\vep>0$ fixed, 
$\cI^{\vep,0}_{(t,\ol{t})}$ is a subset of 
$\cI_{(t,\ol{t})}$, hence from (\ref{eq:tight-sets}) 
and the assumed local weak convergence along the 
sub-sequence $n_\ell$, we deduce that 
\beq\label{eq:supp-Ieps}
1- \vep \le \limsup_{\ell \ra \infty} 
{\frac{1}{n_\ell} \sum_{i=1}^{n_\ell} 
\bI \Big\{  \oPb_{n_\ell}^{\ol{t}} (i) \in \cI^{\vep,0}_{(t,\ol{t})}
 \Big\} } \;
\le \ol{\gm}^{\ol{t}} (\cI^{\vep,0}_{(t,\ol{t})}) 
\le \ol{\gm}^{\ol{t}} (\cI_{(t,\ol{t})}) 
\,. 
\eeq
Upon considering $\vep \downarrow 0$, we conclude that
(\ref{eq:supp-ising}) holds in this case.    

\noindent
(b). For {odd} $n \in \N$ and $i \in [n]$, let
$$
Z^{0,t}_{n,i} = Z^{0,t}_{n,i} 
(\ol{t},\Ball_i(\ol{t}),\ul{x}_{\Ball_i(\ol{t})}) 
:= \nu_n \big( \, \sum_{j=1}^n x_j {\ge} 0 \,|\, 
\ul{x}_{\Ball_i(t,\ol{t})} \, \big) \,, 
$$
adopting also the notation $Z^0_{n,i}:=Z^{0,-1}_{n,i}$.
While due to conditioning on $\{\sum_j x_j {\ge} 0\}$ 
the measures $\nu_{n,+}$ are not Ising Gibbs measures, 
it is not hard to verify that for any $i \in [n]$ and
finite $\ol{t}>t$,
\beq
\nu_{n,+}\big(\ul{x}_{\Ball_i(t)}\,|\, \ul{x}_{\Ball_i(t,\ol{t})}\big)
= \frac{Z^{0}_{n,i}}{Z^{0,t}_{n,i}} \; \nu_{n,\Ball_i(\ol{t})} 
\big(\ul{x}_{\Ball_i(t)}\,|\, \ul{x}_{\Ball_i(t,\ol{t})}\big) 
\label{eq:ident_pm}
\eeq
(for clarity of presentation we ignore the 
slight modification of $Z^{0,t}_{n,i}$ which 
is required for $n$ even, in accordance with Remark \ref{rmk:even}).  
{The conditioning effect eventually 
washes away, since setting} 
$$
Z^{\pm}_{n,i} :=
\nu_n \big( \, 
\sum_{j \notin \Ball_i(\ol{t})} x_j > \pm |\Ball_i(\ol{t})| 
\;\; \big| \; \ul{x}_{\Ball_i(t,\ol{t})} \, \big) \,,
$$
{which are independent of $t<\ol{t}$, and fixing} $\vep,\delta>0$ we 
show 
that for all $n$ large enough and $i \in [n]$, 
\beq\label{eq:bd-z-pm}
\Ball_i(\ol{t}) \in \cG^+_{\vep}(\ol{t})
\quad \Longrightarrow \quad
\inf_{\ul{x}_{\Ball_i(\ol{t})}} \; \Big\{ \; 
\frac{Z^+_{n,i}}{Z^-_{n,i}} \Big\} \ge 1 - \delta \,.
\eeq
{Indeed, clearly
$Z^{+}_{n,i} \le Z^{0,t}_{n,i} \le Z^{-}_{n,i}$ and so 
by 
(\ref{eq:ident_pm})
the right-side of (\ref{eq:bd-z-pm})} 
yields that the probability measures $\P^{\ol{t}}_n(i)$
corresponding to $\nu_{n,+}$ are then in 
$\cI^{\vep,\delta}_{(t,\ol{t})}$. Consequently,
following the derivation of (\ref{eq:supp-Ieps}) 
we find that 
$1-\vep \le 
{\ol{\gm}_+^{\ol{t}}}
\big(\cI^{\vep,\delta}_{(t,\ol{t})}\big)$
for any sub-sequential limit {$\ol{\gm}_+$} of $\{\nu_{n,+}\}$
and all $\vep,\delta>0$. Since 
$$
\cI^{\vep,0}_{(t,\ol{t})} = \bigcap_{\delta>0}
\cI^{\vep,\delta}_{(t,\ol{t})} \,,
$$
considering $\delta \downarrow 0$ followed by
$\vep \downarrow 0$ completes the 
proof of (\ref{eq:supp-ising}).
{As for (\ref{eq:bd-z-pm}),  
necessarily,}
$$
\kappa := \sup_{\Graph \in \cG_{\vep}(\ol{t}+1)} \, |E(\Graph)| < \infty 
$$
(since $\cG_\vep(\ol{t}+1)$ is a 
finite collection of finite graphs).
Thus, assuming hereafter that 
$\Ball_i(\ol{t}+1) \simeq \Graph$ for some 
$\Graph \in \cG_{\vep}(\ol{t}+1)$,
at most $\kappa$ edges of $\Graph_n$ touch $\Ball_i(\ol{t})$.
{Hence, by the invariance with respect to a global 
sign change of the Ising measure
$\nu_{E_n \setminus E(\Ball_i(\ol{t}+1))}$
on the sub-graph of $\Graph_n$ in which all
edges within $\Ball_i(\ol{t}+1)$ have been deleted,
we conclude} that  
\beq\label{eq:z-minus-bdd-away}
Z_{n,i}^{-} \ge
\nu_n \big( \, \sum_{j \notin \Ball_i(\ol{t})} x_j \ge 0 
\, | \, \ul{x}_{\Ball_i(\ol{t})} \, \big)
\ge 
e^{-2 \be \kappa} 
\nu_{E_n \setminus E(\Ball_i(\ol{t}+1))} 
\big( \, \sum_{j \notin \Ball_i(\ol{t})} x_j \ge 0
\,\big)
\ge \frac{1}{2} e^{-2 \be \kappa} \,. 
\eeq
Further, 
$|\Ball_i(\ol{t})| \le \kappa$ and by the assumed uniform 
sparseness of $\{\Graph_n\}$, there exists $k \in \N$
and $n_0 \ge 3\kappa$ large enough so that
$$
\sum_{i=1}^n 
\Delta_i(\Graph_n) \bI (\Delta_i(\Graph_n) \ge k) \le \frac{n}{3}
\qquad  \quad \forall n \ge n_0  
$$
(see (\ref{eq:unif-sparse})). Consequently, for 
any $n \ge n_0$ there are at least $n/3$ vertices 
in $\Graph_n \setminus \Ball_i(\ol{t})$ of degree at most
$k-1$, out of which collection one can extract  
an independent set $S$ of $\Graph_n$ whose size is 
at least $n/(3k)$.
Thereby, one has as in the proof of  
\cite[Lemma 4.1]{mossel_mont_sly} that 
under $\nu_n$ and conditional on the 
values of $\ul{x}_{S^c}$, 
the $\pm$-valued $\{x_j\}_{j \in S}$ 
are mutually independent, each having  
expectation within $(-\eta,\eta)$ for some 
$\eta=\eta(\be,k)<1$ and all $n$. As explained there, 
the Berry-Esseen theorem then implies that for
some $C=C(k,\eta)$ finite and all $n \ge n_0$, 
$$
\sup_{r} \nu_n 
\big( \, \sum_{j \notin \Ball_i(\ol{t})} x_j = r 
\, | \, \ul{x}_{\Ball_i(\ol{t})} \, \big) \le C n^{-1/2} \,,
$$
from which it follows that uniformly in $\ul{x}_{\Ball_i(\ol{t})}$,
$$
0 \le Z^{-}_{n,i} - Z^{+}_{n,i} \le 2 |\Ball_i(\ol{t})| 
C n^{-1/2} \le 2 \kappa C n^{-1/2} \,. 
$$ 
Combining this bound 
with (\ref{eq:z-minus-bdd-away}), we 
conclude that (\ref{eq:bd-z-pm}) holds
for all $n \ge n_\delta$ sufficiently large.
\qed

\section{Identifying the limit Gibbs measure} \label{section:identify_limit}
It helps to consider in the course of our proofs 
vertex dependent magnetic fields $B_i$. That is, 
to replace the model (\ref{eq:Ising_model_descrip}) by
\beq \label{eq:model_vertex_dependent}
\nu(\ul{x})= \frac{1}{Z(\beta,\ul{B})} \exp \Big\{ \beta \sum_{(i,j) \in E} x_i x_j + \sum_{i \in V} B_i x_i\Big\}.
\eeq
In this context, we often take advantage of 
Griffith's inequality for ferromagnetic Ising models
(which for completeness we state next,
see also \cite[Theorem IV.1.21]{liggett}).
\begin{prop} {\em [Griffith's inequality]}
Consider two Ising models $\nu(\cdot)$ and $\nu'(\cdot)$ on {finite} graphs $\Graph=(V,E)$ and $\Graph'=(V,E')$, inverse temperatures $\beta$ and $\beta'$, and magnetic fields $\{B_i\}$ and $\{B_i'\}$, respectively. If $E \subseteq E'$, $\beta \le \beta'$ and $0 \le B_i \le B_i'$, for all $i \in V$, then 
$$
0 \le  \nu \Big[ 
\prod_{i \in W} x_i \Big] \le  \nu' \Big[
 \prod_{i \in W} x_i \Big] \,,
\qquad \qquad 
\forall 
\, W \subseteq V \,.
$$ 
\end{prop}
As we are having locally tree-like graphs, yielding 
local weak limit points supported on Ising Gibbs 
measures on trees, we often rely on the following 
representation for marginals of Ising measures on 
finite trees. 
\begin{prop} \label{prop:marginal_ising} \cite[Lemma 4.1]{dembo_mont}
For a subtree $\Tree'$ of a finite tree 
{$\Tree$}, let $\partial_\star \Tree'$ denote the subset of vertices 
$\Tree'$ connected by an edge to ${\sf W} := \Tree \setminus \Tree'$ and for each
$u \in \partial_\star \Tree'$ let $\langle x_u \rangle_{\sf W}$ denote the root magnetization of the Ising model on the maximal subtree $\Tree_u$ 
of ${\sf W} \cup \{ u\}$ rooted at $u$. The marginal on $\Tree'$ 
of an Ising measure $\nu$ on $\Tree$, denoted $\nu_{\Tree'}^\Tree$ is then an Ising measure on $\Tree'$ with magnetic field 
$B_u'=\atanh (\langle x_u \rangle_{\sf W}) \ge B_u$ for 
$u \in \partial_\star \Tree'$ and $B_u'=B_u$ for 
$u \in \Tree' \setminus\partial_\star \Tree'$.
\end{prop}

Adopting hereafter the notation $\Tree_{x \to y}$ for the
connected component of the sub-tree of $\Tree$ rooted at $x$,
after the path between $x$ and $y$ has been deleted, we
start by relating $\U(\be,0)$ to
the limiting correlation $x_i x_j$ 
across a uniformly chosen edge $(i,j) \in E_n$, under 
the measures $\nu_{n,\pm}$ and $\nu_n$.  
\begin{lem} \label{lem:pair_exp_conv}
Suppose
$\Graph_n \stackrel{\text{\abbr{LWC}}}{\Lra} \mu$ 
for some $\mu \in \cUs$. {Then, $(\be,B) \mapsto \U(\be,B)$ 
is bounded, non-decreasing, right-continuous at $\be,B \ge 0$,
continuous at any $B>0$, and} 
\beq \label{eq:firstequality}
\lim_{n \ra \infty}\frac{1}{n}\sum_{(i,j) \in E_n} 
\nu^{\be,0}_{n,+} \langle x_i  x_j \rangle 
= \lim_{n \ra \infty}\frac{1}{n}\sum_{(i,j) \in E_n} 
\nu^{\be,0}_{n} \langle x_i  x_j \rangle 
=  \U (\be, 0) \,,
\eeq
{at any continuity point of $\be \mapsto \U(\be,0)$.}
\end{lem}

\noindent
{\emph{Proof:}}
Since 
$\nu^{\be,0}_n = \f{1}{2} \nu^{\be,0}_{n,+} + \f{1}{2} \nu^{\be,0}_{n,-}$
and $\nu^{\be,0}_{n,-} (\ul{x}) = \nu^{\be,0}_{n,+}(-\ul{x})$ for all $\ul{x}$, clearly
$\nu^{\be,0}_{n,\pm} \langle x_i x_j \rangle
=\nu^{\be,0}_{n} \langle x_i  x_j \rangle$ for any 
$(i,j) \in E_n$ and all $n$. It thus suffices to establish 
(\ref{eq:firstequality}) in case of $\nu^{\be,0}_n$,
which since 
$$
\frac{\partial}{\partial \beta} \phi_n (\be, B)=
{\f{1}{n}}\sum_{(i,j)  \in E_n} \nu^{\be,B}_n \langle x_i  x_j \rangle \,,
$$
for all $n$, $\be$ and $B$, amounts to proving that
\beq
\lim_{n \ra \infty} \frac{\partial}{ \partial \beta} \phi_n (\be, B)=
\U(\beta, B) \,,
\label{eq:partial_beta}
\eeq
for $B=0$ and {any $\beta \ge 0$ at which 
$\U(\be,0)$ is continuous}. 
To this end, we first establish \eqref{eq:partial_beta} for all $\beta \ge 0$ and $B>0$.

\noindent
{Applying \cite[Lemma 2.12]{dembo_montanari_survey} for $A \equiv \{i,j\}$ and ${\sf U} \equiv \Ball_i(t)$, using Griffith's inequality and local weak convergence, we obtain that 
per $\be, B \ge 0$ and $t \ge 2$,} 
\begin{align}
\mu \Big[
\frac{1}{2} \sum_{i \in \partial o} \nu_{\free, \Tree}^{\be,B,t} 
\langle x_o  x_i \rangle 
\Big] 
&\le \liminf_{n \ra \infty} \frac{\partial}{\partial \beta} \phi_n (\beta, B)  \nonumber
\\
&\le \limsup_{n \ra \infty} \frac{\partial}{\partial \beta} \phi_n (\beta, B) 
\le 
\mu \Big[ 
\frac{1}{2} \sum_{i \in \partial o} \nu_{+, 
\Tree}^{\be,B,t} \langle x_o  x_i \rangle 
\Big]
\,, 
\label{eq:eq_new}
\end{align}
{where $\nu_{\free, \Tree}^{\be,B,t}$ is the Ising measure on $\Tree(t)$ with {\em free} boundary 
condition on $\partial \Tree(t){:=\Tree(t) \setminus 
\Tree(t-1)}$ 
(for more details, see \cite[pp. {163-164}]{dembo_montanari_survey}).}
Next, {for probability measures 
\beq\label{eq:edge-marginal-Ising}
\widehat{\nu} (x_1, x_2) = z^{-1}
\exp \{\beta x_1 x_2 + H_1 x_1 +H_2 x_2 \}\,, 
\eeq
on $\{-1,+1\}^2$
it is easy to check that 
\beq
\label{eq:corr-formula}
\widehat{\nu} \langle x_1 x_2 \rangle 
= F(\tanh(\beta), m_1 m_2) \,,
\eeq 
with $m_j=\tanh(H_j)$, $j=1,2$ and
\beq\label{eq:edge-corr}
F(\gth , r) := \frac{\gth + r}{1+\gth r} \,.
\eeq
Setting
$m^{\ell,\ddagger}(\Tree'):=\nu^{\be,B,\ell}_{\ddagger,\Tree'}
\langle x_{o'} \rangle$ 
for $\ddagger \in \{\free, +\}$ and the corresponding 
\emph{root-magnetization} of the Ising measure on 
$(\Tree'(\ell),o') \in \cT_*(\ell)$, we note that 
for any $i \in \partial o$, 
the marginal on ${\sf U'}=(o,i)$ of the
Ising measures $\nu^{\be,B,t}_{\ddagger,\Tree}$  
is by Proposition \ref{prop:marginal_ising}
of the form (\ref{eq:edge-marginal-Ising}),
with $m_1=m^{t,\ddagger}(\Tree_{o \ra i})$
and $m_2=m^{t-1,\ddagger}(\Tree_{i \ra o})$.
Consequently, 
$\nu_{\ddagger, \Tree}^t \langle x_o  x_i \rangle
=F(\tanh(\beta), r_{\ddagger}(t))$ is a continuous function of 
$r_{\ddagger}(t) := m^{t,\ddagger}(\Tree_{o \ra i})
m^{t-1,\ddagger}(\Tree_{i \ra o})$.
In case $B>0$, upon applying \cite[Lemma 3.1]{dommers_giardina_hofstad} ({which only requires local finiteness of the tree}), first for $\Tree=\Tree_{o \ra i}$ and then 
for $\Tree=\Tree_{i \ra o}\,$, we deduce that $r_{+}(t)-r_{\free}(t) \to 0$ and hence
$\nu_{+, \Tree}^{\be,B,t} \langle x_o  x_i \rangle -
\nu_{\free, \Tree}^{\be,B,t} \langle x_o x_i \rangle \to 0$
when $t \to \infty$. This holds for all $i \in \partial o$,
so recalling that $\odeg$ is finite
(by uniform sparseness of $\{\Graph_n\}$),}
we get by dominated convergence (\abbr{DCT}),
that 
\beq\label{eq:free-plus}
\lim_{t \ra \infty} 
\mu \Big[ 
\, 
\frac{1}{2} \sum_{i \in \partial o} \nu_{\free, \Tree}^{\be,B,t} \langle x_o x_i \rangle \,
\Big] 
= \lim_{t \ra \infty}  
\mu \Big[ 
\frac{1}{2} \sum_{i \in \partial o} \nu_{+, \Tree}^{\be,B,t} \langle x_o x_i \rangle 
\,\Big]\,,
\eeq
for any $\be \ge 0$ and $B >0$. Now using (\ref{eq:eq_new}), and recalling the definition of $\U(\be,B)$, we note that (\ref{eq:partial_beta}) holds, 
at any $B>0$ and $\be \ge 0$.

\noindent
 While (\ref{eq:free-plus})
is typically false at $B=0$ and $\be$ large enough, 
clearly for any {$\Tree \in \cT_*$ and} 
finite $t \ge 0$, the function   
$\nu_{+,\Tree}^{\beta,B,t} \langle x_o  x_i\rangle$ 
is jointly continuous in $\beta$ and $B$. 
{These Ising measures of plus boundary 
condition correspond to taking  $B_i \uparrow \infty$
at all $i \in \Tree \setminus \Tree(t-1)$ (see
Definition \ref{dfn:ising-pm}). Hence, by Griffith's inequality} 
we have that 
${\frac{1}{2}} 
\nu_{+,\Tree}^{\beta,B,t} \langle x_o x_i\rangle$ 
is non-increasing in $t$ and 
non-decreasing in $\be,B$ for $\be,B \ge 0$.
The same monotonicity properties apply for 
the sum of such functions 
over $i \in \partial o$ and in so far as $(\beta,B)$ are concerned, retained by the expectation $\U(\be,B)$
with respect to the law $\mu$ of $\Tree$, of its
limit as $t \uparrow \infty$. Since 
$\odeg$ is finite we further
deduce by \abbr{DCT} the joint continuity of 
$$
(\be,B) \mapsto 
\mu \Big[ 
\,  
\sum_{i \in \partial o} \nu_{+, \Tree}^{\be,B,t}
\langle x_o  x_i \rangle \, \Big] \,,
$$
which upon interchanging limits in $t$ and $\be,B$, yields the right-continuity of $\U(\be, B)$
at all $\be, B \ge 0$.

\noindent
We denote hereafter by 
$f_n (\cdot  )  \stackrel{{\Q}^c}{\ra} f ( \cdot  )$ 
the convergence of $f_n$ to $f$ on 
some co-countable set, and
$f( \cdot  )\stackrel{\Q^c}{=}g( \cdot  )$ 
when $f$ and $g$ agree on a co-countable set.
Since $\beta \mapsto \phi_n(\beta,B)$ are convex functions, 
so is their limit $\phi(\beta,B)$ (see 
{\cite[Theorem 1.8]{dembo_montanari_sun}} 
for existence of such limit at any 
$\beta \ge 0$, $B \in \R$ fixed).
{Such pointwise convergence of $\R$-valued convex functions 
yields that} 
$\frac{\partial}{\partial \beta} \phi_n(\beta,B) 
\stackrel{\Q^c}{\ra} \frac{\partial} {\partial \beta} \phi(\beta,B)$ per fixed $B \ge 0$, and consequently  
$\frac{\partial}{\partial \beta} \phi(\beta, B) \stackrel{\Q^c}{=} \U(\beta, B)$ at any given $B >0$. 
Fixing a sequence $B_m \decto 0$, by the convexity 
of $\be \mapsto \phi(\beta,B)$ and the continuity of
$B \mapsto \phi(\beta,B)$ we have $\frac{\partial}{\partial \beta} \phi(\beta,B_m) \stackrel{\Q^c}{\to} \frac{\partial}{\partial \beta} \phi(\beta,0)$. Further, $B \mapsto \U(\beta,B)$ 
is right continuous, hence $\U(\beta,B_m) \to \U(\beta,0)$.
From these two convergences we deduce that 
$\U(\beta,0)\stackrel{\Q^c}{=}  
\frac{\partial}{\partial \beta} \phi(\beta,0)$. 
We have seen already that 
$\frac{\partial}{\partial \beta} \phi_n(\beta,0) \stackrel{\Q^c}{\ra}  \frac{\partial}{\partial \beta} \phi(\beta,0)$, hence also $\frac{\partial}{\partial \beta} \phi_n(\beta,0) \stackrel{\Q^c}{\ra} \U(\beta,0)$. Since 
{$\frac{\partial}{\partial \beta} \phi_n(\beta,0)$
are non-decreasing} continuous functions,
this convergence extends to {\em all continuity points 
of $\be \mapsto \U(\be,0)$}. 
\qed

\medskip
The following extension of \cite[Lemma 3.2]{mossel_mont_sly}
to arbitrary $\Tree \in \cT_*$ allows us to utilize 
Lemma \ref{lem:pair_exp_conv} for restricting the 
weak limit points of $\nu_{n,+}$ and $\nu_n$, to convex combinations of $\nu_{\pm,\Tree}$.  
\begin{lem}\label{lem:ineq_extremal}
For any Ising Gibbs measure $\nu_\Tree$ on some 
$\Tree \in \cT_*$ and all $i \in V(\Tree)$,
\beq
\sum_{j \in \partial i} \nu_\Tree \langle x_i  x_j \rangle \le 
\sum_{j \in \partial i} \nu_{+,\Tree} 
\langle x_i  x_j \rangle 
=
\sum_{j \in \partial i} \nu_{-,\Tree} \langle x_i  x_j \rangle \,,
\label{eq:ineq_eq}
\eeq
with strict inequality \emph{for some $i \in V(\Tree)$} 
unless $\nu_\Tree$ is a convex combination of $\nu_{+,\Tree}$ and $\nu_{-,\Tree}$.
\end{lem}

\noindent
{\em Proof.} The equality in (\ref{eq:ineq_eq}) 
is an immediate consequence of 
the fact that under $\nu_{+,\Tree}$ the random 
vector $-\ul{x}_{\Tree}$ admits the law $\nu_{-,\Tree}$.
Further, due to uniqueness of the Ising Gibbs measure 
for a finite $\Tree$, we may and shall consider hereafter 
a fixed {\em infinite tree} $\Tree$. There are
only countably many edges in $\Tree$ and the non-empty 
collection of Ising Gibbs measures on $\Tree$ is convex, 
with each Ising Gibbs measure on $\Tree$ being a 
mixture of the extremal Ising Gibbs measures on 
$\Tree$ (see \cite[Chapter 7]{georgii}). Consequently, 
it suffices to fix an {\em extremal} Ising
Gibbs measure $\nu_\Tree \neq \nu_{\pm,\Tree}$ 
and show that for every edge $(i,j) \in E(\Tree)$,
\beq\label{eq:ineq_edge}
\nu_\Tree \langle x_i  x_j \rangle
\leq \nu_{+,\Tree} \langle x_i x_j \rangle 
\eeq
with a strict inequality for 
{\em at least one} $(i,j) \in E(\Tree)$.
To this end, 
for each $(i,j) \in E(\Tree)$ let 
$m^{\nu}_{i \to j} := \nu_\Tree^{(ij)} \langle x_i \rangle$ 
for the probability measure $\nu_\Tree^{(ij)}$  
whose Radon-Nikodym 
derivative with respect to $\nu_{\Tree}$ 
is proportional to $e^{-\be x_i x_j}$.
That is, 
$$
m_{i \to j}^\nu  
= \frac{\nu_\Tree \langle x_i e^{-\be x_i x_j} \rangle}
{\nu_\Tree \langle e^{-\beta x_i x_j} \rangle}
= \lim_{l \to \infty}
\frac{\nu_\Tree \langle x_i e^{-\be x_i x_j} | 
\ul{x}_{\Ball_i(l)^c} \rangle}
{\nu_\Tree \langle e^{-\beta x_i x_j} | 
\ul{x}_{\Ball_i(l)^c} \rangle} \,,
$$ 
where the limit exists by backward martingale convergence 
theorem and is a.e. constant by the tail triviality of 
the extremal measure $\nu_\Tree$ (see \cite[Chapter 7]{georgii}). Using the \abbr{DLR} condition (\ref{eq:dlr}) 
for $\nu_\Tree$ and the tree structure of $\Tree$, 
we deduce that 
$\nu_\Tree$-a.e. 
\beq\label{eq:mij}
m_{i \ra j}^\nu = \lim_{l \ra \infty} \widetilde{\nu} 
{\langle x_i} \,|\, \ul{x}_{\Tree_{i\ra j}(l,l+1)}, 
\Tree_{i \to j}(l+1) \rangle \,.
\eeq 
By the \abbr{DCT}, the \abbr{DLR} condition (\ref{eq:dlr}) 
for $\nu_\Tree$, Proposition \ref{prop:marginal_ising} 
and (\ref{eq:mij}), for each $t \in \N$ the 
marginal law of $\ul{x}_{\Tree(t)}$ under $\nu_\Tree$ 
is completely determined by 
${\{m_{i \to j}^{\nu}, i \in \partial \Tree(t),
j \in \partial \Tree(t-1) \}}$.
In particular, considering the formula (\ref{eq:corr-formula}), we get by the same line of reasoning that 
\beq\label{eq:Ising-corr-formula}
\nu_\Tree \langle x_i x_j \rangle =
F(\tanh(\beta), m_{i \ra j}^\nu m_{j \ra i}^\nu) \,,
\eeq
for $F(\gth,r)$ of (\ref{eq:edge-corr}) 
and any $(i,j) \in E(\Tree)$,
with the analogous expression in case of 
$\nu_{+,\Tree} \langle x_i  x_j \rangle$.
Denoting by $m_{i \ra j}^-$ and $m_{i \ra j}^+$ the values 
of $m_{i \ra j}^\nu$ for Ising Gibbs measures
$\nu_{-,\Tree}$ and $\nu_{+,\Tree}$, respectively, from (\ref{eq:mij}) and Griffith's inequality 
we know that $|m_{i \to j}^\nu| \leq m_{i \to j}^+$ 
for all $(i,j) \in E(\Tree)$, out of which we get 
the inequality (\ref{eq:ineq_edge}) by the strict monotonicity of $r \mapsto F(\gth,r)$ on $[-1,1]$ 
(when $|\gth| < 1$). Turning to prove that having 
equality in (\ref{eq:ineq_edge}) for all $(i,j) \in E(\Tree)$ implies either $\nu_\Tree = \nu_{+,\Tree}$ or 
$\nu_\Tree = \nu_{-,\Tree}$, note that by the preceding 
{such an equality in \eqref{eq:ineq_edge} translates} into 
\beq
m_{i \ra j}^\nu m_{j \ra i}^\nu = m_{i \ra j}^+ m_{j \ra i}^+ \qquad 
\forall (i,j) \in E(\Tree) \,. \label{eq:a_crucial_step}
\eeq
From (\ref{eq:mij}) one also have by an 
explicit calculation for Ising 
measures on trees, that 
\beq\label{eq:recursion_m}
m_{i \to j}^\nu = \tanh \big[ 
\sum_{k \in \partial i \setminus \{j\}} \atanh(\tanh(\beta) 
m_{k \to i}^\nu) \big] 
\qquad 
\forall (i,j) \in E(\Tree) \,,
\eeq
with the same recursion holding for the 
collections $\{m_{i \to j}^{\pm}, (i,j) \in E(\Tree)\}$.
Suppose now that some $(i,j) \in E(\Tree)$ is a 
\emph{plus edge}, namely both
$m^{\nu}_{i \to j} = m^{+}_{i \to j}$ and
$m^{\nu}_{j \to i} = m^{+}_{j \to i}$. 
Out of (\ref{eq:recursion_m}) we have that 
$m_{i \to j}^\nu$ is \emph{strictly} increasing in 
each $m_{k \to i}^\nu$, $k \in \partial i \setminus \{j\}$,
so with $|m^{\nu}_{k \to i}| \le m^+_{k \to i}$, the 
assumed equality {$m_e^\nu=m_e^+$ at both 
directed edges $e = \{i \to j\}$ and
$e=\{ j \to i\}$,} implies the same at all directed 
edges $k \to i$,
$k \in \partial i$. Further, by (\ref{eq:recursion_m}) 
the values of $m^\nu_{i \to k}$ and $m^+_{i \to k}$ 
are given by the same function of $\{m^\nu_e\}$
and $\{m^+_e\}$ respectively, whose arguments are   
directed edges $e$ where we already have $m_e^\nu=m_e^+$.
Hence, that equality holds also for  
all directed edges of the form $e=\{i \to k\}$. 
That is, every 
edge of {$\Ball_i(1)$} is a plus edge. This property 
extends in the same manner to {$\Ball_i(t)$}, 
$t=2,3,\ldots$, and so we conclude that a single plus 
edge in $\Tree$ results with each edge being plus 
edge, and thereby with $\nu_\Tree=\nu_{+,\Tree}$.
By the same line of reasoning, a single 
\emph{minus edge} $(i,j)$ where both 
$m^{\nu}_{i \to j} = m^{-}_{i \to j}$ and
$m^{\nu}_{j \to i} = m^{-}_{j \to i}$ yields
that all edges of $\Tree$ are minus edges and 
thereby $\nu_\Tree=\nu_{-,\Tree}$. Suppose now 
that there are neither plus nor minus edges 
in $\Tree$. We then have by (\ref{eq:a_crucial_step}) 
that at each edge $(i,j)$ either $m^+_{i \to j}>0$ 
and $m^+_{j \to i}=0$, or the same applies upon 
reversing the roles of $i$ and $j$, and we thus 
complete the proof by ruling out the possibility 
of $\nu_{+,\Tree}$ having the latter property.
Indeed, by (\ref{eq:recursion_m}) if some $m^+_{i \to j}>0$ 
then $m^+_{l \to i}$ is strictly positive for at least
one edge $(l,i)$ of $\Tree$. The latter is neither 
plus nor minus edge, so $m^+_{i \to l}=0$, which 
with $m^+$ everywhere non-negative, implies by (\ref{eq:recursion_m}) that $m^+_{k \to i} = 0$ 
at all $k \in \partial i \setminus \{l\}$. 
That is, having $m^+_{i \to j}>0$ results with
{$m_e^+$} strictly positive at exactly one edge 
{$e$} directed into $i$. Continuing in this 
manner we find an infinite directed ray 
$\{i_s \to i_{s-1}: (i_s,i_{s-1}) \in E(\Tree), 
s \in \N \}$ (ending at $i_1=i$ and $i_0=j$), 
with $m^+_{i_s \to i_{s-1}}>0$ while 
$m^+_{k \to i_{s-1}}=0$ for all $k \ne i_s$, $s \ge 1$. 
That is, again 
by (\ref{eq:recursion_m}), $m^+_{i_s \to i_{s-1}} = 
\tanh(\beta) m^+_{i_{s+1} \to i_s}$ for all $s \ge 1$.
With $\tanh(\beta) < 1$ it is obviously impossible to 
have such an infinite sequence of strictly positive 
$m^+_{i_s \to i_{s-1}} \le 1$.
\qed

\begin{rmk}
Unlike the case of $k$-regular trees $\Tree_k$ 
considered in \cite[Lemma 3.2]{mossel_mont_sly},
we may have
$$
\sum_{i \in \partial o} \nu_\Tree \langle x_o  x_i \rangle = \sum_{i \in \partial o} 
\nu_{+,\Tree} \langle x_o  x_i \rangle 
$$
for some $\Tree \in \cT_*$ and an extremal Ising 
Gibbs measure $\nu_\Tree \ne \nu_{\pm,\Tree}$ on it.
Indeed, as the proof of Lemma \ref{lem:ineq_extremal} shows,
this happens whenever {$\be>0$ is such that 
for some $i \in \partial o$ 
there is a unique Ising Gibbs measure on the 
sub-tree $\Tree_{o \to i}$ while}
$\Tree' := \Tree_{i \to o}$ admits an extremal 
Ising Gibbs measure other than $\nu_{\pm,\Tree'}$
{(e.g. when $\Tree_{i \to o}$ is $k_2$-regular, while
$\Tree_{o \to i}$ is finite or $k_1$-regular and 
$\be_c(k_2) < \be < \be_c(k_1)$).}
Nevertheless, our next lemma utilizes the unimodularity of 
$\mu$ to circumvent this problem.
\end{rmk}
\begin{lem} \label{lem:ineq}
{Fixing 
$\mu \in {\cUs}$, for any 
$\ol{\gm}$ supported on the collection $\cI_*$ of 
Ising Gibbs measures $\ol{\nu} = 
\delta_{\Tree} \otimes \nu_\Tree$ and having the 
law $\mu$ for $\Tree$,}
\beq 
\ol{\gm} \Big[ \sum_{ i \in \partial o} \nu_\Tree \langle
x_o x_i \rangle \Big]  \le \mu\Big[\sum_{ i \in \partial o} \nu_{+,\Tree} \langle x_o  x_i \rangle \Big] = \mu \Big[\sum_{ i \in \partial o} \nu_{-,\Tree} \langle x_o x_i \rangle\Big], \label{eq:ineq}
\eeq
with strict inequality 
unless $\ol{\gm}$ is supported on the sub-collection 
$\cI_\pm \subset \cI_*$ of 
$\ol{\nu} = \delta_\Tree \otimes \nu_\Tree$ where
\beq
\nu_\Tree = \alpha_\Tree \nu^\be_{+,\Tree} + (1 - \alpha_\Tree) \nu^\be_{-, \Tree} 
\label{eq:limit_pm_combination}
\eeq
for some Borel measurable function 
$\alpha: \cT_* \mapsto [0,1]$. Further, w.l.o.g. we take hereafter 
$\alpha_\Tree=\f{1}{2}$ on the set $\{\Tree \in \cT_*: \, \nu^\be_{+,\Tree} = 
\nu^\be_{-,\Tree} \}$.  
\end{lem}

\begin{rmk}\label{rw-ergodicity}
In the proof of Lemma \ref{lem:ineq} we take advantage of 
the {$\cT_*$-valued Markov chain \blue{$\{\wt{Y}_\ell\}$} 
commonly known as 
``walk from the point of view of the particle'',
induced by setting the root of $\Tree$ to follow 
the path of
\blue{discrete time simple random 
walk (\abbr{DSRW})} \blue{$\{Y_\ell\}$} of law $\blue{\wh{\P}_o^\Tree}$
on $(\Tree,o) \in \cT_*$,
starting at $Y_0=o$. Specifically, associating with
each $\mu \in \cUs$ for which $\odeg>0$,
the ``size-biased-root'' probability measure 
$\blue{\widehat{\mu}}:= \frac{\Delta_o}{\odeg} \mu$
and choosing 
$\wt{Y}_0 \in \cT_*$ 
according to $\widehat{\mu}$, yields the stationary and 
reversible joint law 
$\blue{\widehat{\mu} \otimes \widehat{\P}_o^\Tree}$ 
for the trajectory $\{\wt{Y}_\ell\}$ 
(c.f. \cite[Theorem 4.1]{aldous_lyons}).} 
\end{rmk}

\noindent\emph{Proof}: We get (\ref{eq:ineq}) 
by considering the expectation of (\ref{eq:ineq_eq}) 
for $i = o$, over the law $\ol{\gm}$ of $\Tree$ and 
the Ising Gibbs measure $\nu_\Tree$ on it. Further, 
there is only one Ising Gibbs measure on $\Tree=\{o\}$.
So, our claim about strictness of the 
inequality in (\ref{eq:ineq}) trivially holds in case $\odeg=0$,
and assuming hereafter that $\odeg>0$, 
we consider the $\cT_*$-valued stationary Markov chain 
$\{\wt{Y}_\ell\}$, as in Remark \ref{rw-ergodicity}.
Let $\ol{\nu}_\Tree$ denote the  
expected value of $\nu_\Tree$ under the
probability measure $\ol{\gm}$ conditional
upon $\Tree \in \cT_*$, which 
up to some $\mu$-null set $\cN \subset \cT_*$ 
is a uniquely defined Ising Gibbs measure 
on $\Tree$ (due to convexity of the 
latter collection).
Equality in (\ref{eq:ineq}) thus 
amounts to $\E [f(\wt{Y}_0)] = 0$
for the $\cT_*$-measurable, uniformly bounded and 
non-negative (see (\ref{eq:ineq_eq})),    
\beq
f((\Tree,o)) := \f{1}{\Delta_o}
\sum_{j \in \partial o} 
\big[ \nu_{+,\Tree} \langle x_o x_j \rangle
-  \ol{\nu}_{\Tree} \langle x_o  x_j \rangle \big] \,, 
\notag
\eeq
which by the stationarity of $\{\wt{Y}_\ell\}$ 
(c.f. \cite[Theorem 4.1]{aldous_lyons}), implies that 
\beq
\E[f(\wt{Y}_\ell)] = 0 \qquad \forall \ell \in \N\,. 
\label{eq:mc_eq1}
\eeq
Conditional on $\wt{Y}_0=(\Tree,o)$,
the probability of $\wt{Y}_\ell=(\Tree,i)$ is strictly positive
for each $\ell \in \N$ and $i \in \partial \Tree(\ell)$, 
hence with $f(\cdot)$ non-negative, it follows from  
(\ref{eq:mc_eq1}) that 
$$
\widehat{\mu} \big( (\Tree,o) \in \cT_*:
\exists i \in V( \Tree), \;\; f((\Tree,i))>0 \big) = 0 \,.
$$
We thus conclude that for $\mu$-a.e. $\Tree$,  
equality holds in (\ref{eq:ineq_eq})  
for $\ol{\nu}_\Tree$ and {\em all $i \in V(\Tree)$},
so by Lemma \ref{lem:ineq_extremal} 
the Ising Gibbs measure $\ol{\nu}_\Tree$ must 
then be a convex combination of $\nu_{+,\Tree}$ 
and $\nu_{-,\Tree}$.  
{Now recall that to any Ising Gibbs measure 
$\nu_\Tree$ on $\Tree$ corresponds a unique 
probability measure 
$\Theta_{\nu_\Tree}$ supported on the collection 
$\{\nu^e_\Tree\}$ of extremal Ising Gibbs measures  
on $\Tree$, such that 
$\nu_\Tree(\cdot) = \int \nu^e_{\Tree}(\cdot) d \Theta_{\nu_\Tree}$ (c.f. \cite[Theorem 7.26]{georgii}). Therefore, by its definition, $\mu$-a.e.  
$\ol{\nu}_\Tree (\cdot) = \int \nu_{\Tree}^e(\cdot) d\Theta_\Tree$ for the expected value 
$\Theta_\Tree(\cdot)$ of $\Theta_{\nu_\Tree}(\cdot)$ 
under the probability measure $\gm$ conditional 
upon $\Tree \in \cT_*$. We have just shown that $\mu$-a.e. $\Theta_\Tree (\{\nu_{+,\Tree}, \nu_{-,\Tree}\}^c)=0$,
hence $\gm$-a.e. this holds for $\Theta_{\nu_\Tree}$. 
That is, up to some $\gm$-null set, $\nu_\Tree$ is 
of the form (\ref{eq:limit_pm_combination}), as claimed.}
\qed

\vskip10pt

\noindent
The following lemma completes the proof of Theorem \ref{thm:average_free}. 
\begin{lem}\label{lem:olnu+}
Under the conditions of Lemma \ref{lem:pair_exp_conv}
we have that:

\noindent
(a). Any sub-sequential local weak limit $\ol{\gm}_+$ 
of $\{\nu_{n,+}\}$ is supported on the collection
$\cI_\pm$,
with $\Tree$ distributed according to $\mu$.

\noindent
(b). Any sub-sequential local weak limit of $\{\nu_n\}$ 
equals $\ol{\gm}= \mu \circ \opsi^{-1}$ for 
$\opsi(\Tree)= \delta_\Tree \otimes
(\f{1}{2} \nu^\be_{+,\Tree} + \f{1}{2} \nu^\be_{-, \Tree})$.
\end{lem}

\noindent
{\em Proof}: (a). {Recall Lemma \ref{lem:pair_exp_conv},
that  
\beq\label{eq:conv-corr1}
\frac{2}{n} \sum_{(i,j) \in E_n} \nu_{n,+}^{\be,0} 
\langle x_i  x_j \rangle 
= \mu_n \big[ F(\oPb^{{2}}_n(I_n)) \big] \to 
\mu \big[ 
\sum_{i \in \partial o} \nu_{+,\Tree} \langle x_o x_i \rangle \big] \,,
\eeq
for $\oPb_n^t(i)$ corresponding to $\nu_{n,+}$ and 
the function 
$
F(\ol{\nu}) := \ol{\nu} \big\langle 
\sum_{i \in \partial o} x_o x_i \big\rangle$
on $\cP(\ol{\cG}_*{(2)})$, which 
is bounded by $\ol{\nu}(\Delta_o)$ and continuous with 
respect to weak convergence.
By   assumption, 
under $\mu_n$ the law of $\oPb_{n}^{2}(I_{n})$ converges weakly to 
$\ol{\gm}^2_+$ along some sub-sequence $n_\ell \to \infty$.
Hence, by \abbr{DCT} and the
uniform integrability of $\{\Delta_{I_n}\}$, 
\beq\label{eq:conv-corr2}
\lim_{\ell \ra \infty} \mu_{n_\ell} \big[ F(\oPb^{{2}}_{n_\ell}(I_{n_\ell}))\big]
= \ol{\gm}_+ \big[ F(\ol{\nu}^{2}) \big]. 
\eeq
Recall part (b) of Lemma \ref{lem:conditioned_gibbs}
that $\ol{\gm}_+$ is supported on the collection $\cI_*$
of Ising Gibbs measures of the form 
$\delta_\Tree \otimes \nu_\Tree$, having the 
law $\mu \in \cU$ for $\Tree \in \cT_*$.
Thus, comparing the \abbr{RHS} of (\ref{eq:conv-corr1}) with 
the \abbr{RHS} of (\ref{eq:conv-corr2}), we deduce that 
\beq
\ol{\gm}_+ \big[ 
\sum_{i \in \partial o} \nu_\Tree \langle 
x_o x_i \rangle \big] = 
\mu 
\big[\sum_{i \in \partial o} \nu_{+,\Tree} \langle x_o x_i \rangle \big] 
\notag
\eeq
out of which it follows
by Lemma \ref{lem:ineq} that $\gm_+$ is supported on
the sub-collection $\cI_\pm$. 

\noindent
(b). Considering now part (a) of Lemma \ref{lem:conditioned_gibbs} we get by the preceding 
argument that any sub-sequential weak limit $\ol{\gm}$ 
of $\{\nu_n\}$ is supported on $\cI_\pm$ with $\Tree$
distributed according to $\mu$. In particular, $\ol{\gm}$-a.e. 
$$
|\nu_\Tree \langle x_o \rangle| = 
|2\alpha_\Tree-1| \nu_{+,\Tree} \langle x_o \rangle \,.
$$
As in the proof of Lemma \ref{lem:ineq_extremal}, if 
$\nu_{+,\Tree} \langle x_o \rangle = 0$, then necessarily 
$m^+_{i \to j}=0$ for all $(i,j) \in E(\Tree)$, hence
$\nu_{+,\Tree} = \nu_{-,\Tree}$ and {by our convention}  
$\alpha_\Tree = \frac{1}{2}$. More generally, 
the bounded function $\wt{F}(\ol{\nu}) := 
|\ol{\nu} \langle x_o \rangle|$ on $\cP(\ol{\cG}_*{(1)})$
is 
continuous with respect to weak convergence.
{Since} 
$$
0 = n^{-1} \sum_{i=1}^{n} |\nu_{n} \langle x_i \rangle| 
= \mu_n \big[ \wt{F} (\oPb^{{1}}_n(I_n)) \big]
$$
{for all $n$, it thus follows that}
for any local weak limit point $\ol{\gm}$ of $\{\nu_n\}$,
\beq
0 = \ol{\gm} [\wt{F}(\ol{\nu}^{{1}})] =
\ol{\gm} \big[ \, |\, \nu_\Tree \langle x_o \rangle \, | \, \big] = 
\ol{\gm} \Big[ \, 
| 2 \alpha_\Tree - 1 | \nu_{+,\Tree} \langle x_o \rangle \, \Big] \,, 
\notag
\label{eq:nu_mod}
\eeq 
thereby forcing $\ol{\gm}$-a.s. $\alpha_\Tree = \f{1}{2}$.}  
\qed

\section{Proof of Theorem \ref{thm:average_+}}\label{section:mu_n_+}

\noindent
Given part (a) of Lemma \ref{lem:olnu+} it 
remains only to show that $\ol{\gm}_+$-a.s., we may take 
{$\alpha_\Tree = 1$}
for any sub-sequential local weak limit point 
$\ol{\gm}_+$ of $\{\nu_{n,+}\}$. To this end, 
we make use of the following definition.
\begin{dfn} \label{dfn:conv_average} 
{Given graphs $\{\Graph_n\}_{n \in \N}$ having
vertex sets $V_n=[n]$ and probability measures 
$\ze_n$ on ${\cX^{[n]}}$, let 
$\oPb_n^t \in \cP(\ol{\cG}_*(t), \sC_{\ol{\cG}_* (t)})$ 
denote the average over a uniformly chosen $I_n \in [n]$,
of the law of $(\Ball_{I_n}(t),\ul{x}_{\Ball_{I_n}(t)})$, 
for a positive integer $t$ and
$\ul{x}$ drawn according to $\ze_n$}
(i.e. 
$\oPb_n^t = \mu_n (\oPb_n^t(I_n))$
for $\oPb_n^t(i)$ of Definition \ref{dfn:conv_weakly}).
We say that $(\Graph_n,\ze_n)$, or 
in short, that $\ze_n$, 
converge {\em on average} to $\ol{\nu}$, 
a probability measure on $(\ol{\cG}_*, \sC_{\ol{\cG}_*})$, 
if for any fixed positive integer $t$,
\beq
\oPb_n^t \Ra \ol{\nu}^t,\quad \text{ as } n \ra \infty\,.
\eeq
\end{dfn}

\begin{rmk}\label{rmk:conv-av}
Note that if $\{\ze_n\}$ converges locally weakly to 
$\ol{\gm}$ then it also converges on average to 
$\ol{\nu} = \int \nu \ol{\gm} (d\nu)$.
In particular, if $\ol{\gm}$ is supported on 
the subset $\cI_\pm$
of Ising Gibbs measures 
then 
it follows by linearity of the 
conditional expectation that the corresponding 
limit on average $\ol{\nu}$ of $\{\ze_n\}$, is 
itself an Ising Gibbs measure,
with $\Tree$ distributed according to the $\cP(\cT_*)$-marginal of 
{$\ol{\gm}$} and
$\ol{\nu}_\Tree$ of the form (\ref{eq:limit_pm_combination}) for some measurable $\alpha : \cT_* \mapsto [0,1]$.
\end{rmk}

\noindent
Given $\Graph \in \cG_*$ and $X_0\in V(\Graph)$, let \blue{$\{X_s\}$} denote the 
\blue{variable speed continuous time simple random walk (\abbr{VSRW})} on $\Graph$, 
i.e. the Markov jump process of state space $V(\Graph)$,
which upon arriving at any $j \in V(\Graph)$, jumps with unit rate to each 
possible $j' \in \partial j$.
Now, for $r \in \N$, $l >0$ 
and $i \in V(\Graph)$, let $a_{i,j}^{l,r,\Graph}$ denote 
the expected {relative to $l$} occupation 
time at $j \in V(\Graph)$ by {such} 
\abbr{VSRW} $\{X_s\}$ on $\Graph$ {which starts} 
at $X_0=i$ and run till $\min(l,\theta_r)$ for
$\theta_r := \inf \{s \ge 0 : X_s \not\in \Ball_i(r)\}$. 
That is, 
with $\P_i^{\Graph}$ denoting  
the law of \abbr{VSRW} on the fixed $\Graph$, starting at $X_0=i$,
\beq\label{eq:adef}
a_{i,j}^{l,r,\Graph} := \frac{1}{l} \int_0^l 
\P_i^{\Graph} (X_s =j, s \le \theta_r) \ ds\,.
\eeq
These non-negative weights induce for every 
$\ul{x} \in \cX^{V(\Graph)}$ the weighted averages
\beq\label{eq:ydef}
y_i^{l,r,\Graph}(\ul{x}):=\sum_j x_j a_{i,j}^{l,r,\Graph} \,,
\eeq
{having mean value 
\beq\label{eq:rho-def}
\blue{\rmg_i^{l,t,r,\Graph}} :=
\nu_{+,\Ball_i(t)} \langle y_i^{l,r,\Graph} \rangle \,,
\eeq
under the Ising measure $\nu_{+,\Ball_i(t)}$ 
on $(\Graph,i)$, at parameters $(\be,0)$, 
conditioned to 
$\ul{x}_{\Ball_i(t)^c} = (+)_{\Ball_i(t)^c}$.}
Our proof is based on analyzing per {$\eta \in (0,1)$} and $t \in \N$,
the functionals
\begin{align}\label{eq:FA-def}
\blue{\Fsp_i^{l,t,r,\Graph}(\ul{x},\eta)}
&:= \blue{\Jsp_i^{l,r,\Graph}(\ul{x},\eta) \Amg_i^{l,t,r,\Graph}(\eta)}
\,,
\\
\blue{\Jsp_i^{l,r,\Graph}} (\ul{x},\eta) &:= \bI \Big\{ 
y_i^{l,r,\Graph}(\ul{x}) \le - {\eta} \Big\} \,, 
\qquad
\label{eq:fn2_choice}
\blue{\Amg_i^{l,t,r,\Graph} (\eta) := \bI
\Big\{
\rmg_i^{l,t,r,\Graph}  
\ge 
2\eta \Big\}}\,.
\end{align}
In doing so, we use {$a_{i,j}^{l,r,n}$, 
$y_i^{l,r,n}$, $\rmg_i^{l,t,r,n}$ (\blue{$\Jsp_i^{l,r,n}$, $\Amg_i^{l,t,r,n}$, 
$\Fsp_i^{l,t,r,n}$}), when $\Graph=\Graph_n$ 
and similarly $a_{i,j}^{l,r,\Tree}$,
$y_i^{l,r,\Tree}$, $\rmg_i^{l,t,r,\Tree}$ ($\Jsp^{l,r,\Tree}_i$, 
$\Amg_i^{l,t,r,\Tree}$, $\Fsp_i^{l,t,r,\Tree}$)} when
$\Graph=\Tree \in \cT_*$,
omitting $r$ and $t$ in case $r=\infty$ (respectively, 
$t=\infty$, which for $\Amg_i^{l,r,\Tree}$ means using $\nu_{+,\Tree}$), 
and arguments $\eta$, $\Graph$, $\ul{x}$ whose value is 
clear from the context. 

\smallskip
\blue{To explain the role of the various quantities 
introduced in \eqref{eq:ydef}--\eqref{eq:fn2_choice}, 
recall that for $k$-regular graphs \cite{mossel_mont_sly} 
fix $\eta>0$ small so 
the indicators $\Jsp_i^{l,n}$ identify 
vertices $i \in \Graph_n$ in the ``- state'' of each configuration 
$\ul{x}$, while the conditioning inherent to $\nu^\be_{n,+}$ keeps  
at least $\frac{\eta}{2} n$ vertices $i \in \Graph_n$ out of this state. 
If we take $a_{i,j} = |\Ball_i(l)|^{-1} \bI_{\{j \in \Ball_i(l)\}}$ 
{in \eqref{eq:ydef}, 
as \cite{mossel_mont_sly} do, then due to the variability 
of ball sizes $|\Ball_i(l)|$ across $i \in \Graph_n$,
we would no longer find a clear relation between $\sum_i y_i$ 
and the value $\sum_j x_j$ on which we conditioned.
We resolve this problem by using instead 
the weights of \eqref{eq:adef} and taking advantage of 
the reversibility of the \abbr{vsrw}. Indeed, as we show 
next, then within the support 
of $\nu_{n,+}$ one has at least $\frac{\eta}{2} n$ 
vertices} $i \in \Graph_n$ for which $\Jsp_i^{l,n}(\ul{x},\eta)=0$
(and hence $\Fsp_i^{l,t,n}(\ul{x},\eta)=0$).}
\begin{lem} \label{lem:u_n_small}
For any {$\eta \in (0,1)$,} 
$l \ge 0$, $n \in \N$ and $\ul{x}$ such that $\sum_j x_j \ge 0$, 
\beq
{\mu_n} \big[ 1 - \Jsp_{I_n}^{l,n}(\ul{x},\eta) \big] 
\ge \frac{\eta}{2} \,.
\notag
\eeq
\end{lem}

\noindent
\emph{Proof:} Since $\sum_k a_{j,k}^{l,n}=1$ for any $n,l$, 
we have that $\Jsp_i^{l,n}(\ul{x},\eta) = 
\bI_{\{z_i \ge 1 + {\eta} \}}$ for the
non-negative $z_i := \sum_j (1-x_j) a_{i,j}^{l,n}$.  
Further, due to reversibility of the \abbr{VSRW},
$a_{i,j}^{l,n}= a_{j,i}^{l,n}$ 
for all $i,j \in V_n$. Hence, 
by our assumption that $\sum_j x_j \ge 0$,
\beq
\mu_n [z_{I_n}] = \frac{1}{n}
\sum_{k,j=1}^n (1-x_j) a_{k,j}^{l,n} 
= \frac{1}{n} \sum_{j=1}^n (1-x_j) \le 1 \,. \notag
\eeq
Thus, applying Markov's inequality to $z_{I_n}$
completes the proof.
\qed

\blue{In the regular case \cite{mossel_mont_sly} show that for $\be > \be_c$, 
if $\ell \gg 1$ then $\Jsp_i^{l,n} = \Jsp_j^{l,n}$ for most $(i,j) \in E_n$
which by the assumed edge-expander properties of $\Graph_n$
forces every limit point of $\nu_{n,+}^\be$ to have 
$\frac{\eta}{2} \ge 1 - \alpha_{\Tree_k}$ (so taking $\eta \to 0$ 
completes their proof). {To make this argument work, one needs
that as $l \to \infty$ the means $\rmg_i^{l,\Tree}$ 
be uniformly bounded away from zero, for $\mu$-a.e. $\Tree$.
We have the latter property for $\be>\be_c$,
provided that $\mu$ is an {\em extremal} element of $\cUs$, 
since then $\rmg_o^{l,\Tree}$ converges as $l \to \infty$ to the strictly 
positive expected magnetization
\beq\label{eq:rho-lim-def}
\rmg_\mu := \mu [\nu_{+,\Tree} \langle x_o \rangle] 
\eeq
(see Lemma \ref{lem:boundary_ineq}
and
Remark \ref{rem:yuval}).
However, for general $\mu \in \cUs$ we have no non-trivial uniform 
asymptotic lower bound on $\rmg_i^{l,\Tree}$, 
so use the indicators $\Amg_i^{l,n}$ for masking out in \eqref{eq:FA-def} 
those $i \in \Graph_n$ for which $\Ball_i(t)$ converges
to a tree $\Tree$ of too small mean (and we later dispense of this
masking effect by taking $\eta \to 0$).}
 
Both for utilizing the reversibility of \abbr{VSRW} and 
for masking the noise by $\Amg_o^{l,\Tree}$ we needed
non-local functionals, so we in turn approximate 
these in \eqref{eq:fn2_choice} by the local functions 
corresponding to $r,t \in \N$. Indeed, our next
order of business is to use such approximations in relating 
{the relevant functions of $\Fsp_{I_n}^{l,t,n}$ to 
those of} $\Fsp_o^{l,\Tree}$ (when $t,n \to \infty$).}
\begin{lem}\label{lem:loc_fns_conv}
Suppose 
$\Graph_n \stackrel{\text{\abbr{LWC}}}{\Lra} \mu$ for 
some $\mu \in {\cUs}$,
and $\{\nu_{n,+}\}$ converges locally weakly to 
some $\ol{\gm}_+$ supported on $\cI_*$. Then, with 
$\ol{\nu}_+ \in \cP(\ol{\cT}_*)$ denoting the
corresponding limit on average of $\{\nu_{n,+}\}$,
for any fixed $l$ and except for at most countably 
many $\eta >0$,
\begin{align}\label{eq:limit1}
\lim_{t \to \infty} \displaystyle{\lim_{n \ra \infty}}  {\mu_n} [\nu_{n, +}  
\langle \Fsp_{I_n}^{l,t,n} \rangle] &= 
\lim_{t \to \infty} \ol{\nu}_+  [\Fsp_o^{l,t,\Tree}] 
= \ol{\nu}_+[\Fsp_o^{l,\Tree}]\,,\\
\label{eq:limit2}
\lim_{t \to \infty}
\displaystyle{\lim_{n \ra \infty}}  {\mu_n} \Big[ \displaystyle{\sum_{i \in \partial I_n}} \nu_{n,+}( \Fsp_{I_n}^{l,t,n} \ne \Fsp_i^{l,t,n} ) 
 \Big] &= \lim_{t \to \infty} 
 \ol{\nu}_+\Big[ \sum_{i \in \partial o} 
 \bI( \Fsp_{o}^{l,t,\Tree} \ne \Fsp_i^{l,t,\Tree} ) \Big] 
\notag \\
&=
\ol{\nu}_+\Big[ \sum_{i \in \partial o} 
 \bI( \Fsp_{o}^{l,\Tree} \ne \Fsp_i^{l,\Tree} ) \Big] 
\,.
\end{align}
\end{lem}

\noindent
\emph{Proof:} We show that all functions considered 
here can be approximated well by local functions, 
upon which our conclusions follow from
the local weak convergence of $\{\nu_{n,+}\}$. 
Indeed, 
with $|x_j| \le 1$, 
for any graph $\Graph$, {positive} $l,r$, and $i \in V(\Graph)$,
\beq
|y_i^{l,\Graph}(\ul{x}) - y_i^{l,r,\Graph}(\ul{x})| \le 
\frac{1}{l} \int_0^l 
\P_i^{\Graph} 
(\theta_r \le s) ds
\le 
\P_i^{\Graph} 
(\theta_r < l)
=: \ola_i^{l,r,\Graph} \,.
\label{eq:nonlocaltolocal}
\eeq
\blue{In particular, for all $t \in \N$,
\beq\label{eq:nonlocaltolocalA}
\big|
\rmg_i^{l,t,\Graph} - \rmg_i^{l,t,r,\Graph} 
\big| \le \ola_i^{l,r,\Graph} 
\eeq
}
and by \eqref{eq:nonlocaltolocal}--\eqref{eq:nonlocaltolocalA},
for any graph $\Graph$, all $i \in V(\Graph)$, {positive}
$l$, $t$, $r$ and $\eta {> \vep_r \ge 0}$,
\beq \label{eq:ineq1}
\Fsp_i^{l,t,r,\Graph} (\eta+\vep_r) - 
\bI(\ola_i^{l,r,\Graph} \ge \vep_r) \le \Fsp_i^{l,t,\Graph}(\eta) \le 
\Fsp_i^{l,t,r,\Graph}(\eta-\vep_r) + 
\bI(\ola_i^{l,r,\Graph} \ge \vep_r)\,.
\eeq
Further, if the balls $\Ball_i(t \vee r + 1)$ of $\Graph_1$ 
and $\Graph_2$ are isomorphic then 
$a_{i,j}^{l,r,\Graph_1}=a_{i,j}^{l,r,\Graph_2}$ and, restricted to $\Ball_i(t)$, 
the Ising measures $\nu_{+,\Ball_i(t)}$ coincide for both graphs. Consequently,
\beq \label{eq:ineq2}
\Fsp_i^{l,t,r,\Graph_1} (\ul{x},\eta) = \Fsp_i^{l,t,r,\Graph_2} (\ul{x},\eta)\,, \qquad 
\ola_i^{l,r,\Graph_1} = \ola_i^{l,r,\Graph_2}\,. 
\eeq 
{
Choosing $\vep_r^2 = \mu [\ola_o^{l,r,\Tree}]$ 
we get from Markov's inequality that 
\beq
\label{eq:bd-vepr}
\mu(\ola_o^{l,r,\Tree} \ge \vep_r) \le \vep_r \,.
\eeq 
{Recall that for $\ze_n = \nu_{n,+}$, {as in Definition \ref{dfn:conv_average}}} we assumed that 
${\oPb}_{n}^s \Ra \ol{\nu}_+^s$ for any fixed $s > t \vee r$, 
hence by \eqref{eq:ineq1}--\eqref{eq:bd-vepr}, {for $\eta > 2 \vep_r$}
\begin{align}
\notag
\ol{\nu}_+\big[ \Fsp_o^{l,t,\Tree}(\eta+2\vep_r)\big] 
-2\vep_r \le
\ol{\nu}_+^{{s}}\big[ \Fsp_o^{l,t,r,\Tree}(\eta+\vep_r)\big]
&-\mu (\ola_o^{l,r,\Tree} \ge \vep_r)\\
\le \liminf_{n \ra \infty} {\mu_n} \Big[
\nu_{n,+} \big\langle \Fsp_{I_n}^{l,t,n}(\eta) \big\rangle \Big] \label{eq:ineq_local_1}
&\le \limsup_{n \ra \infty} {\mu_n} \Big[ \nu_{n,+} \big\langle 
\Fsp_{I_n}^{l,t,n}(\eta) \big\rangle \Big] 
\\
\le 
\ol{\nu}_+^{{s}} \big[ \Fsp_o^{l,t,r,\Tree}(\eta -\vep_r)\big]
&+
\mu(\ola_o^{l,r,\Tree} \ge \vep_r)
\le
 \ol{\nu}_+\big[ \Fsp_o^{l,t,\Tree}(\eta -2\vep_r)\big]
 +2\vep_r \,. \notag
\end{align}
Proceeding to show that $\vep_r \downarrow 0$, recall that 
$\theta_r \ge \tau_r$, the time of the $r$-th jump 
made by the \abbr{VSRW} $\{X_t\}$ on $\Tree$. 
With $\odeg < \infty$, we have by \cite[Corollary 4.4]{aldous_lyons} 
that this continuous time Markov chain is a.s. non-explosive. That is,
$\tau_r \uparrow \infty$ a.s. and hence for $r \to \infty$,
\[
\vep_r^2 = \mu[
\P_o^{\Tree}
(\theta_r < l)] \le 
\mu[\P_o^{\Tree}
(\tau_r \le l)] \ra 0\,.
 \]
} 
Taking $r \to \infty$ and excluding for $\eta>0$ 
the union over $t \in \N \cup \{\infty\}$ of
the countably many points of discontinuity for 
the $[0,1]$-valued, non-increasing, left-continuous 
$\ol{\nu}_+[\Fsp_o^{l,t,\Tree}(\eta)]$,
we deduce that both lower and upper bounds 
in (\ref{eq:ineq_local_1}) converge
to $\ol{\nu}_+[\Fsp_o^{l,t,\Tree}(\eta)]$,
thus establishing the left identity of (\ref{eq:limit1}),
as well as the bounds
\beq\label{eq:bd-rhs-limit1}
\ol{\nu}_+\big[ \Fsp_o^{l,t,r,\Tree}(\eta+\vep_r)\big]
-\vep_r \le \ol{\nu}_+\big[ \Fsp_o^{l,t,\Tree}(\eta)\big] \le
\ol{\nu}_+\big[ \Fsp_o^{l,t,r,\Tree}(\eta-\vep_r)\big]
+\vep_r\,,
\eeq
for all $t,r \in \N$. Further, 
$\Amg_i^{l,t,r,\Tree}(\eta) = \bI\big\{ \rmg_o^{l,t,r,\Tree} \ge 2\eta \big\}$
and for any fixed $i \in V(\Tree)$ the Ising measures
$\nu^{\be,0,t}_{+,(\Tree,i)}$ of Definition \ref{dfn:ising-pm}  
converge locally to $\nu_{+,\Tree}$ when $t \to \infty$.
Consequently, upon taking $t \to \infty$ followed 
by $r \to \infty$, and further excluding for $\eta>0$ 
the countable collection of points of discontinuity 
for any of 
$\{\ol{\nu}_+[\Fsp_o^{l,r,\Tree}(\eta \pm \vep_r)], r \in \N\}$, we deduce 
that  
$$
\ol{\nu}_+ \big[\Fsp_o^{l,t,r,\Tree}(\eta \pm \vep_r)\big] \to 
\ol{\nu}_+ \big[ \Fsp_o^{l,\Tree}(\eta)\big]\,, 
$$
which by (\ref{eq:bd-rhs-limit1}) gives the \abbr{RHS} of
(\ref{eq:limit1}).
Turning to prove the left identity in 
(\ref{eq:limit2}), since
\begin{align*}
{\mu_n} \Big[
 \sum_{i \in \partial I_n} \nu_{n,+}
 (\Fsp_{I_n}^{l,t,n} =1, \Fsp_i^{l,t,n}=0 )
 \Big] &={\mu_n} \Big[\sum_{i \in \partial I_n} 
 \nu_{n,+}( \Fsp_{I_n}^{l,t,n} =0, \Fsp_i^{l,t,n}=1) \Big]
\end{align*}
it suffices to prove that
\beq
\lim_{n \ra \infty} {\mu_n} \Big[\sum_{i \in \partial I_n} \nu_{n,+} \langle 
\Fsp_{I_n}^{l,t,n} \Fsp_i^{l,t,n} \rangle \Big]= \ol{\nu}_+ \Big[ \sum_{i \in \partial o} \Fsp_o^{l,t,\Tree} 
\Fsp_i^{l,t,\Tree} \Big] =: \HF^{{t}}_{11}(\eta) , \label{eq:toprove1}
\eeq
and
\beq
\lim_{n \ra \infty}{\mu_n}\Big[\Delta_{I_n} \nu_{n,+}\big\langle \Fsp_{I_n}^{l,t,n} 
\big\rangle \Big]= \ol{\nu}_+\Big[\Delta_o \Fsp_o^{l,t,\Tree} \Big] =: \HF^{{t}}_1(\eta) . \label{eq:toprove2}
\eeq
Both $\HF^{t}_{11}(\eta)$ and $\HF^t_1(\eta)$ are
bounded (by $\ol{\nu}_+ [\Delta_o] = \odeg$),
non-negative, left-continuous, non-increasing 
functions of $\eta$.
Thus, excluding the 
at most countably many points of discontinuity
of $\eta \mapsto (\HF^t_{11}(\eta), \HF^t_1(\eta))$
{over all choices of $t \in \N \cup \{\infty\}$,}
we establish (\ref{eq:toprove1}) and (\ref{eq:toprove2}) 
upon deriving 
inequalities analogous to (\ref{eq:ineq_local_1})
for 
$\sum_{i \in \partial I_n} \Fsp_{I_n}^{l,t,n}(\eta) 
\Fsp_i^{l,t,n} (\eta)$
and $\Delta_{I_n} \Fsp_{I_n}^{l,t,n}(\eta)$, respectively.
These in turn also provide the analogs of 
(\ref{eq:bd-rhs-limit1}) with $\ol{\nu}_+[\Fsp_o^{t,r}(\eta \pm \vep_r)]$ 
replaced by the non-increasing in $\eta$ 
and uniformly bounded $\HF_{11}^{t,r}(\eta \pm 2\vep_r)$, 
$\HF_1^{t,r}(\eta \pm 2\vep_r)$, respectively, out of 
which we get the \abbr{RHS} of (\ref{eq:limit2}) along the same
lines we used for deriving the \abbr{RHS} of (\ref{eq:limit1}).
\qed

\smallskip
\blue{In Lemma \ref{lem:plusminusdiff} we show that for generic $\eta>0$,
as $l \to \infty$ the \abbr{rhs} of \eqref{eq:limit2} goes to zero, whereas 
the \abbr{rhs} of \eqref{eq:limit1} has the limit point 
$$
g(\eta) := \mu[(1-\alpha_{\Tree}) \liminf_{l \to \infty} \Amg_o^{l,\Tree}(\eta)] \,.
$$  
Utilizing the edge-expander property of $\Graph_n$ to control the 
\abbr{lhs} of the corresponding identities allows us to then deduce
that $g(\eta) \to 0$ when $\eta \to 0$ 
out of which we reach the
stated conclusion that $\mu$-a.e. $\alpha_\Tree = 1$. To be able 
to carry this out, we next show that $\Amg_i^{l,\Tree}$ is 
sufficiently regular for $i \in \partial o$, and that
$\liminf_l \{\rmg_o^{l,\Tree}\}$ is uniformly (in $\Tree$), 
bounded away from zero, at least $\mu^e$-a.e. for each $\mu^e$ 
which is an extremal element of $\cUs$. For proving the 
latter result, we recall \cite[Corollary 4.4]{aldous_lyons}
that every $\mu \in \cUs$ is invariant for the $\cT_*$-valued 
Markov process $s \mapsto \wt{X}_s$, where $\wt{X}_s=(\Tree,X_s)$
(for the \abbr{vsrw} $\{X_s\}$ on $\Tree$, starting at $X_0=o$), 
and say that such $\mu$ is \abbr{vsrw}-ergodic if all the (continuous)-shift 
invariant events for $\wt{X}_\cdot$ are $\mu \otimes \P_o^\Tree$-trivial.
}

%
\begin{lem}\label{lem:boundary_ineq}
If $\mu \in \cUs$ is \abbr{vsrw}-\emph{ergodic} then  
\beq\label{eq:ergodic-rho-lim}
\rmg_o^{l,\Tree} \ra \rmg_\mu \;\; \text{ as } \;\; l \ra \infty\,,
\qquad \text{ for 
 $\mu$-a.e. } \Tree \in \cT_*\,.
\eeq
Further, for all $\mu \in \cUs$ and any fixed $\vep >0$,
\beq
{\mu} \Big[  \sum_{i \in \partial o} \bI \Big(|
\rmg_o^{l,\Tree} - \rmg_i^{l,\Tree} | > \vep \Big)\Big] 
\ra 0 \;\;\; \text{ as } \;\;\; l \ra \infty\,. 
\label{eq:boundary_ergodic}
\eeq
\end{lem}

\blue{
\begin{rmk}\label{rem:yuval}
From \cite{aldous_lyons} we easily deduce that all
extremal elements $\mu^e$ of the convex set $\cUs$
are \abbr{vsrw}-ergodic. Indeed, this trivially holds 
if $\mu^e(\Delta_o=0)=1$. Otherwise, by 
extremality $\mu^e(\Delta_o=0)=0$, in which 
case by \cite[Theorem 4.6, Theorem 4.7]{aldous_lyons}
the ``size-biased-root'' $\widehat{\mu}^e$ is 
\abbr{dsrw}-ergodic (i.e. all shift invariant events
are $\widehat{\mu}^e \otimes \widehat{\P}_o^\Tree$-trivial
for the corresponding stationary $\cT_*$-valued 
Markov chain $\{\wt{Y}_\ell\}$ of Remark \ref{rw-ergodicity}).
Now if $\mu^e$ is not \abbr{vsrw}-ergodic then 
the corresponding stationary $\cT_*$-valued Markov chain 
$\{\wt{X}_\ell\}_{\ell \in \N}$ must be non-ergodic, hence
has some $\cS \subseteq \cT_*$ with $\mu^e(\cS) \in (0,1)$, 
as a trap set (namely, starting from $\wt{X}_0 \in \cS$, 
w.p.1. $\{ \wt{X}_\ell \} \subseteq \cS$,
c.f. \cite[Proposition 1.8]{krengel}).
Clearly 
also $\widehat{\mu}^e(\cS) \in (0,1)$,
and by the same reasoning, due to the \abbr{dsrw}-ergodicity 
of $\widehat{\mu}^e$, with positive 
$\widehat{\mu}^e \otimes \widehat{\P}_o^\Tree$-probability 
$\wt{Y}_0 \in \cS$ and the first exit time $\tau$ of $\cS$
by $\{\wt{Y}_\ell\}$ is finite. Recall that the chain 
$\{\wt{Y}_\ell\}$ is embedded at the jump-times of 
$\{\wt{X}_s\}$, so applying the 
strong Markov property of $\{\wt{X}_s\}$ at the stopping time
$\tau$, we have that $\wt{X}_s = \wt{X}_\tau \notin \cS$ 
for all $s \in [\tau,\tau+1]$
with positive $\mu^e \otimes \P^\Tree_o$-probability,
in contradiction to $\cS$ being a trap set for $\{\wt{X}_\ell\}$.  
\end{rmk}
}

\noindent
{\em Proof}: 
By definition of $a_{i,j}^{l,\Tree}$ we have the representation,
\begin{eqnarray}
\rmg_o^{l,\Tree} = \frac{1}{l}\int_0^l \sum_j \nu_{+,\Tree} \langle x_j \rangle 
\P_o^{\Tree} (X_t =j ) dt = \P^{\Tree}_o \Big[\frac{1}{l} \int_0^l 
\nu_{+,\Tree} \langle x_{X_t} \rangle dt \Big]\,. \label{eq:int_exp}
\end{eqnarray}
Further, if $\mu \in \cUs$ is \abbr{vsrw}-ergodic then
$\mu \otimes \P_o^\Tree$-a.s. 
\beq
\frac{1}{l} \int_0^l \nu_{+,\Tree} \langle x_{X_t} \rangle dt \lra \rmg_\mu 
\notag
\eeq
(cf. ~{\cite[Pages 10-11]{krengel}}), which by (\ref{eq:int_exp}) and 
\abbr{DCT} for conditional expectation,
yields the $\mu$-a.e. convergence (\ref{eq:ergodic-rho-lim}).
Turning to \eqref{eq:boundary_ergodic}, we assume w.l.o.g. that 
$\odeg>0$ and 
setting 
$$
\hf^{l,\Tree} := \sum_{i \in \partial o} \bI \Big(|\rmg_o^{l,\Tree} - \rmg_i^{l,\Tree} | > 2\vep \Big)\,,
$$ 
note that by the triangle inequality, 
for any $l \in \N$, $\vep>0$,
\begin{align}
\mu \langle \hf^{l,\Tree} \rangle 
&\le  {\mu} \Big[\Delta_o \bI \Big(|\rmg_o^{l,\Tree} - \rmg_\mu | > \vep \Big)\Big] + {\mu} \Big[ \sum_{i \in \partial o} \bI \Big(|\rmg_i^{l,\Tree} - \rmg_\mu| > \vep \Big)\Big] \notag\\
&= \frac{1}{\odeg}
\Big[ {\widehat{\mu}} \big[ \bI \big(|\rmg_o^{l,\Tree} 
- \rmg_\mu | > \vep \big)\big] + {\widehat{\mu}} \big[ \frac{1}{\Delta_o}\sum\limits_{i \in \partial o} \bI \big(|\rmg_i^{l,\Tree} - \rmg_\mu | > \vep \big)\big] \Big] 
\notag\\
&= \frac{2}{\odeg} {\widehat{\mu}}
\Big(|\rmg_o^{l,\Tree} - \rmg_\mu | > \vep \Big). \label{eq:boundary_is_0_simplify}
\end{align}
For \abbr{vsrw}-ergodic $\mu$ we have, in view of (\ref{eq:ergodic-rho-lim}),
the convergence to zero of the bound (\ref{eq:boundary_is_0_simplify}). Hence, 
$\mu \langle \hf^{l,\Tree} \rangle \to 0$, namely
(\ref{eq:boundary_ergodic}) holds for \abbr{vsrw}-ergodic measures, 
and in particular for all extremal elements of $\cUs$ (by Remark \ref{rem:yuval}).
Recall that any fixed $\mu \in \cUs$ can be written as a 
\emph{Choquet integral} of extremal measures 
\cite[Lemma 6.8]{aldous_lyons}. So, we have a probability
measure $\Theta$ on the collection of extremal measures of $\cUs$
such that $\mu \langle \hf^{l,\Tree} \rangle 
=\int \mu^e \langle \hf^{l,\Tree} \rangle d\Theta(\mu^e)$ for all $l$. 
The non-negative $\hf^{l,\Tree}$ are bounded by $\Delta_o$ 
hence $0 \le \mu^e \langle \hf^{l,\Tree} \rangle \le \odege$
for all $l$. Further, 
$\int \odege d\Theta(\mu^e) = \odeg$
is finite, so by \abbr{DCT} we 
deduce from the fact that 
$\mu^e \langle \hf^{l,\Tree} \rangle \to 0$ for $\Theta$-a.e. $\mu^e$ 
that $\mu \langle \hf^{l,\Tree} \rangle \to 0$. That is,
(\ref{eq:boundary_ergodic}) holds for all 
$\mu \in \cUs$.
\qed

\bigskip
Equipped with Lemma \ref{lem:boundary_ineq} we 
proceed to identify the limit as $l \to \infty$ of the 
relevant functionals from 
Lemma \ref{lem:loc_fns_conv}.
\begin{lem}\label{lem:plusminusdiff}
Suppose probability measure 
$\ol{\nu}_+=\mu 
\otimes \ol{\nu}_{+,\Tree}$, 
with $\Tree$ distributed according to $\mu \in {\cUs}$
and 
$\ol{\nu}_{+,\Tree} = 
\alpha_\Tree \nu_{+,\Tree} + (1-\alpha_\Tree) 
\nu_{-,\Tree}$ for some 
fixed, measurable 
$\alpha :\cT_* \mapsto [0,1]$, 
with $\alpha_\Tree=1$
whenever $\nu_{+,\Tree}=\nu_{-,\Tree}$.
Then, for any $\eta>0$, 
\beq\label{eq:Bo-tree-lim}
\lim_{l \to \infty} {\Big|
\ol{\nu}_+ \big[ \Fsp_o^{l,\Tree}\big] -
\mu \big[(1-\alpha_\Tree) \Amg_o^{l,\Tree}\big] \Big|} = 0 \,.
\eeq
Furthermore, {for Lebesgue a.e. $\eta \in (0,1)$,} 
\beq
\liminf_{l \to \infty}  
{\ol{\nu}_+}
\Big[\sum_{i \in \partial o} \bI\big(\Fsp_o^{l,\Tree} \ne \Fsp_i^{l,\Tree}  \big)\Big] = 0 
\,. 
\label{eq:bdry_is_0}
\eeq
\end{lem}

\begin{rmk} Recall the {\em branching number} of a rooted 
tree $\Tree \in \cT_*$, 
\beq
\br \Tree := \Big\{ \lambda > 0: \inf_\Pi \sum_{j \in \Pi} \lambda^{-|j|}=0\Big\}, \notag
\eeq
where {$\Pi \subseteq V(\Tree)$} is a 
cutset (i.e. a finite set of vertices 
that every infinite path from the root intersects), 
and {$|j|$ denotes the 
distance in $\Tree$ between $j$ and the root.}
Our proof of Lemma \ref{lem:plusminusdiff} relies on connections 
between $\br \Tree$ and recurrence/transience of 
the \abbr{VSRW} or phase transitions for Ising models on 
$\Tree$ (c.f. \cite{lyons_ising, lyons_randomwalk}). 
\end{rmk}

\noindent
\emph{Proof:} For any $\Tree$, $l$, $\eta>0$ and $i \in V(\Tree)$,
\begin{align*}
1- \nu_{-,\Tree} \langle \Jsp_i^{l,\Tree}(\eta) \rangle &= 
\nu_{-,\Tree} (y_i^{l,\Tree} > -\eta) \\
&= \nu_{+,\Tree} (y_i^{l,\Tree} < \eta)
\ge \nu_{+,\Tree} (y_i^{l,\Tree} \le -\eta) 
= \nu_{+,\Tree} \langle \Jsp_i^{l,\Tree}(\eta) \rangle 
\end{align*}
and consequently
\beq\label{eq:Ddef-minus-bigger-plus}
D_i^{l,\Tree} (\eta) := 
\Amg_i^{l,\Tree}(\eta) \nu_{+,\Tree} (y_i^{l,\Tree} < \eta)
= \Amg_i^{l,\Tree}(\eta) \max \big\{ 
\nu_{+,\Tree} \langle \Jsp_i^{l,\Tree} \rangle, 
1- \nu_{-,\Tree} \langle \Jsp_i^{l,\Tree} \rangle \big\}  \,.
\eeq
Next, recall that 
$\ol{\nu}_{+}[\Fsp_i^{l,\Tree}] = \mu [\Amg_i^{l,\Tree} 
\ol{\nu}_{+,\Tree}\langle J_i^{l,\Tree} \rangle]$,
for all $l \in \N$ and $i \in V(\Tree)$. So, with 
$\alpha_\Tree \in [0,1]$ and
$\Jsp_o^l= \Jsp_o^{l,\Tree} \in \{0,1\}$, fixing $\eta>0$ we get 
(\ref{eq:Bo-tree-lim}) \blue{by showing that 
\beq\label{eq:D-lim}
\lim_{l \to \infty} \mu [D_o^{l,\Tree}(\eta)] =  0 \,.
\eeq
}
To this end, with
$\Amg_i^{l,\Tree} = \bI(\rmg_i^{l,\Tree} \ge 2 \eta)$ we get by Markov's inequality,
\begin{align}
\mu [D_o^{l,\Tree}(\eta)] &
\le 
{\mu} \Big[\nu_{+,\Tree} \big(  
y_{o}^{l,\Tree} - \rmg_o^{l,\Tree} < - \eta \big) \Amg_o^{l,\Tree} \Big] 
\notag \\
\label{eq:var_simplify1}
& \le  \eta^{-2} {\mu} \Big[\Var_{\nu_{+,\Tree}}
(y_o^{l,\Tree}) \Amg_{o}^{l,\Tree} \Big] =
\eta^{-2} 
\mu \Big[ \sum_{j} \Cov_{\nu_{+,\Tree}} (x_o,x_j) \sum_i a_{o,i}^{l, \Tree} 
a_{i,j}^{l, \Tree} \Amg_{i}^{l,\Tree} \Big],
\end{align}
with the latter identity obtained by expanding the 
variance of $y_o^{l,\Tree}=\sum_j x_j a_{o,j}^{l,\Tree}$, then using unimodularity of $\mu$ as well as 
$a_{o,i}^{l,\Tree} = a_{i,o}^{l,\Tree}$
(by reversibility of the \abbr{VSRW} 
on $\Tree$). 

\noindent
Fixing $r \in \N$, we partition the sum 
over $j$ in the \abbr{RHS} of (\ref{eq:var_simplify1})
into Term I consisting of sum over all $j \in \Ball_o(r)$, 
and Term II for the sum over $j \notin \Ball_o(r)$. We then control 
Term II by confirming for 
$\gth : =\tanh(\beta) \in (0,1)$
and all $\Tree \in \cT_*$ the uniform correlation decay
\beq
 0 \le \Cov_{\nu_{+,\Tree}} (x_o, x_j)  \le  \gth^{|j|} \,.
\label{eq:cov_decay}
\eeq
Indeed, it follows
from (\ref{eq:cov_decay}), 
by non-negativity of $\{a_{i,j}^{l,\Tree}\}$ and the fact $\sum_{i,j} a_{o,i}^{l,\Tree} a_{i,j}^{l,\Tree}=1$, that 
\begin{align}
\text{Term II}  \le & \sum_{k =r+1}^\infty  \gth^{k} \mu \Big[ \sum_{j \in {\Ball_o(k-1,k)}}\sum_i  a_{o,i}^{l, \Tree} a_{i,j}^{l, \Tree} \Big]  
\label{eq:TermII_bound}
 \le \gth^r \,.
\end{align}
Turning to prove (\ref{eq:cov_decay}), note that for
any tree $\Tree$ the marginal of $\nu_{+,\Tree}$ on
$\ul{x}_{\Tree'}$ with $\Tree'=(v_0,v_1,\ldots,v_k)$ 
a finite path in $\Tree$, is an Ising measure 
on $\Tree'$ or in turn a
Markov chain of state space $\{-1,1\}$
(for finite $\Tree$ this follows by summation 
over all possible values of 
$\ul{x}_{\Tree\setminus \Tree'}$, 
hence holding also for infinite trees due to (\ref{eq:dlrm})). While this 
\blue{tree-indexed} Markov chain is in general  
non-homogeneous, recall \cite[Lemma 4.1]{berger_kenyon_mossel_peres} that 
for any $v \ne w \in V(\Tree')$ and Ising measure 
$\nu$ on finite $\Tree'$ with $\beta \ge 0$
and {\em any} external magnetic field 
parameters, 
the value of 
$$
\Phi[\nu](v,w):= \nu[x_w | x_v =1] - \nu[x_w | x_v=-1]
= \nu[x_v=-1]^{-1} \big(\nu[x_w|x_v=1] - \nu[x_w]\big)
$$ 
is non-negative (by Griffith's inequality at
$0=B_v \le B_v' \uparrow \infty$), and
maximal at the measure $\nu_f$ of zero 
external magnetic fields. Now, since $x_v \in \{-1,1\}$, 
we get that
\beq
\Cov_\nu (x_v, x_w)  = 2 \nu[x_v=1] \nu[x_v =-1]
\Phi[\nu] (v,w) \le \frac{1}{2} \Phi[\nu_f] (v,w) 
= \Cov_{\nu_f} (x_v, x_w) \,.
\label{eq:cov_compare}
\eeq
The {tree-indexed} Markov chain corresponding to $\nu_f$ is homogeneous, 
of zero-mean and non-degenerate transition probabilities 
$\pi(y|x)
=\frac{1}{2}(1+ xy \gth)$
on $\{-1,1\}$, 
from which we get by direct computation that 
$\Cov_{\nu_f} (x_{v_0},x_{v_k})=\gth^k$, 
and (\ref{eq:cov_decay}) follows from (\ref{eq:cov_compare}).

\noindent
As for Term I, recall that if
$\nu_{+,\Tree} \langle x_o \rangle =0$, then   
$\nu_{+,\Tree}=\nu_{-,\Tree}$ and 
$\Amg_{i}^{l,\Tree} \equiv 0$ for all $i \in V(\Tree)$ 
and $l \in \N$. Therefore, 
\beq\label{eq:bd-term1}
0 \le \text{Term I} \le \mu\Big[ \Big(\sum_{j \in \Ball_o(r)}\sum_i a_{o,i}^{l,\Tree}a_{i,j}^{l,\Tree}\Big) \bI\{\nu_{+,\Tree} \langle x_o \rangle >0\}\Big] 
\eeq
(the non-negativity of Term I is due to 
$\Cov_{\nu_{+,\Tree}}(x_o,x_j) \ge 0$, 
per (\ref{eq:cov_decay})). 
It is further known that for Ising model on tree $\Tree$
with zero external magnetic field, one has
$\nu^{\be,0}_{+,\Tree} \langle x_o \rangle >0$
only for $\beta \ge \beta_c$, where 
$[\br \Tree] \tanh(\beta_c)=1$ 
(see \cite[Theorem 1.1]{lyons_ising}). In particular,
we bound 
$\bI \{\nu_{+,\Tree} \langle x_o \rangle >0\}$ in 
(\ref{eq:bd-term1}) by $\bI\{[\br \Tree] >1\}$, 
and note that
\begin{align*}
\sum_{j \in \Ball_o (r)} \sum_i a_{o, i}^{l , \Tree} a_{i , j}^{l, \Tree} & = \sum_{j \in \Ball_o (r)} \frac{1}{l^2}\int_0^l \int_0^l \sum_i \P^\Tree_o (X_t =i) 
\P^\Tree_i (X_s = j) \ dt \ ds \notag \\
& = \frac{1}{l^2} \int_0^l \int_0^l \P_o^\Tree (X_{t+s} \in \Ball_o (r))\ dt \ ds\,. 
\end{align*}
In case $[\br \Tree] >1$, the  
\abbr{DSRW} on $\Tree$ is transient (see 
\cite[Theorem 4.3]{lyons_randomwalk}). Consequently, for
such a tree {also} $\{X_t\}_{t \ge 0}$ is transient and in particular 
$1 \ge \P^\Tree_o(X_t \in \Ball_o(r)) \ra 0$ as $t \ra \infty$ for any fixed $r \in \N$. By bounded convergence it thus follows that Term I goes to zero as $l \to \infty$,
for arbitrarily large (fixed) value of $r \in \N$.
Taking $r \to \infty$ we conclude from 
(\ref{eq:TermII_bound}) and
(\ref{eq:var_simplify1}) that 
$
\mu [D_o^{l,\Tree}]
\to 0$ as $l \to \infty$,
thereby establishing (\ref{eq:Bo-tree-lim}).

\vskip 10pt
\noindent
Moving now to the proof of (\ref{eq:bdry_is_0}),
for $\{0,1\}$-valued random variables 
$\Amg_o=\Amg_o^{l,\Tree}$, $\Amg_i = \Amg_i^{l,\Tree}$,
$\Jsp_o=\Jsp_o^{l,\Tree}$ and $\Jsp_i=\Jsp_i^{l,\Tree}$,
we clearly have
per $\Tree$, $l \in \N$ and $i \in \partial o$, that
\begin{align*}
\Amg_o \Amg_i \nu_{+,\Tree} (\Jsp_o \ne \Jsp_i) 
& \le \Amg_o \nu_{+,\Tree} \langle \Jsp_o \rangle + \Amg_i 
\nu_{+,\Tree} \langle \Jsp_i \rangle 
\\
\Amg_o \Amg_i \nu_{-,\Tree} (\Jsp_o \ne \Jsp_i) 
& \le 
\Amg_o (1- \nu_{-,\Tree} \langle \Jsp_o \rangle ) + \Amg_i 
(1- \nu_{-,\Tree} \langle \Jsp_i \rangle )
\,.
\end{align*}
Consequently, with $\alpha_\Tree \in [0,1]$ and 
each $\Fsp_j = \Jsp_j \Amg_j$, we have per $\Tree$, $l$, $\eta>0$
and $i \in \partial o$ that 
\begin{align*}
\ol{\nu}_{+,\Tree}(\Fsp_o \ne \Fsp_i) 
&\le \bI(\Amg_o \ne \Amg_i) + 
\alpha_\Tree \Amg_o \Amg_i \nu_{+,\Tree} (\Jsp_o \ne \Jsp_i)
+ (1-\alpha_\Tree) \Amg_o \Amg_i \nu_{-,\Tree} (\Jsp_o \ne \Jsp_i) \\
&\le \bI(\Amg_o \ne \Amg_i) + D_o + D_i \,,
\end{align*}
for $D_i={D}_i^{l,\Tree}$ 
of \eqref{eq:Ddef-minus-bigger-plus}. 
Taking the expectation with respect to {$\Tree$
of unimodular law $\mu$ we thus get that, 
\beq\label{eq:bd-B0-ne-Bi}
\ol{\nu}_+ \Big[ \sum_{i \in \partial o}   
\bI \big( \Fsp_o^{l,\Tree} \ne \Fsp_i^{l,\Tree}  \big)\Big]
\le \mu \Big[\sum_{i \in \partial o} 
\bI (\Amg_o^{l,\Tree} \ne \Amg_i^{l,\Tree})\Big] 
+ 2 \mu \Big[\Delta_o {D}_o^{l,\Tree} \Big]
\,.
\eeq
Since ${D}_o^{l,\Tree} \in [0,1]$ and $\odeg$ finite, 
we have from \eqref{eq:D-lim} that $\mu[\Delta_o {D}_o^{l,\Tree}] \to 0$. 
Turning to deal with the other term on 
the \abbr{RHS} of (\ref{eq:bd-B0-ne-Bi}),
note that for any $\eta,\vep>0$, if
$\Amg_o^{l,\Tree}(\eta) \ne \Amg_i^{l,\Tree}(\eta)$,
then either 
$|{\rmg}_o^{l,\Tree} - {\rmg}_i^{l,\Tree}| > \vep$
or $\rmg_o^{l,\Tree} \in [2\eta-\vep, 2\eta+\vep)$.
{Further, with $\odeg$ finite,
integrating the non-negative
\beq
\gE(\eta):= \liminf_{\vep \ra 0} \liminf_{l \ra \infty}
\, {\mu} \Big[ \Delta_o \bI \Big(
\rmg_o^{l,\Tree} \in [2\eta-\vep,2\eta+\vep) \Big) \Big]
\notag
\eeq
over $\eta$, we get by Fatou's lemma and Fubini's theorem that
$$
\int_0^1 \gE(\eta) \mathrm{d}\eta \le {
\liminf_{\vep \ra 0} \liminf_{l \ra \infty}
\, {\mu} \Big[ \Delta_o 
\int_0^1 
\bI \Big(
\rmg_o^{l,\Tree} \in [2\eta-\vep,2\eta+\vep) 
\Big) \mathrm{d}\eta \Big] = 0 \,.}
$$
Consequently, $\gE(\eta)=0$ for a.e. $\eta \in (0,1)$, in which case 
the identity \eqref{eq:boundary_ergodic} of Lemma \ref{lem:boundary_ineq} 
completes} the proof of (\ref{eq:bdry_is_0}).
\qed

\vskip 10pt
\noindent
\emph{Proof of Theorem} \ref{thm:average_+}.
Recall part (a) of Lemma \ref{lem:olnu+} that 
any sub-sequential local weak limit point $\ol{\gm}_+$ of $\{\nu_{n,+}\}$, is effectively a distribution over random 
$\alpha : \cT_* \mapsto [0,1]$. From Remark \ref{rmk:conv-av} we 
know that to such $\ol{\gm}_+$ corresponds
$\ol{\nu}_+ = \mu \otimes \ol{\nu}_{+,\Tree}$ with
$\ol{\nu}_{+,\Tree} = \alpha_\Tree \nu_{+,\Tree}
+ (1-\alpha_\Tree) \nu_{-,\Tree}$ for some
fixed measurable $\alpha : \cT_* \to [0,1]$, 
where without loss of generality $\alpha_\Tree=1$ 
whenever $\nu_{+,\Tree} = \nu_{-,\Tree}$ (i.e. 
$\nu_{+,\Tree} \langle x_o \rangle = 0$), as done
in Lemma \ref{lem:plusminusdiff}. 
In particular, it suffices to show that 
the assumed edge-expansion 
property of $\{\Graph_n\}_{n \in \N}$ yields 
\beq
{\mu} \big[ 
(1-\alpha_\Tree) 
\bI\{\nu_{+,\Tree} \ne \nu_{-,\Tree} \}
\big] =0, \label{eq:conclude_2}
\eeq
for then also $\ol{\gm}_+$-a.e. $\alpha_\Tree=1$, as claimed.
To this end, recall
{Lemma \ref{lem:boundary_ineq} (and Remark \ref{rem:yuval}),} that 
for any extremal element $\mu^e$ of $\cUs$ and for 
$\mu^e$-a.e. $\Tree$,
$$
{\mu^e} [ \nu_{+,\Tree} \langle x_o \rangle] =
\lim_{l \to \infty} \rmg_o^{l,\Tree}  \,.
$$
In particular, setting 
$$
\sS_\pm := \{ \Tree: \nu_{+,\Tree} \ne \nu_{-,\Tree}, \; \liminf_{l \ra \infty} 
\rmg_o^{l,\Tree} = 0 \}\,,
$$
we have that $\mu(\sS_{\pm})=0$ for each 
extremal $\mu \in \cUs$ and thus for all $\mu \in \cUs$. 
{Consequently, \eqref{eq:conclude_2} holds as soon as 
\beq\label{eq:conclude3}
{\mu} \Big[
(1- \alpha_\Tree) 
\bI \Big\{ \liminf_{l \ra \infty} 
\rmg_o^{l,\Tree} >0\Big\}  \Big] = 0 \,.
\eeq
Now for any $l$, {$t$}, $n$, $\eta$ and $\ul{x}$, let
$$
\olS_n := n^{-1} \sum_{i=1}^n \Fsp_i^{l,t,n}\,.
$$
That is, $\olS_n = n^{-1} |W^{l,t,n}|$ for the subset of vertices} 
\beq
W^{l,t,n} (\ul{x},\eta) :=
\{ i \in V_n : \Fsp_i^{l,t,n}(\ul{x},\eta) =1 \} \,.
\notag
\eeq
{Setting $\delta:=\eta/2$ for $\eta>0$ such that  
both Lemma \ref{lem:loc_fns_conv} and Lemma \ref{lem:plusminusdiff} hold,}
recall Lemma \ref{lem:u_n_small} that 
whenever $\sum_j x_j \ge 0$
$$
1-\olS_n \ge 1 - {\mu_n} [\Jsp_{I_n}^{l,n}] \ge 
\delta \,.
$$ 
Further, since 
$\{\Graph_n\}_{n \in \N}$ are  $(\delta,1/2, \lambda_\delta)$ 
edge-expanders, we have for such $\ul{x}$ that  
\begin{align*}
\{\olS_n \ge \delta\} \;\; \Longrightarrow \;\;
\frac{1}{n} \sum_{(i,j) \in E_n} \bI (\Fsp_i^{l,t,n} \ne \Fsp_j^{l,t,n})& \ge  
\lambda_\delta \min \big\{ \olS_n, 1 - \olS_n \big\} 
\ge \delta \lambda_\delta \,. \notag
\end{align*}
Taking the expectation with respect to $\nu_{n,+}$ we find that 
\beq
\f{1}{2} {\mu_n} \Big[ 
\sum_{i \in \partial I_n} \nu_{n,+}( \Fsp_{I_n}^{l,t,n} \ne \Fsp_i^{l,t,n} ) \Big] \ge \delta \lambda_\delta
\nu_{n, +} \big( \olS_n \ge \delta \big)
\ge 
\delta \lambda_\delta
\Big( {\mu_n} \big[ \nu_{n,+} \langle \Fsp_{I_n}^{l,t,n} \rangle \big]- \delta  \Big), \notag
\eeq
since 
{$\P (\olS \ge \delta) \ge \E [\olS] - \delta$
for any random variable $\olS \le 1$} and $\delta >0$.
Considering {first the limit over the 
sub-sequence $n_\ell$ such that $\nu_{n_\ell,+}$
converges locally weakly to $\ol{\gm}_+$, followed by
the limit $t \to \infty$}, we deduce from 
Lemma \ref{lem:loc_fns_conv} that 
\beq
{\ol{\nu}_+} \Big[ 
\sum_{i \in \partial o}
\bI ( \Fsp_{o}^{l,\Tree} \ne \Fsp_i^{l,\Tree}) \Big] \ge 2 \delta \lambda_\delta \Big(\ol{\nu}_+ \langle \Fsp_{o}^{l,\Tree} \rangle - \delta  \Big). \notag
\eeq
Hence, considering $l \ra \infty$, 
by Lemma \ref{lem:plusminusdiff} and Fatou's lemma we get that,
\begin{eqnarray}
\delta \ge {
\mu \big[(1-\alpha_\Tree)
\liminf_{l \ra \infty} \Amg_o^{l,\Tree}(\eta) \big]} 
\ge
\mu \Big[ (1-\alpha_\Tree) 
\bI  \big\{\liminf_{l \to \infty} \rmg_o^{l,\Tree} > {2} \eta \big\} \Big] \,.
\label{eq:conclude5}
\end{eqnarray}
Taking now $\eta \to 0$ along suitable sub-sequence, we arrive at
\eqref{eq:conclude3} and complete the proof. 
\qed

\section{Continuity of $\U( \cdot, 0)$ in $\be$ and edge-expander property}\label{section:U_beta_continuous}

{With continuity of $\be \mapsto \U(\be,0)$ 
at $\be < \be_c$ being a consequence of 
uniqueness of the corresponding Ising Gibbs 
measure on $\Tree$, we prove here such continuity 
for any $\mu \in \cU$ supported on trees 
of minimum degree at least three and  
all $\be > \be_\star$, and also at
$\be=\be_c$ for all \abbr{UMGW} measures,
concluding the section with the proof of edge-expander property of the corresponding configuration models.}  
\begin{lem} \label{lem:U_cont}
Suppose $\mu \in \cUs$ such that
$\mu$-a.e. the tree
$\Tree$ has minimum degree {at least} 
$d_\star > 2$ and set $\be_\star:=\atanh[(d_\star -1)^{-1}]$. Then, $\be \mapsto \U (\be, 0)$ is 
continuous on $(\be_\star, \infty)$.
\end{lem}

\noindent
In the next lemma we provide sufficient condition for 
continuity of $\U(\be,0)$ at $\be=\be_c$, in case
$\be_c(\Tree)=\be_c$ is constant for 
$\mu$-a.e. {infinite} $\Tree$.

{
\begin{lem}\label{lem:U_cont_critical}
Suppose $\mu \in \cUs$ and 
$\be_c(\Tree)=\be_c$ finite, for $\mu$-a.e. infinite $\Tree$. If
\beq\label{eq:sT-t-def}
S_\Tree(t) := \sum_{k=1}^t (\br \Tree)^{2k} |
\partial \Tree(k)|^{-2}
\eeq
diverges 
for $\mu$-a.e. infinite $\Tree$, then 
$\be \mapsto \U(\be,0)$ is continuous at $\be=\be_c$.
\end{lem}
\begin{rmk}\label{rmk:T-t}
 Same applies if 
$|\partial \Tree(k)|$
 in (\ref{eq:sT-t-def})
taken for
size of
subset of $\partial \Tree(k)$
connected 
to $\partial \Tree(t)$. 
\end{rmk}
}

We defer the proof of these two lemmas to 
the sequel, proving first  
Lemma \ref{lem:assumption_gw_mgw} 
by verifying that
\abbr{UMGW} measures
satisfy the assumptions of Lemma \ref{lem:U_cont_critical}.

\vskip 10pt
\noindent
{\em Proof of Lemma} \ref{lem:assumption_gw_mgw}:  
Since on any finite tree $\Tree$ there is only one Ising Gibbs measure, $\be\mapsto \U(\be,0)$ is continuous 
for unimodular measures supported on finite trees. It 
thus suffices to prove the continuity of $\U(\cdot,0)$ 
for super-critical \abbr{UMGW} measures conditioned 
on non-extinction. Hence we merely need to verify the assumptions of Lemma \ref{lem:U_cont_critical} for such
\abbr{UMGW} measures conditioned on non-extinction. 
To this end, assume first that all entries of the mean 
matrix $\wh{\AM}$ of Definition \ref{dfn:umgw} are finite.

\noindent
$\bullet$ {\em Branching number:} We need to show that, 
for super-critical \abbr{UMGW}
conditioned on non-extinction, $\be_c(\Tree)=\be_c$
for almost every $\Tree$.  
By the one to one relation between $\br \Tree$ and $\be_c(\Tree)$ (c.f. \cite[Theorem 1.1]{lyons_ising}), 
it suffices to show that conditioned on non-extinction, 
$\br \Tree$ is constant \abbr{UMGW}-a.e.
{This follows from having 
$\br \Tree_{v \ra o}$ 
constant, conditional on non-extinction of 
$\Tree_{v \ra o}$, for \abbr{UMGW} almost 
every $\Tree$ and $v \in \partial o$
(since $\br \Tree = \max_{v \in \partial o} 
\{\br \Tree_{v \ra o}\}$, with zero branching number for 
finite trees and the non-extinction of $\Tree$
equivalent to non-extinction of some $\Tree_{v \ra o}$).
Each $\Tree_{v \ra o}$ has the same super-critical \abbr{MGW} law corresponding to probability kernels $\wh{P}_{i,j}$ over the extended type space $\cQ_\AM$, 
so our claim follows from \cite[Proposition 6.5]{lyons_randomwalk} 
which says that for any super-critical, positive regular,
non-singular \abbr{MGW} law of finite mean matrix $M$, regardless of the type of its root-vertex, conditional 
on its non-extinction the branching number of such 
\abbr{MGW} tree is a.s. the spectral radius $r(M)$
of $M$.  
}

\noindent
{
$\bullet \ S_\Tree$ {\em diverges a.s.:} 
Having finite, positive
regular and non-singular mean matrix $\wh{\AM}$, 
recall the Kesten-Stigum 
characterization of the a.s. finite limit 
of $r(\wh{\AM})^{-k} |\partial \Tree_{v \ra o} (k)|$
conditional on non-extinction of 
$\Tree_{v \ra o}$ (generated according 
to the \abbr{MGW} law with 
probability kernels $\wh{P}_{i,j}$ and 
type space $\cQ_\AM$, for example, see
\cite[Theorem 1]{kurtz_lyons_pemantle_peres} 
).
With $\Delta_o$ finite a.s., by the preceding argument 
it follows that $S_{\Tree}(t) \to \infty$ a.s.
conditional on non-extinction of the \abbr{UMGW} tree.}

{
Turning to the case where some entry 
of $\wh{\AM}$ is infinite, consider the
following truncation of $\wh{P}_{i,j}$,  
$$
\wh{P}_{i,j}^\ell(\ul{k})
:= \wh{P}_{i,j}(\ul{k}) \bI_{\{\|\ul{k}\| \le \ell\}}
+ \bI_{\ul{k}=\ul{0}} \sum_{\| \ul{k'} \|>\ell} 
\wh{P}_{i,j}(\ul{k'}) \,. 
$$ 
For all $\ell$ large enough, both 
positive regularity and non-singularity of $\wh{\AM}$ are inherited by the 
finite mean matrices $\wh{\AM}^\ell$ corresponding to the kernels $\wh{P}^\ell$. 
Further, positive regularity of 
the matrix $\wh{\AM}$ having some infinite entries 
implies that $r(\wh{\AM}^\ell) \to \infty$
as $\ell \to \infty$. Hence, by the preceding proof, 
upon choosing $\ell$ large enough, one can make 
$\br \Tree_{v \ra o}$ under the 
kernels $\wh{P}_{i,j}^\ell$ uniformly arbitrarily large,
conditioned on non-extinction of $\Tree_{v \ra o}$. 
Since $\br \Tree_{v \ra o}$ under 
kernels $\wh{P}_{i,j}^\ell$ is stochastically dominated 
by that for kernels $\wh{P}_{i,j}$, it follows that 
conditioned on non-extinction of $\Tree_{v \ra o}$, 
almost surely $\br \Tree_{v \ra o}=\infty$.
Therefore, a.s. $\br \Tree = \infty$ conditional 
on non-extinction, and all assumptions of Lemma 
\ref{lem:U_cont_critical} are satisfied.
\qed

\vskip 5pt
To prove Lemma \ref{lem:U_cont} we identify functions $\U_{\ell}(\be) \le \U(\be,0)$ that are
non-decreasing in $\ell \in \N$ and $\be \ge 0$, so 
the left continuity of 
$\U(\be,0)$ follows by interchanging 
the order of limits in $\be$ and $\ell$, provided that
\beq\label{eq:u-ell-conv}
\U(\be,0) = \lim_{\ell \to \infty} \U_\ell(\be) \,.
\eeq 
Indeed, for $\Tree \in \cT_*$, non-negative
$\be$, $\ell$ and $\{H_v, v \in V(\Tree)\}$, 
consider the Ising model 
$\nu_{\Tree(\ell)}^{\be,\{H_v\}}$
of (\ref{eq:model_vertex_dependent}),
for graph $\Tree(\ell)$, inverse temperature
parameter $\be$ and external field 
$B_v=H_v \bI_{v \in \partial \Tree(\ell)}$,
with 
$m_{\ell}(\{H_v\}) 
=\nu_{\Tree(\ell)}^{\be,\{H_v\}}\langle x_o \rangle$
denoting its root magnetization. 
Key to the proof of (\ref{eq:u-ell-conv}) is 
the joint continuity property (\ref{eq:joint-cont}) of 
$(\be,\ell) \mapsto m_{\ell}(\{h^{\be'}_v\})$, 
where
$$
h^{\be'}_v := \atanh\big( \nu^{\be',0}_{+,\Tree_{v \to o}} 
\langle x_v \rangle \big) \,, \quad v \in V(\Tree) 
$$
and $\Tree_{v \ra o}$ denotes the connected 
component of the sub-tree of $\Tree$ rooted at $v$, 
after the path between $v$ and $o$ has been deleted
(so $\Tree_{o \ra o}=\Tree$).
\begin{lem}\label{lem:h_cont}
If $\be > \be_0$ such that 
$(d_\star-1) \tanh(\be_0) > 1$, 
then there exists ${\kappa=\kappa(\beta,\beta_0,d_\star)}$ 
finite such that for any $\Tree \in \cT_*$ 
of minimum degree at least $d_\star > 2$
and all $\ell \ge 1$,
\beq\label{eq:joint-cont}
0 \le \ell \Big[ 
m_{\ell} (\{h^{\be}_{v} \}) 
- m_{\ell}(\{h^{\be_0}_{v} \}) 
\Big] \le {\kappa} \,.
\eeq
\end{lem}

\noindent
{\em Proof:} Fixing $\be>\be_0>0$, let
$\gth:=\tanh(\beta)$,
$\gth_0:=\tanh(\beta_0)$. 
Using
$v \hookrightarrow w$ 
to denote that $v$ is the parent 
of $w$ in $\Tree \in \cT_*$,
the identity 
(\ref{eq:recursion_m}) becomes
\beq
h^{\be}_{v} = \sum_{\{w : v \hookrightarrow w\}}  
f_{\gth}(h^{\be}_{w}) \,, 
\label{eq:message_recursion}
\eeq
for $f_\gth (h):=\atanh(\gth \tanh(h))$. 
Since $g: [0,1] \to (1,\infty)$ given by 
$$
g(0)=\frac{\gth}{\gth_0}\,, \quad
g(r) = 
\frac{\atanh(\gth r)}{\atanh(\gth_0 r)}\,,
\quad \forall r \in (0,1] \,,
$$
is continuous, necessarily  
$g(r) \ge 1+ \vep$ for some 
$\vep=\vep(\be,\be_0) >0$ and all $r \in [0,1]$. 
Hence, by Proposition \ref{prop:marginal_ising},
Griffith's inequality and our uniform lower 
bound on $g(\cdot)$, for any $k \ge 0$ we have 
\begin{align}
m_{k+1}(\{h^{\be_0}_{w}\}) 
& = m_{k} 
\Big(\big\{ \sum_{\{w : v \hookrightarrow w\}}
f_\gth (h^{\be_0}_{w}) \big\}\Big) 
= 
 m_{k} 
\Big(\big\{ \sum_{\{w: v \hookrightarrow w\}}
g(\tanh(h^{\be_0}_{w}))
f_{\gth_0} (h^{\be_0}_{w}) \big\}
\Big)
\nonumber \\
&\ge 
m_{k} 
\Big(\big\{ \sum_{\{w: v \hookrightarrow w\}} (1+ \vep) 
f_{\gth_0} (h_{w}^{\be_0}) \big\}\Big)  =  
m_{k} (\{(1+ \vep) h_{v}^{\be_0}\}) \,, \label{eq:ineq+f3}
\end{align} 
with the last equality due to (\ref{eq:message_recursion}).
The minimum degree of $\Tree$ is at least $d_\star$, 
so we have by Griffith's inequality that 
$h_{w}^{\be_0} \ge h^{\be_0}_\star$ for
all $w \in V(\Tree)$ and 
$h^{\be_0}_\star :=\atanh(r^{\be_0}_\star)$
with $r_\star^{\be_0}$ 
the positive root magnetization 
for Ising plus measure on the $(d_\star-1)$-ary tree, at parameter $\be_0$
(which by assumption exceeds the critical parameter 
for Ising measure on the regular tree $\Tree_{d_\star}$).
It then follows 
from (\ref{eq:message_recursion}) that moreover 
$h_{v}^{\be_0} \ge \xi {\Delta}_v$,
with $\xi := \frac{1}{2} f_{\gth_0}(h^{\be_0}_\star)$ 
strictly positive. 
Using (\ref{eq:message_recursion}) once more, we see 
that $h_{v}^{\be} \le f_\gth (1) {\Delta}_v = \be \Delta_v$ for all $v \in V(\Tree)$. Thus, by Griffith's inequality,  
\beq \label{eq:ineq+f1}
m_{k+1} (\{h_{w}^{\be}\}) 
=m_{k}(\{h_{v}^\be\}) 
\le m_{k} (\{\be {\Delta}_v\}) \le 
m_{k} (\{(\be/\xi) h_{v}^{\be_0}\}).
\eeq 

\noindent
Choosing $\vep>0$ small enough, we have
$\beta/\xi = 1 + \kappa \vep$ with $\kappa > 1$ finite,
hence by the concavity on $\R_+$ of 
$\lambda \mapsto m_{k} (\{\lambda H_v\})$, 
for each $k \ge 0$ and non-negative $\{H_v\}$
(which is a special case of the GHS inequality, 
see \cite{griffiths_hurst_sherman}),
%
we get the inequality,
\beq
m_{k} (\{(\beta/\xi) h_{v}^{\be_0}\}) - m_{k} (\{h_{v}^{\be_0}\}) \le \kappa 
\Big[ m_{k} (\{(1+\vep) h_{v}^{\be_0}\}) - m_{k} (\{h_{v}^{\be_0}\}) \Big]\,.\label{eq:ineq+f2}
\eeq
Combining (\ref{eq:ineq+f3}), (\ref{eq:ineq+f1})
and (\ref{eq:ineq+f2}) we deduce that
\beq
m_{k+1} (\{h^{\be}_{w}\}) - 
m_{k+1} (\{h^{\be_0}_{w}\}) 
\le \kappa [m_{k+1}(\{h^{\be_0}_{w}\}) 
- m_{k}(\{h^{\be_0}_{w}\})]\,.
\nonumber
\eeq
Recall, for example from (\ref{eq:ineq+f3}), that 
$k \mapsto m_{k} (\{h^{\be_0}_{v}\}) \in [0,1]$ 
is non-decreasing, and bounded above by
$m_{k}(\{h_v^\be\})$ which is   
independent of $k$. Hence, summing the latter 
inequality over $k=0,\ldots,\ell-1$ results with
\beq
0 \le \ell \Big[ 
m_{\ell} (\{h^{\be}_{v} \}) 
- m_{\ell}(\{h^{\be_0}_{v} \}) 
\Big] 
\le
\sum_{k=1}^{\ell}
\Big[
m_{k} (\{h^{\be}_{w} \}) 
- m_{k}(\{h^{\be_0}_{w} \})
\Big]
\le \kappa m_\ell(\{h_w^{\be_0}\})
\le \kappa \,,
\nonumber
\eeq 
as claimed. \qed

\vskip 5pt
\begin{rmk}\label{rmk:h-cont-ext}
Fixing $i \in \partial o$ and keeping same 
choices of external field, the argument we
used in proving Lemma \ref{lem:h_cont} also
establishes (\ref{eq:joint-cont}) when
$m_\ell(\cdot)$ is replaced by the Ising 
root magnetization on 
$\Tree(\ell) \cap \Tree_{o \to i}$, 
as well as when it is replaced 
by the magnetization at $i$ 
for such Ising models on 
$\Tree(\ell+1) \cap \Tree_{i \to o}$.
Hereafter, we denote the former by 
$m_{\ell,o \to i}(\cdot)$ and the 
latter by  $m_{\ell+1,i \to o}(\cdot)$.
\end{rmk}

\begin{rmk} 
For \abbr{UGW} measure $\mu$ the variables
$\{h^{\be'}_v, v \ne o \}$ are identically distributed,
each having the law we called $h^{\be',+}$ in Lemma \ref{lem:recursion}. Starting the 
recursion (\ref{eq:recursion_h}) with 
$h^{(0)} \stackrel{d}{=} h^{\be_0,+}$ yields 
the sequence $h^{(\ell)}$ having the laws of 
$\atanh\big(m_{\ell+1,i \ra o}(\{h^{\be_0}_v\})\big)$.
We have just coupled these with 
$\atanh\big(m_{\ell+1,i \ra o}(\{h^{\be}_v\})\big)$
whose law equals $h^{\be,+}$, establishing the 
convergence in law of Lemma \ref{lem:recursion}
(and by Griffith's inequality this 
extends to starting laws which stochastically 
dominate $h^{\be_0,+}$).
\end{rmk}

\noindent
{\em Proof of Lemma} \ref{lem:U_cont}: 
As mentioned in Remark \ref{rmk:left-cont},
fixing $\be > \be_0 > \be_\star$
it suffices to show that $\U(\be,0)$ is 
left continuous at $\be$. To this end, for 
any infinite $\Tree$ and integer $\ell \ge 1$, 
using the Ising model $\nu_{\Tree(\ell)}^{\be, \{h_{w}^{\be_0}\}}$ 
on $\Tree(\ell)$ with positive external field 
only at $\partial \Tree(\ell)$, as in Lemma \ref{lem:h_cont}, we define
\beq\label{eq:U-ell-exp}
\U_\ell (\be) = \f{1}{2}
\E_\mu \Big[ \sum_{i \in \partial o} 
\nu_{\Tree(\ell)}^{\be, \{h_{v}^{\be_0}\}}
\langle x_o  x_i \rangle \Big].
\eeq
With $\Tree(\ell)$ a finite graph, fixing 
$\be_0$ and $\ell$, the function 
$\be \mapsto \U_\ell(\be)$ is continuous and 
non-decreasing (by Griffith's inequality). 
By Proposition \ref{prop:marginal_ising} we 
further have that
\beq
\U_{\ell+1} (\be) = 
\f{1}{2} \mu \Big[ \sum_{i \in \partial o} 
\nu_{\Tree(\ell)}^{\be,\{H_v\}}
\langle x_o x_i \rangle \Big]\,,
\nonumber
\eeq
and since $\be > \be_0$, it follows from 
(\ref{eq:message_recursion}) and the monotonicity of 
$\gth \mapsto f_\gth (h)$, that for any 
$v \in \partial \Tree(\ell)$,
$$
H_v := \sum_{\{w: v \hookrightarrow w\}} f_{\gth} (h_w^{\be_0})
\ge 
\sum_{\{w: v \hookrightarrow w\}} 
f_{\gth_0} (h_w^{\be_0})
= h_v^{\be_0} \,.
$$
By yet another appeal to Griffith's inequality we deduce
that $\ell \mapsto \U_\ell(\be)$ is also non-decreasing.
Recall that $h_v^\be \ge h_v^{\be_0}$ for all $v \in V(\Tree)$, so by similar reasoning, $\U_\ell(\be) \le 
\U(\be,0)$ and as explained before it remains only to 
establish (\ref{eq:u-ell-conv}). To this end,
in view of (\ref{eq:Ising-corr-formula}), we have that
for any $i \in \partial o$ and $\{H_v, v \in V(\Tree)\}$,
$$
\nu_{\Tree(\ell)}^{\be, \{H_{v}\}}
\langle x_o  x_i \rangle = 
F\Big(\gth, m_{\ell,i \ra o}(\{H_v\}) m_{\ell,o \ra i}(\{H_v\})\Big)  
\,,
$$
where $F(\gth,r)$ of (\ref{eq:edge-corr}) 
is continuous and bounded on $[0,1]^2$. Thus, 
with $\Psi(\gth,\delta):= 
\sup \{ |F(\gth,r) - F(\gth,r')|$
over $r,r' \in [0,1]$ such that $|r-r'| \le \delta \}$
and $\delta_{\ell}:=2\kappa/(\ell-1)$, clearly 
$\Psi(\gth,\delta_\ell) \to 0$ as $\ell \to \infty$.  
Now, in view of Remark \ref{rmk:h-cont-ext},
the expression (\ref{eq:U-ell-exp}) for $\U_{\ell}(\be)$
and the corresponding expression for $\U(\be,0)$, we 
deduce that
$$
|\U(\be,0)-\U_{\ell}(\be)| \le \frac{1}{2} 
\Psi(\gth,\delta_\ell) 
\odeg
\,, 
$$
from which (\ref{eq:u-ell-conv}) follows. 
\qed

\begin{rmk}\label{rmk:U_cont_difficulty}
It is easy to see that the proof of 
Lemma \ref{lem:U_cont} applies at any 
$\be \ge 0$ and $\mu \in \cUs$ such that for some 
$\be_0<\be$ one has a bound of the type 
(\ref{eq:joint-cont}). That is, as soon as 
$m_{\ell} (\{h^{\be}_{v}\})- m_{\ell}(\{h^{\be_0}_{v}\}) 
\to 0$ in probability, when $\ell \to \infty$.  
Further, the proof of (\ref{eq:joint-cont}) is 
completely general, {\em except} for requiring in 
(\ref{eq:ineq+f1}) that $h_v^{\be}/h_v^{\be_0}$ 
(alternatively, $r_v^\be/r_v^{\be_0}$),
be uniformly bounded over $v \in V(\Tree)$.
Unfortunately, while 
$h_v^{\be_0}$ is strictly positive as soon as
$\be_0 > \be_c(\Tree)$, even for \abbr{UGW} $\mu$, 
when $\be \in (\be_c,\be_\star)$ such ratios may 
be arbitrarily large (with small $\mu$-probability,
but nevertheless, they appear at some $v$
and a.e. infinite tree $\Tree$). We did not 
find a way to by-pass this technical difficulty,
hence our requirement of $\be > \be_\star$.
\end{rmk}

\begin{rmk}
Lemma \ref{lem:h_cont} and Lemma \ref{lem:recursion} 
are the analogs of \cite[Lemma 4.3]{dembo_mont}
and \cite[Lemma 2.3]{dembo_mont}, respectively, 
in case of zero external field and low temperature 
(i.e. $\be>\be_\star$). While we do not pursue 
this here, utilizing the former one can establish 
similar conclusions as done in \cite{dembo_mont} 
based on \cite[Lemma 2.3 and Lemma 4.3]{dembo_mont}.
\end{rmk}

\medskip
The proof of Lemma \ref{lem:U_cont_critical}
builds on results from \cite{pemantle_peres},
to which end we introduce few 
relevant definitions and notations. 
First, for any finite $(\Tree,o) \in \cT_*$ let 
$\partial_\star \Tree$ denote the collection 
of {\em rays} emanating from $o$, namely
finite non-backtracking paths in one-to-one 
correspondence with the leaves of $\Tree$ 
other than $o$ (where each such ray terminates).
Next, a {\em flow} $\varpi$ on such $(\Tree,o)$ 
is a non-negative function on $E(\Tree)$, of {\em strength}
$|\varpi| :=\displaystyle{\sum_{y: o \hookrightarrow y}} \varpi (oy)$, such that 
$\varpi(vw) = \displaystyle{\sum_{y: w \hookrightarrow y}} \varpi(wy)$, whenever $v \hookrightarrow w$ 
and $w \notin \partial_\star \Tree$.  
Any given collection of {\em resistances}
$\{R(e) \ge 0: e \in E(\Tree)\}$, 
induces the functional
$$
V_\varpi := \sup \Big\{ \,
\sum_{e \in y} (\varpi(e) R(e))^2
 : y \in \partial_\star \Tree \Big\}\,,
$$
over flows $\varpi$ on $\Tree$, in terms of which 
we define   
\begin{align*}
\capa_3(\Tree) := \sup \{ \, |\varpi|\, : 
\varpi \text{ a flow on } \Tree \text{ with } V_\varpi =1\} \,.
\end{align*}

\noindent
{\em Proof of Lemma \ref{lem:U_cont_critical}}:
For any $(\Tree,o) \in \cT_*$ 
and $e = vw \in E(\Tree)$
let $|e|=|v| \vee |w|$ where $|v|$ denotes the graph distance between $v \in V(\Tree)$ and $o$.
From \cite[Lemma 4.2]{pemantle_peres} we know that 
for any $\gth>0$ there exists $\kappa>0$ such that 
\beq
f_\gth (h) \le \f{\gth h}{(1+(\kappa h)^2)^{1/2}}
\label{eq:f_relation}
\eeq
for $f_\gth (\cdot)$ of (\ref{eq:message_recursion}) 
and all $h \ge 0$. Futher, recall that for 
any finite $t \ge 1$, $\gth =\tanh(\beta) >0$ 
and infinite tree 
$(\Tree,o) \in \cT_*$ without leaves, the positive
$$
h_v^{(t)}(\Tree) := \atanh\big(
\nu^{\be,0,t}_{+,\Tree(t)_{v \ra o}} \langle x_v \rangle )
\,,
$$
satisfies the system of equations (\ref{eq:message_recursion}) at all $|v| < t$, 
starting with $h_w^{(t)}(\Tree) = \infty$ when $|w|=t$
(i.e. $w \in \partial \Tree(t)$).} 
{More generally, in case $(\Tree,o)$ has leaves, 
let $\Tree_t \subseteq \Tree(t)$ denote the union of
all vertices and edges along rays of $\Tree(t)$ 
of length $t$, emanating from $o$. All 
non-root leaves of $\Tree_t$ are at distance $t$ 
from $o$ and it is easy to verify that 
$h_v^{(t)}(\Tree) = h_v^{(t)}(\Tree_t)$ satisfy 
for $v \in \Tree_t$ the corresponding 
equations (\ref{eq:message_recursion}) on $\Tree_t$,
starting with $h_w^{(t)}(\Tree_t)=\infty$ at 
$w \in \partial \Tree(t)$.} 
{
In view of (\ref{eq:f_relation}), it then follows from
\cite[Theorem 3.2]{pemantle_peres} that
\beq\label{eq:bd-capa3-t}
h_o^{(t)} (\Tree) \le \kappa^{-1} \capa_3(\Tree_t) \,,
\eeq
for $\capa_3(\Tree_t)$ corresponding to 
resistances $R(e) = \gth^{-|e|}$ on $(\Tree_t,o)$.
Set $\gth=\tanh(\be)$ for $\be=\be_c(\Tree)$ finite, namely $\gth =1/(\br \Tree)$  (see \cite[Theorem 1.1]{lyons_ising}). If such 
$\capa_3(\Tree_t) \ra 0$ for $t \ra \infty$, 
then by (\ref{eq:bd-capa3-t}) we deduce that 
$$
\nu^{\be,0}_{+,\Tree} \langle x_o \rangle =
\lim_{t \to \infty} 
\tanh\big( h_o^{(t)}(\Tree) \big) \le 0 \,,
$$
so at $\be=\be_c(\Tree)$ there is then a {\em unique}
Ising Gibbs measure on $(\Tree,o)$.
Now, should this happen for $\mu$-a.e.
{infinite} $\Tree$
at the same $\be_c(\Tree)=\be_c$, then 
necessarily $\U(\be_c,0)=0$ and in particular 
$\be \mapsto \U(\be,0)$ is continuous 
at $\be=\be_c$. 
{
With $\Tree_t \subseteq \Tree(t)$, clearly
$$
S_{\Tree_t} := \sum_{k=1}^t \gth^{-2k} 
|\partial \Tree_t (k)|^{-2} \ge S_\Tree (t)
$$ 
of (\ref{eq:sT-t-def}), so it suffices to 
confirm that $\capa_3(\Tree_t) \le S_{\Tree_t}^{-1/2}$
(see also Remark \ref{rmk:T-t}).}
{
To this end, fixing $t \ge 1$ let $\varpi$ be any 
flow on $\Tree_t$ of strength $|\varpi|=1$. 
Then, by the definition of $V_\varpi$,
for any probability measure $\sfp (\cdot)$ 
on $\partial_\star \Tree_t$,
\beq
V_\varpi \ge \sum_{y \in \partial_\star \Tree_t} \Big[ \sum_{e \in y} \varpi^2(e) \gth^{-2|e|} \Big] \sfp (y) = \sum_{k=1}^t \gth^{-2k} \sum_{|e|=k} \varpi^2(e) 
\sum_{y \ni e} \sfp (y). 
\label{eq:capa_1}
\eeq
With slight abuse of notation, set $\sfp (e):= \sum_{y \ni e} \sfp (y)$. 
Note that the thus defined $\{\sfp (e),  
e \in E(\Tree_t)\}$, constitutes a flow of 
strength $|\sfp |=1$. Further, 
$\sum_{|e|=k} \sfp (e)=1$ for any $1 \le k \le t$ 
since all non-root leaves of 
$\Tree_t$ are at $\partial \Tree_t (t)$.
Applying Cauchy-Schwarz inequality and choosing 
$\sfp = \varpi$, we find that 
\beq
\Big[\sum_{|e|=k} \varpi^2(e) \sfp (e)\Big] \ge \Big( \sum_{|e|=k} \varpi(e)
\sfp (e) \Big)^2 = \Big(\sum_{|e|=k} \varpi^2(e) \Big)^2. \label{eq:capa_2}
\eeq
Using Cauchy-Schwarz inequality once more,
\beq
 \Big(\sum_{|e|=k} \varpi^2(e)\Big) \ge \f{1}{|\partial \Tree_t (k)|}\Big(\sum_{|e|=k} \varpi(e) \Big)^2= 
 \f{1}{|\partial \Tree_t (k)|}\,. \label{eq:capa_3}
\eeq
Thus, from (\ref{eq:capa_1}), (\ref{eq:capa_2}) and (\ref{eq:capa_3}), we see that}
{ 
$V_\varpi \ge S_{\Tree_t}$
for any flow $\varpi$ on $\Tree_t$ 
such that $|\varpi|=1$.
%
By simple scaling, it then follows that
$\capa_3(\Tree_t) \le S_{\Tree_t}^{-1/2}$, 
as claimed.}
\qed
}

\vskip10pt
\noindent
{\em Proof of Lemma} \ref{lem:expander_mgw}: 
For each $i \in \cQ$ and $\ul{k} \in \Z_\ge^{|\cQ|}$ let 
$\alpha_{i,\ul{k}}= \ptheta(i)P_i(\ul{k})$, viewed 
as coordinates of the collection 
$$
\bm{\alpha}= (\alpha_{i,\ul{k}})_{i \in \cQ, \ul{k}\in \Z_\ge^{|\cQ|}} 
$$
(which is finite by assumption of bounded support for all $P_i(\cdot), i \in \cQ$).
Fixing $\delta_0 \le 1/2$, 
{for any 
vector $\bm{\delta}=(\delta_{i,\ul{k}})_{i \in \cQ, \ul{k} \in \Z_\ge^{|\cQ|}}$ 
such that $\|\bm{\delta}\| \in (\delta_0,1/2)$, let $W_{\bm{\delta}}$ 
denote a subset of $[n\|\bm{\delta}\|]$ vertices from $V_n$ where
for each $i$ and $\ul{k}$, about 
$n\delta_{i,\ul{k}} (1+o(1))$ of the vertices of 
$W_{\bm{\delta}}$ are of type $i$ and off-springs configuration $\ul{k}$.
For any $n$ and $\vep \ge 0$ denote by $\cSG^{\vep,n}_{\bm{\delta}}$ the event that
within $\Graph_n$ there exists some $W_{\bm{\delta}}$ having precisely 
$[n \vep]$ edges between $W_{\bm{\delta}}$ to $W_{\bm{\delta}}^c$.
\blue{By Definition \ref{dfn:config}, with high probability,
for all large $n$ and each $i$, $\ul{k}$, there are 
$n\alpha_{i, \ul{k}}(1+o(1))$ vertices of type $i\in \cQ$ and 
off-springs configuration $\ul{k}$ in the random graph $\Graph_n$. 
In particular, with high probability only events $\cSG^{\vep,n}_{\bm{\delta}}$ 
having 
\beq\label{eq:del-cond}
\delta_{i,\ul{k}} \le \alpha_{i,\ul{k}}\,, \qquad \qquad \forall i, \ul{k}
\eeq
occur. We have the stated edge-expansion property upon   
the existence of $\vep_0:=\vep_0(\delta_0)>0$ such that the probability 
of the union of all such $\cSG^{\vep,n}_{\bm{\delta}}$ for which 
\eqref{eq:del-cond} holds, $\|\bm{\delta}\| \in (\delta_0,1/2)$ and
$\vep \le \vep_0$, goes to zero as $n \to \infty$. 
Vertex types and edge counts are integer valued, so with both
the length of $\bm{\delta}$ and $\vep \le n^{-1} |E_n|$ 
uniformly bounded, we have at most $n^C$ such events to rule out.}
Consequently, it suffices to show that for any 
$\bm{\delta} \in (\delta_0,1/2)$ satisfying \eqref{eq:del-cond}
and $\vep \le \vep_0$,
}
\beq \label{eq:expander_vep}
\f{1}{n} \log \P(\cSG^{\vep,n}_{\bm{\delta}}) < -\vep < 0,
\eeq
for all large $n$, uniformly over all such 
choices of $\bm{\delta}$ and $\vep$. 
To this end, we first note that for $\vep=0$,
\begin{align}
\f{1}{n} \log \P(\cSG^{0,n}_{\bm{\delta}})  &= \f{1}{n} \log \# \Big\{ \text{choices possible for } W_{\bm{\delta}} \Big\} +\f{1}{n}\log \P \Big\{ \text{such choice matches with itself} \Big\} \notag\\
& =: N_{\bm{\delta}}+Q_{\bm{\delta}}. \notag
\end{align}
We further define  
$\alpha_{i,j}:= \sum_{\ul{k}} k_j \alpha_{i,\ul{k}}$ 
and
$\delta_{i,j}:= \sum_{\ul{k}} k_j \delta_{i,\ul{k}}$
for each $i,j \in \cQ$. Using the approximations,
\beq
 \f{1}{n} \log n ! = \log \big(\f{n}{e}\big)+ o(1) \text{ and } \f{1}{n} \log n !! = \f{1}{2} \log \big( \f{n}{e} \big) + o(1), \notag
\eeq
we have for
$H(q):= - q \log q - (1-q) \log (1-q)$, $q \in [0,1]$, that
\begin{align}
 N_{\bm{\delta}} &\approx \sum_{i, \ul{k}} \alpha_{i, \ul{k}} H\Big( \f{\delta_{i,\ul{k}}}{\alpha_{i,\ul{k}}}\Big)= \sum_{i} \sum_j \sum_{ \ul{k}} \f{k_j{\alpha}_{i,\ul{k}}}{\|\ul{k}\|} H\Big( \f{\delta_{i,\ul{k}}}{\alpha_{i,\ul{k}}}\Big) \label{eq:approx1}
\\
Q_{\bm{\delta}} &\approx -\f{1}{2} \sum_{i \in \cQ} {\alpha}_{i,i} H\Big(\f{{\delta}_{i,i}}{{\alpha}_{i,i}}\Big) - \sum_{i\ne j \in \cQ} {\alpha}_{i,j} H\Big(\f{{\delta}_{i,j}}{{\alpha}_{i,j}}\Big)\,. \label{eq:approx2}
\end{align}
By concavity of $H(\cdot)$, upon noting that 
$\|\ul{k}\| \ge 3$ we have for any $i,j \in \cQ$, that
 \beq
 \sum_{ \ul{k}} \f{k_j{\alpha}_{i,\ul{k}}}{\|\ul{k}\|} H\Big( \f{\delta_{i,\ul{k}}}{\alpha_{i,\ul{k}}}\Big) - \f{1}{2} {\alpha}_{i,j} H\Big(\f{{\delta}_{i,j}}{{\alpha}_{i,j}}\Big) \le \f{1}{3} \sum_{\ul{k}} k_j \alpha_{i, \ul{k}} H\Big( \f{\delta_{i,\ul{k}}}{\alpha_{i,\ul{k}}}\Big) - \f{1}{2} {\alpha}_{i,j} H\Big(\f{{\delta}_{i,j}}{{\alpha}_{i,j}}\Big) \le -\f{1}{6}  {\alpha}_{i,j} H\Big(\f{{\delta}_{i,j}}{{\alpha}_{i,j}}\Big). \notag
 \eeq
With
$\|\bm{\delta}\| \le 1/2 <\| \bm{\alpha}\|=1$ for $\bm{\delta}$ 
{satisfying \eqref{eq:del-cond},} we must have 
${\delta}_{i,j} < {\alpha}_{i,j}$
for at least one pair $(i,j)$. We thus get from 
(\ref{eq:approx1}) and (\ref{eq:approx2}) that
\beq
\limsup_{n \ra \infty}\f{1}{n}\log \P (\cSG_{\bm{\delta}}^{0,n}) \le - \f{1}{6}\sum_{i,j \in \cQ} {\alpha}_{i,j} H \Big(\f{{\delta}_{i,j}}{{\alpha}_{i,j}}\Big) \,,
\label{eq:ubd_entropy}
\eeq
with the \abbr{RHS} strictly negative
(since $H(q)=0$ only for $q \in \{0,1\}$).
Further, the approximations in 
(\ref{eq:approx1}) and (\ref{eq:approx2}) are uniform 
over $\bm{\delta}$, because  
\beq
\sqrt{2 \pi} \le \f{n !}{n^{n+1/2}e^{-n}} \le e, \text{ for all } n. \label{eq:finite_stirling}
\eeq
The supremum of the upper bound of 
(\ref{eq:ubd_entropy}), over the compact set 
of all possible choices of $\bm{\delta}$ is 
strictly negative, yielding (\ref{eq:expander_vep}) for $\vep=0$. Similar rational applies also for all 
$\vep$ small enough. For example, in case $|\cQ|=1$ we have
for $\delta:=\sum_k k \delta_k$ and $\alpha := \sum_k k \alpha_k$,
that $\delta<\alpha$ and 
\beq
\limsup_{n \ra \infty} \f{1}{n} \log \P(\cSG^{\vep,n}_{\bm{\delta}}) \le -\f{{\alpha}}{6} H\Big(\f{{{\delta}}}{{{\alpha}}}\Big) +\f{{\delta}}{2} H\Big(\f{\vep}{{\delta}}\Big) + \f{1}{2} ({\alpha} - {\delta})H\Big(\f{\vep}{{\alpha}-{\delta}}\Big) \le -\f{{\alpha}}{6}H\Big(\f{{{\delta}}}{{{\alpha}}}\Big) +\f{{\alpha}}{2} H\Big(\f{2\vep}{{\alpha}}\Big)\,.
\notag
\eeq
The preceding bound is continuous in $\vep$ and strictly 
negative at $\vep=0$. Consequently, there exists 
$\vep_0>0$ small enough such that this bound is strictly negative at all $\vep \le \vep_0$. Further, from 
(\ref{eq:finite_stirling}) we get uniformity
of the convergence in $n$, over all relevant 
$\bm{\delta}$ and $\vep \le \vep_0$, yielding 
(\ref{eq:expander_vep}) in case $|\cQ|=1$. 
While we do not detail these, the computations 
in case $|\cQ|>1$ and $\vep>0$ are similar.
\qed

\end{document}